\pdfoutput=1
\documentclass[envcountsect, smallextended]{svjour3}
\smartqed  

\usepackage{amssymb}
\usepackage{amsmath}

\usepackage{amsfonts}
\usepackage{xcolor}
\usepackage{enumitem}
\usepackage{graphicx}
\usepackage{lmodern}

\usepackage[caption=false]{subfig}

\usepackage{makecell}
\usepackage{multirow}
\usepackage{algorithm}
\usepackage{algorithmic}

\usepackage{lineno}

\newcommand{\app}{\mathrm{app}}
\newcommand{\dif}{\mathrm{d}}
\newcommand{\Dim}{{ \textup{D}}}
\newcommand{\Int}{\mathrm{int}}
\newcommand{\Tiny}{\textrm{tiny}}

\newcommand{\cond}{\mathrm{cond}}

\spnewtheorem{thm}{Theorem}[section]{\bf}{\it}
\spnewtheorem{defn}[thm]{Definition}{\bf}{\rm}
\spnewtheorem{lem}[thm]{Lemma}{\bf}{\it}
\spnewtheorem{ntn}[thm]{Notation}{\bf}{\it}
\spnewtheorem{rem}[thm]{Remark}{\bf}{\it}
\spnewtheorem{cor}[thm]{Corollary}{\bf}{\it}

\newcommand{\algorithmicprecondition}{ \textbf{Pre-conditions:}}
\newcommand{\PreConditions}{\item[\algorithmicprecondition]}
\newcommand{\algorithmicpostcondition}{ \textbf{Post-conditions:}}
\newcommand{\PostConditions}{\item[\algorithmicpostcondition]}
\newcommand{\algorithmicComment}{ \hspace{0.3cm}}
\newcommand{\Comment}{\item[\algorithmicComment]}

\newcommand{\revise}[1]{{#1}}

\newcommand{\TIKZDIR}{./tikz/}
\newcommand{\PNGDIR}{./png/}
\newcommand{\TABDIR}{./tab/}

\journalname{Journal of Scientific Computing}

\begin{document}

\title{
  A multiphase cubic MARS method for
  fourth- and higher-order interface tracking of
  two or more materials
  with arbitrary topology and geometry
  \thanks{This work was supported by
    the Fundamental Research Funds
    for the Central Universities 226-2025-00254
    and the National Natural Science Foundation of China
    (\#12272346 and \#11871429).
  Yan Tan and Yixiao Qian are co-first authors.}
} 

\titlerunning{Cubic MARS methods for multiphase IT}

\author{Y. Tan \and Y. Qian \and Z. Li \and Q. Zhang}

\institute{Y. Tan \and Y. Qian \and Z. Li \and Q. Zhang \at
  School of Mathematical Sciences,
  Zhejiang University, Hangzhou, Zhejiang, 310058, China
  \and Q. Zhang (Corresponding author)
  \email{qinghai@zju.edu.cn}
  \at Institute of Fundamental and Transdisciplinary Research,
  Zhejiang University, Hangzhou, Zhejiang, 310058, China}

\date{Received: date / Accepted: date}

\maketitle

\begin{abstract}
For interface tracking of an arbitrary number of materials
 in two dimensions,  
 we propose a multiphase cubic MARS method 
 that  \revise{
(a) represents the topology and geometry of the interface
 via graphs, cycles, and cubic splines,
(b) applies to any number of materials
 with arbitrarily complex topology and geometry, 
(c) maintains an $(r,h)$-regularity of the interface 
 so that the distance between any pair of adjacent markers
 is within a user-specified range,
(d) distributes the markers adaptively 
 along the interface 
 so that arcs with high curvature
 are resolved by densely populated markers, 
 and (e) achieves fourth-, sixth-, and eighth-order accuracy
 both in time and in space.}
In particular, all possible types of junctions, 
 which pose challenges to VOF methods and level-set methods,
 are handled with ease.
Results of a variety of benchmark tests confirm the analysis
 and demonstrate the superior accuracy, efficiency,
 and versatility of the proposed method.



\end{abstract}


\section{Introduction}
\label{sec:intro}
As a complex yet significant topic,
 multiphase flows 
 concern the simultaneous movements and interactions
 of a number of homogeneous \emph{materials} or \emph{phases}
 such as liquids, gases, and solids.
These flows are prevalent in natural and industrial processes
 yet pose major challenges to high-fidelity simulations in applied sciences.
One fundamental problem that accounts for these challenges
 is interface tracking (IT),
 the determination of regions occupied by these phases.

The most popular families of IT methods
 are probably level-set methods
 \cite{osher88:_front_propag_curvat_speed}, 
 front-tracking methods
 \cite{tryggvason01:_front_track_method_comput_multip_flow}, 
 and volume-of-fluid (VOF) methods 
 \cite{Hirt.Nichols_1981_volume}.
In level-set methods, 
 the interface is 
 \emph{implicitly} approximated
 as the zero isocontour of a signed distance function
 while, 
 in front-tracking methods, 
 it is \emph{explicitly} represented
 as a set of connected markers.
In VOF methods,
 the interface is not only implicitly described
 by volume fractions of the tracked phase 
 inside the fixed control volumes
 but also explicitly represented as a cellwise function. 
Within each time step,
 a VOF method consists of two substeps: 
 in the first reconstruction substep
 the explicit representation of material regions is determined
 \emph{solely} from volume fractions
 while in the second advection substep
 the volume fractions are advanced to the end of the time step 
 from the explicit representation and the velocity field. 

Previously,
 many IT methods have been developed for \emph{two-phase} flows, 
 where it is sufficient to track only one phase
 and deduce the region of the other.
Most of the aforementioned IT methods
 are second-order accurate for two-phase flows. 
The cubic MARS method \cite{zhang2018fourth}, 
 which belongs to none of the aforementioned three families, 
 even achieves fourth- and higher-order accuracy.
 
\begin{figure}
  \centering
  \subfloat[Three phases that constitute an elliptic disk]{
    \includegraphics[width=0.47\textwidth]{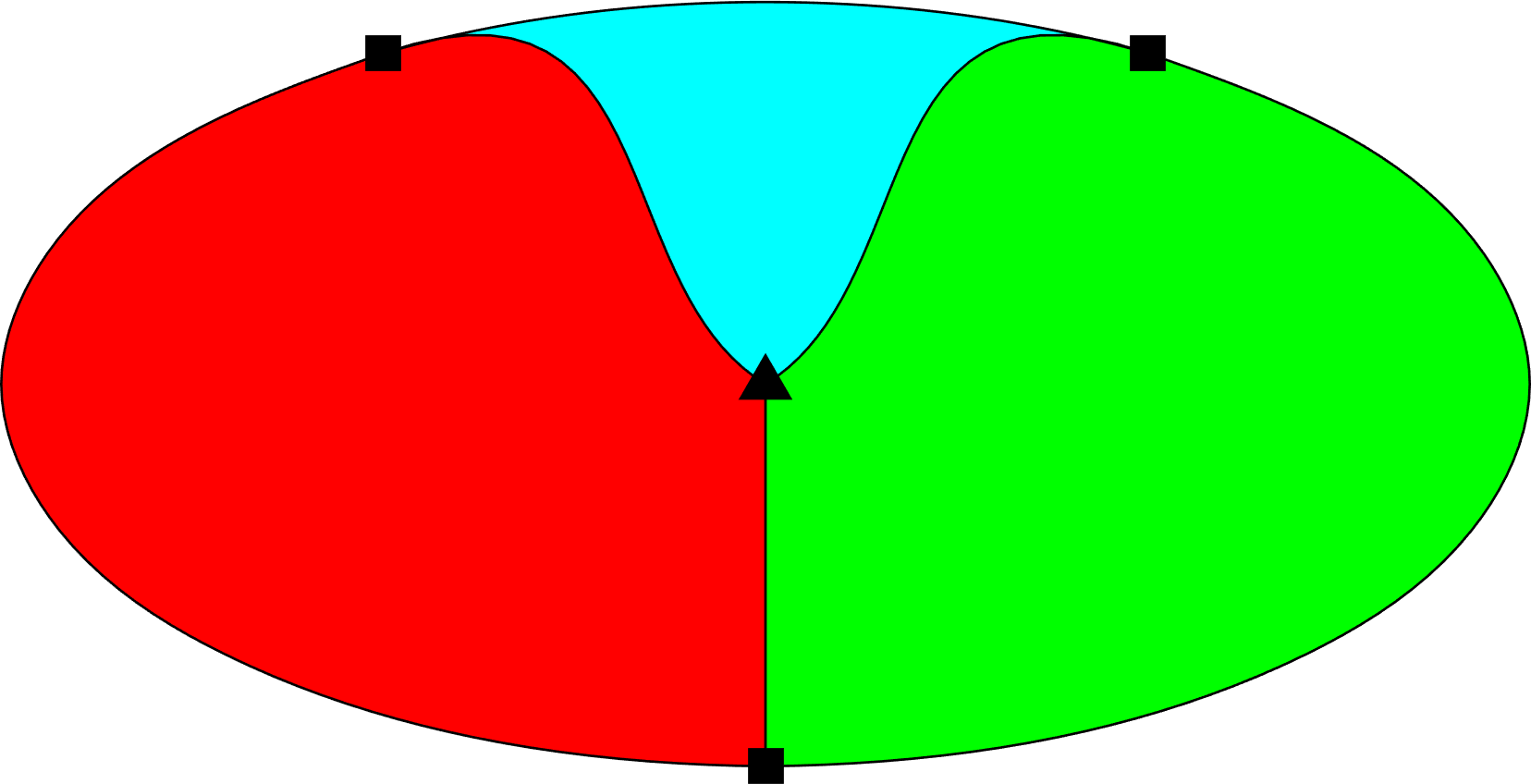}
  }
  \hfill
  \subfloat[Connecting markers on the interface with linear segments]{
    \includegraphics[width=0.47\textwidth]{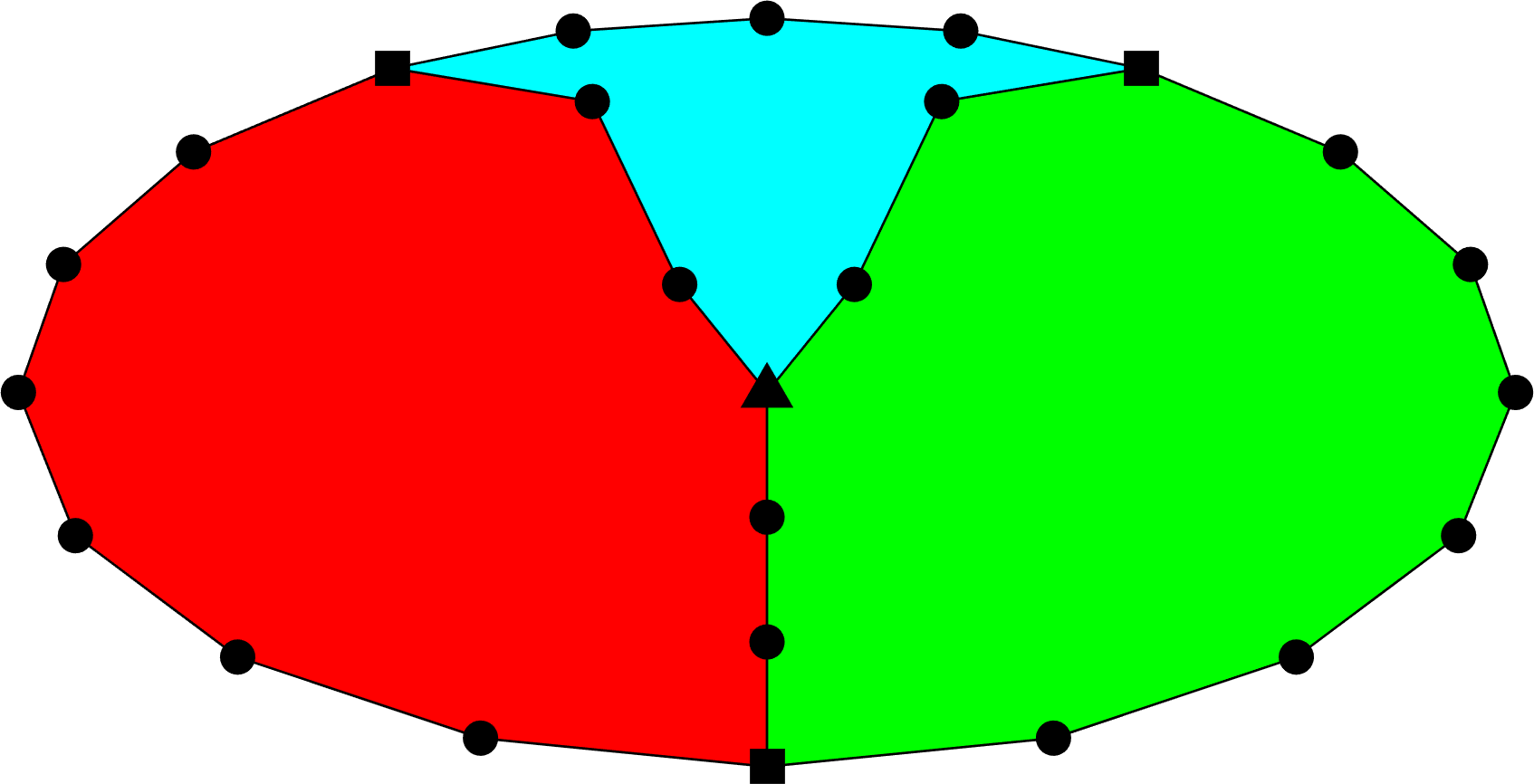}
  }

  \subfloat[Independent spline fitting 
  leads to overlaps and vacuums]{
    \includegraphics[width=0.47\textwidth]{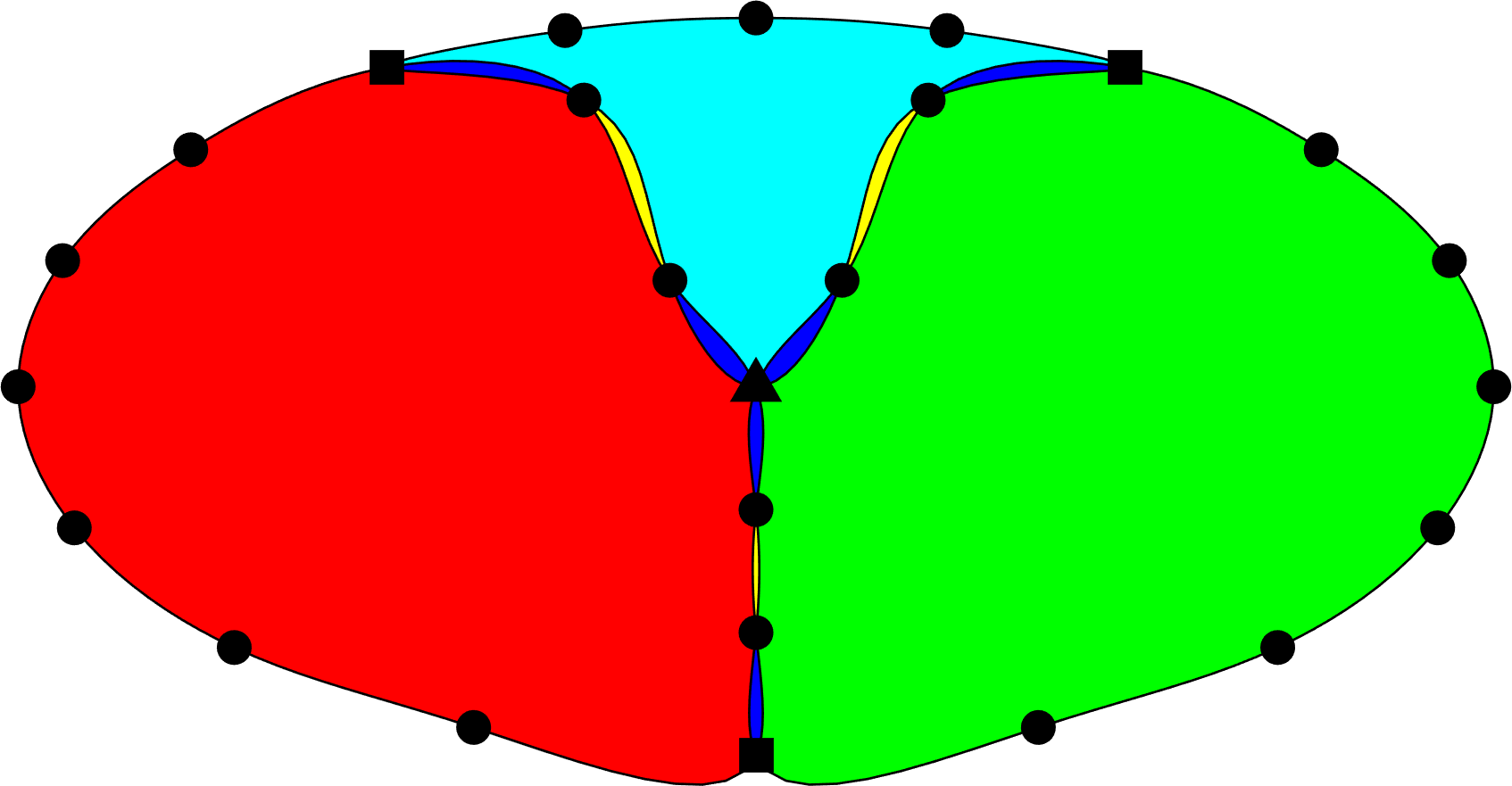}
  }
  \hfill
  \subfloat[Neither overlaps nor vacuums
  exist if each curve is fitted only once]{
    \includegraphics[width=0.47\textwidth]{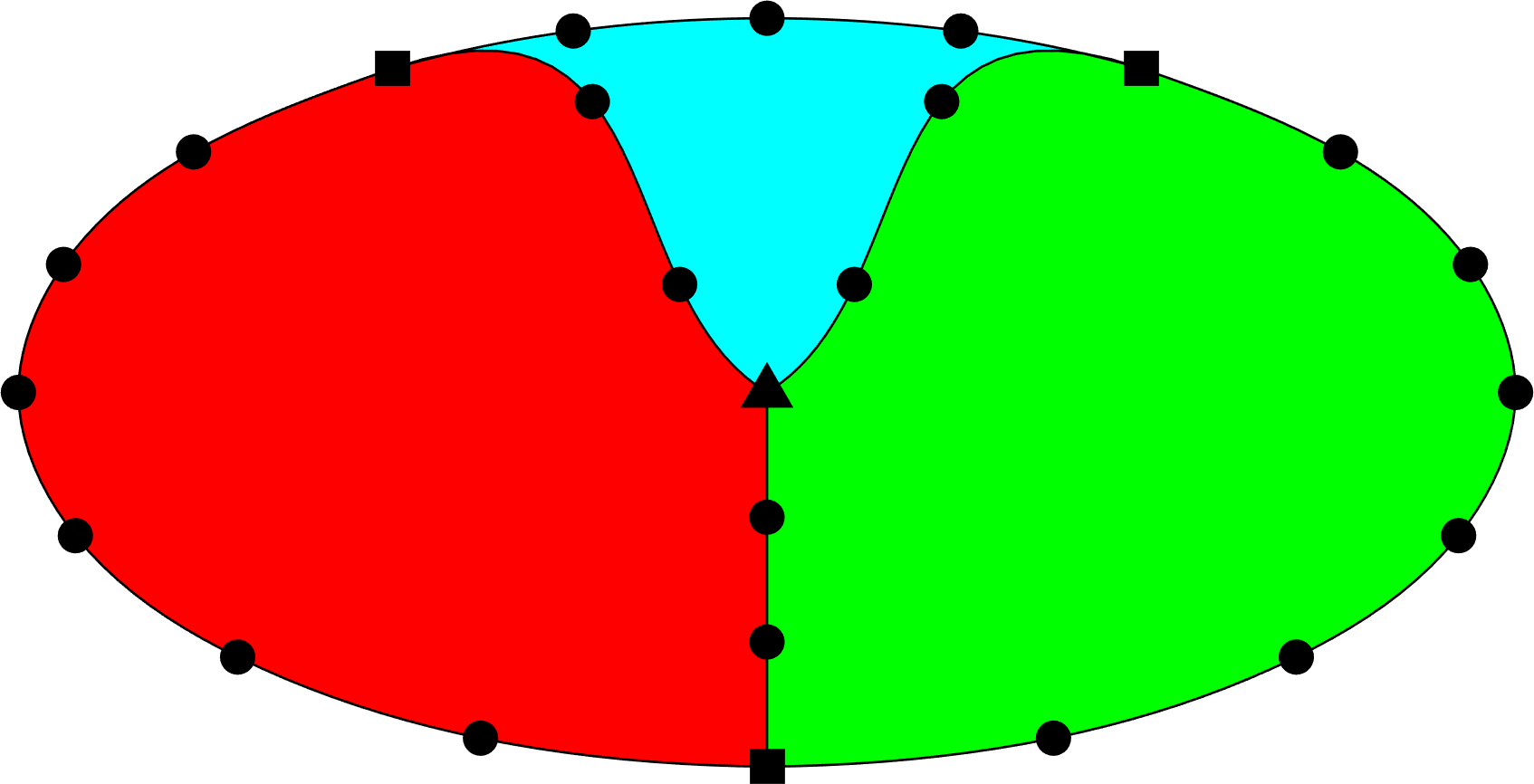}
  }

  \caption{Boundary representations of multiple phases.
    Subplot (a) shows three adjacent phases to be represented
    by a set of interface markers or characteristic points
    (squares, triangles, and dots).
    \revise{At each of the \emph{T junction}s (squares), 
    there exist two curve segments forming a smooth curve
    whereas, at the \emph{Y junction} (the triangle), 
    all curves formed by connecting two radial curve segments
    can only be ${\cal C}^0$.} 
    In subplot (b),
    the markers are connected by linear segments,
    yielding a second-order representation 
    without creating overlaps and vacuums between adjacent phases.
    In subplot (c),
    fitting ${\cal C}^2$ periodic cubic splines 
    independently for each phase
    leads to overlaps (blue areas)
    and vacuums (white areas inside the ellipse).
    \revise{In particular, the bottom T junction
    is qualitatively changed to a Y junction.}
    Therefore, in subplot (d),
    the smooth ellipse is represented
    by a ${\cal C}^2$ periodic cubic spline
    while the three radial curve segments incident to the Y junction
    are approximated separately by not-a-knot cubic splines; 
    see Definition \ref{def:splineTypes}.
  }
  \label{fig:independent-mars}
\end{figure}

In contrast, 
 for IT of three or more phases, 
 the literature is much thinner
 and the accuracy of current IT methods is much lower; 
this case is called the \emph{multiphase IT problem} or
 the \emph{IT problem of multiple materials}, 
 since more than one phase \emph{has to} be tracked.
The core difficulty in tracking multiple phases, however, 
 lies not in the number of phases
 but in their topology and geometry 
 that are fundamentally more complicated than those of two-phase flows.
For example,
 an interface curve might have a \emph{kink},
 i.e., a ${\cal C}^1$ discontinuity of the curve;
 \revise{see Definition \ref{def:kink} for the more general concept of
 non-smooth points}.  
Kinks are problematic for level-set methods and VOF methods: 
 large reconstruction errors at these kinks
 are propagated along the interface in subsequent time steps
 by numerical diffusion, 
 altering geometric features of sharp corners to rounded shapes.
As another example,
 three or more phases might meet at a \emph{junction} 
 (see Definition \ref{def:junctionPoint}
 \revise{and Fig.~\ref{fig:independent-mars}}), 
 where the boundary curve of at least one phase
 contains kinks.
\revise{The tracking of such a junction or kink with low accuracy
 causes more damage to the simulation fidelity
 of multiphase flows than to that of two-phase flows, 
 because
 a junction is usually a place of primary interest 
 where important physics occurs. 
For example,
 accurately capturing the position and dynamics
 of the air-water-solid junctions
 is critical for the study of 
 moving contact lines
 \cite{sui14:_numer_simul_flows_movin_contac_lines,zhao2021}
 and the droplet impacting problem \cite{liu2020,wang2023}.
Yet this remains challenging for conventional level-set and VOF methods.
}
Indeed,
 a core difficulty
 in both two-phase and multiphase IT problems 
 is the handling of these junctions and kinks, 
 for which
 the mathematical modeling
 and highly accurate and efficient algorithms
 are the main focuses of this work.
 
The standard level-set construction is not applicable
 to the local neighborhood of a junction
 because the zero level set of a single $\mathcal{C}^1$
 signed distance function
 is never homeomorphic to the one-dimensional CW complex
 that characterizes the topology at the junction.
To resolve this difficulty,
 Saye and Sethian \cite{saye11:_voron}
 propose the Voronoi implicit interface method
 as a generalization of the level-set method
 for computing multiphase physics,
 via an elegant extension of the Voronoi diagram
 to a set of curves and surfaces.
Unfortunately, 
 the Voronoi diagram that determines the interface
 also limits this method
 to the first-order accuracy at the junctions. 

The extension of VOF methods to multiple phases
 has been primarily focused on the reconstruction substep
 \cite{benson2002volume}.
Piecewise linear VOF reconstruction schemes for junctions 
 are limited to triple points
 \cite{choi2007triplePoint,caboussat2008triplePoint} 
 and it appears that no VOF schemes
 handle junctions with four or more incident edges.
In the ``onion-skin'' model,
 the multiple materials inside a control volume
 are assumed to have a layered topology \emph{with no junctions}. 
Given a material ordering,
 the interface between the $i$th and the $(i+1)$th phases
 is defined by applying a VOF reconstruction scheme
 to the union of materials $1$ through $i$.
When two reconstructed interfaces intersect,
 one either adjusts them to eliminate the intersection
 \cite{deniem2007volume,sijoy2010volume}
 or scales the fluxes to account for volumes of overlapping areas
 \cite[p. 365]{benson1992computational}.
As such, the IT results
 depend substantially on the material ordering.
Youngs \cite{youngs1982:_time} requires the user
 to specify a priority list, 
 which is then utilized as a \emph{static} material ordering for each cell.
Mosso and Clancy \cite{mosso1995geometrically}
 propose to order the materials \emph{dynamically}
 based on estimates of their centroids in each cell.
Benson \cite{benson1998eulerian}
 adds the estimated centroids as solution variables
 and determines the dynamic ordering
 by a least-squares fitting of a line to the centroids
 and then sorting the projected images of centroids along the line. 

For material-order-dependent VOF methods, 
 an incorrect ordering may result in large errors in reconstruction
 and premature/belated advection of multiple phases
 \cite{kucharik2010:_ale}.
In addition, 
 the topology of a junction might be
 changed by the numerical diffusion in these methods; 
 see, e.g., the erroneous alteration of an X junction
 to two T junctions
 illustrated in \cite[Fig. 16]{schofield2009second}.
To alleviate these adverse effects,
 Schofield et al. \cite{schofield2009second} develop
 a power diagram method, 
 a material-order-independent interface reconstruction technique, 
 in which the interface is first reconstructed
 by a weighted Voronoi diagram 
 from material locator points
 and then improved by minimizing an objective function
 that smooths the interface normals.

Another extension of VOF methods  
 is the moment-of-fluid (MOF) method \cite{dyadechko08:_recon}, 
 which reconstructs cellwise materials 
 not only by volume fractions (their 0th moments)
 but also by centroids (their 1st moments). 
Since these two moments already provide
 enough information to construct a linear function, 
 no data from neighboring cells are needed. 
This independence furnishes
 a straightforward generalization
 from two phases to $N_p$ phases, 
 via enumerating all $N_p!$ possible orderings
 to minimize the error norm of the first moment.
Despite being material-order-dependent, 
 the MOF method is second-order accurate
 if the true interface is \emph{${\cal C}^2$-serial},
 i.e., if all phases can be sequentially separated from the bulk
 by ${\cal C}^2$ curves \cite{dyadechko08:_recon}. 
For example,
 the interface in Fig.~\ref{fig:independent-mars}(a)
 is ${\cal C}^2$-serial at the T junctions on the ellipse, 
 but not so at the Y junction inside the ellipse, 
 where the MOF reconstruction is only first-order accurate.
See \cite{kucharik2010:_ale} for an accuracy comparison 
 between MOF and VOF methods. 

Interestingly, multiphase MOF reconstruction 
 is helpful for two-phase flows in capturing filaments,
 thin strands of one material surrounded by another within a cell, 
 e.g., the tail tips
 in Fig.~\ref{fig:test-quartered-circular-disk-T4}(h,k).
Jemison et al. \cite{jemison2015:_filament}
 introduce a fictitious phase
 to reformulate the filament reconstruction
 as three materials separated by two interfaces
 in an onion-skin topology.
Hergibo et al. \cite{hergibo2023moment}
 resolve filaments via a symmetric multi-material approach
 with polygon clipping. 
These multiphase MOF methods reconstruct filaments more accurately
 than the standard MOF for two phases; 
 see Table \ref{tab:comparisonVortex}(b,c).
 
As far as we know, 
 neither level-set methods nor VOF/MOF methods 
 are capable of reconstructing
 \emph{all} types of junctions to second-order accuracy; 
in particular, 
 they all drop to first-order accuracy at Y junctions.
A front-tracking method,
 with the aid of graphs
 such as that in Definition \ref{def:topoGraph}, 
 can represent the interface topology exactly
 and thus achieves full second-order accuracy,
 so long as all junctions and kinks 
 are already selected as interface markers. 
In particular,
 such a front-tracking method is independent of material ordering.
Even if interface markers are tracked twice,
 neither overlaps nor vacuums
 are created between adjacent phases
 provided that the markers are connected with \emph{linear} segments; 
 see Fig.~\ref{fig:independent-mars}(b).
However,
 this statement does not hold for higher-order splines: 
 to achieve an accuracy higher than the second order, 
 one needs the geometric information at each junction 
 on the pairing of smoothly connected curve segments.
As shown in Fig.~\ref{fig:independent-mars}(c), 
 independent approximations
 of the boundary Jordan curves of each phase
 with ${\cal C}^2$ periodic cubic splines
 result in overlaps and vacuums between adjacent phases.
In contrast,
 a blend of periodic and not-a-knot cubic splines
 fitted with due considerations
 of topological structures and geometric features 
 gives satisfactory results; 
 see Fig.~\ref{fig:independent-mars}(d).

For traditional IT methods, 
 why is it so difficult to achieve high-order accuracy
 for multiple materials?
In our humble opinion,
 the reason is that
 \emph{topology and geometry are avoided in these methods
 via converting topological and geometric problems in IT
 to numerical solutions of differential equations}, 
 such as the ordinary differential equations (ODEs)
 of interface markers in front-tracking methods
 and the scalar conservation laws in level-set and VOF methods. 
Being relinquished at the very beginning,
 key topological structures and geometric features
 can hardly be recovered to high-order accuracy
 in subsequent time steps. 
Indeed,
 MOF methods
 \cite{dyadechko08:_recon,jemison2015:_filament,hergibo2023moment}
 are at best second-order accurate, 
 even after utilizing more geometric information.
For front-tracking methods to achieve an accuracy
 higher than the second order, 
 the connected markers
 must be supplemented with additional topological information. 

\revise{Another fundamental challenge in explicit IT
  is how to maintain geometric regularity of interface markers.
  While traditional front-tracking methods 
  avoid this problem 
  by calculating interface-grid intersections, 
  those aiming for subgrid resolutions often 
  suffer from irregular marker distributions under large deformations.
The methods of artificial tangential 
 velocity~\cite{hou1994,ceniceros2003,barrett2007,hu2022evolving}
 achieve uniform distribution of markers 
 by moving them along the tangential direction.
However, they seem to be unable to adapt the marker density
 according to the local curvature of the evolving curve. 
Consequently, a strongly curved arc might cause 
 a much greater IT error than mildly curved arcs. 
}

The above discussions motivate questions as follows. 
\begin{enumerate}[label=(Q-\arabic*), leftmargin=2.7em]
\item Given their physical significance, 
  can kinks and junctions of all possible types
  be faithfully represented and accurately tracked
  without creating overlaps and vacuums between adjacent phases?
  \label{enum:q-kinks}
\item VOF and level-set methods cannot preserve geometric
  features under isometries, 
  nor can they preserve topological structures
  under homeomorphisms. 
  To resolve these difficulties,
  can we develop an IT method
  that respectively preserves,
  under isometric and homeomorphic flow maps, 
  the geometric features and topological structures of each phase?
  \label{enum:q-data-structure}
\item 
  Can the distance between adjacent markers be made 
  a function of the local curvature
  so that the IT error is roughly uniform along the interface? 
  \label{enum:arms-strategy}
\item Can we design an efficient fourth-order IT method
  for tracking 
  two or more phases with arbitrary topology and geometry?
  \label{enum:shared-interfaces}
\end{enumerate}

In this paper,
 we provide positive answers to all above questions. 
Fundamentally different from that of current IT methods,
 our primary principle is to
 \emph{tackle topological and geometric problems in IT 
   with tools in topology and geometry.}
 
Previously, we have 
 proposed a topological space (called the Yin space)
 as a mathematical model
 of two-dimensional continua \cite{zhang2020:_boolean}, 
 analyzed explicit IT methods under the framework
 of mapping and adjusting regular semianalytic sets (MARS)
 \cite{zhang2016mars}, 
 developed a cubic MARS method for two-phase flows \cite{zhang2018fourth},
 and augmented MARS methods to curve shortening flows 
 via the strategy of adding and removing markers on splines (ARMS)
 \cite{hu2025:_arms}. 
As an extension of MARS to multiple materials,
 this work is another manifestation that
 IT methods coupling
 (even elementary) concepts in topology and geometry
 can be highly accurate and highly efficient.

The main contributions of this work are 
\begin{enumerate}[label=(C-\arabic*)]
\item the mathematical models and data structures
  for representing an arbitrary number of materials
  with arbitrary topology and geometry, 
\item the extension of the MARS framework
  \cite{zhang2016mars,zhang2020:_boolean}
  to the general scenario of multiple phases
  with junctions and kinks, 
\item a multiphase cubic MARS method
  for solving the multiphase IT problem in Sec.~\ref{sec:generalIT}
  with fourth- and higher-order accuracy. 
\end{enumerate}
 
The rest of this paper is structured as follows. 
Sec.~\ref{sec:yinsets} is a brief review on the Yin space,
 with all topology of a two-dimensional continuum 
 classified in Theorem \ref{thm:uniqueCases}.
Sec.~\ref{sec:generalIT} is a precise definition
 of the multiphase IT problem.
In Sec.~\ref{sec:representation},
 we answer (Q-1,2)
 by designing concepts and data structures
 for representing \emph{static} multiple phases
 and by separating their topology from the geometry.
In particular,
 the accuracy and stability of periodic splines and not-a-knot splines
 in respectively approximating closed curves and curve segments
 are meticulously examined.
 \revise{A main difficulty in answering (Q-3) by the ARMS strategy
   is how to maintain the permanent status of some special markers
   that characterize the topology and geometry
   of the evolving interface. 
In Sec.~\ref{sec:algorithms},
 we resolve this difficulty 
 by augmenting the original ARMS strategy
 \cite[Def. 3.5]{hu2025:_arms}
 to a more general ARMS strategy in Definition \ref{def:ARMS}.} 
Then 
 we propose in Definition \ref{def:MARS-for-multiple-phases}
 the multiphase cubic MARS method as our answer to \mbox{(Q-4)}.
In Sec.~\ref{sec:tests},
 we demonstrate the fourth-, sixth-, and eighth-order accuracy
 of the proposed method
 by performing various classic and new benchmark tests.
Results of these tests show
 that the proposed method is more accurate
 than state-of-the-art IT methods
 by many orders of magnitude.
Finally, we conclude this paper
 in Sec.~\ref{sec:conclusion}
 with several future research prospects.



\section{Modeling continua by Yin sets}
\label{sec:yinsets}


In a topological space ${\mathcal X}$,
 the \emph{complement} of a subset ${\mathcal P}\subseteq {\mathcal X}$,
 written ${\mathcal P}'$,
 is the set ${\mathcal X}\setminus {\mathcal P}$.
The \emph{closure} of a set ${\mathcal P}\subseteq{\mathcal X}$,
 written $\overline{\mathcal P}$,
 is the intersection of all closed 
 supersets of ${\mathcal P}$.
The \emph{interior} of ${\mathcal P}$, written ${\mathcal P}^{\circ}$,
 is the union of all open subsets of ${\mathcal P}$.
The \emph{exterior} of ${\mathcal P}$,
 written ${\mathcal P}^{\perp}:= {\mathcal P}^{\prime\circ}
 :=({\mathcal P}')^{\circ}$,
 is the interior of its complement.
A point $\mathbf{x}\in {\mathcal X}$ is
 a \emph{boundary point} of ${\mathcal P}$
 if $\mathbf{x}\not\in {\mathcal P}^{\circ}$
 and $\mathbf{x}\not\in {\mathcal P}^{\perp}$.
The \emph{boundary} of ${\mathcal P}$, written $\partial {\mathcal P}$,
 is the set of all boundary points of ${\mathcal P}$.
It can be shown that
 ${\mathcal P}^{\circ}={\mathcal P}\setminus \partial {\mathcal P}$
 and
 $\overline{\mathcal P}= {\mathcal P}\cup \partial {\mathcal P}$.

A \emph{regular open} set 
 is an open set ${\cal P}$ satisfying
 ${\mathcal P}=\overline{\mathcal P}^{\circ}$
 while a \emph{regular closed} set
 is a closed set ${\cal P}$
 satisfying ${\mathcal P}=\overline{\mathcal P^{\circ}}$.
Regular sets, open or closed,
 capture a key feature of continua
 that their regions are free of lower-dimensional elements
 such as isolated points and curves in $\mathbb{R}^2$
 and dangling faces in $\mathbb{R}^3$.
The intersection of two regular sets, however,
 might contain an infinite number of connected components
 \cite[eqn~(3.1)]{zhang2020:_boolean},
 making it difficult to perform Boolean algorithms on regular sets 
 since no computer has an infinite amount of memory. 
This difficulty is resolved by requiring each regular set
 to be simultaneously a \emph{semianalytic} set,
 i.e., a set
  ${\cal S}\subseteq \mathbb{R}^{\Dim}$
 in the universe of 
 a finite Boolean algebra formed from the sets
 \mbox{${\cal X}_i
 =\left\{{\mathbf{x}}\in\mathbb{R}^{\Dim}:g_i({\mathbf{x}})\geq 0\right\}$}
 where each $g_i:\mathbb{R}^{\Dim}\to\mathbb{R}$
 is an analytic function.
Intuitively,
 $\partial {\cal S}$ is piecewise ${\cal C}^{\infty}$ 
 so that ${\cal S}$ can be described by a finite number of entities.

\begin{defn}[Yin space \cite{zhang2020:_boolean}] 
  \label{def:YinSet}
  A \emph{Yin set} $\mathcal{Y}\subseteq \mathbb{R}^{\Dim}$
  is a regular open semianalytic set 
  whose boundary is bounded. 
  All Yin sets form the \emph{Yin space}
  $\mathbb{Y}$.
\end{defn}

\revise{The reader is referred to
  \cite{zhang2020:_boolean} for more physical intuitions
  and geometric illustrations
  \cite[Figures 1,5,8,10]{zhang2020:_boolean} of Yin sets.
}

A \emph{curve (segment)} is (the image of) a continuous map
 $\gamma: [0,1]\to\mathbb{R}^2$;
 it is \emph{closed}
 if its \emph{endpoints} coincide, i.e., $\gamma(0)=\gamma(1)$. 
The \emph{open curve} of a curve segment $\gamma$
 is its restriction  $\gamma|_{(0,1)}$,
 whose endpoints are those of $\gamma$.
 \revise{The curve segment of an open curve
   is obtained by appending the two endpoints to the open curve
   while the open curve of a curve segment
   by removing the two endpoints.
   Hence there is a one-to-one correspondence between
   these two concepts.
 }
An open curve is \emph{simple} if it is injective.
A curve is \emph{Jordan}
 if it is closed and its corresponding open curve is simple.
%
The \emph{interior of an oriented Jordan curve},
 written $\Int(\gamma)$,
 is the component of $\mathbb{R}^2\setminus\gamma$ 
 that always lies to the left 
 of the observer who traverses $\gamma$ according to $\gamma([0,1])$.
A Jordan curve $\gamma$ is \emph{counterclockwise}
 or \emph{positively oriented}
 if $\Int(\gamma)$ is the bounded component
 of $\mathbb{R}^2\setminus\gamma$;
 otherwise it is \emph{clockwise} or \emph{negatively oriented}.

Following \cite[Def. 3.7]{zhang2020:_boolean},
 we call two Jordan curves \emph{almost disjoint}
 if they have no proper intersections (i.e., crossings) 
 and the number of their improper intersections is finite. 
A Jordan curve $\gamma_k$
 is said to \emph{include} another Jordan curve $\gamma_{\ell}$,
 written $\gamma_{k}\ge \gamma_{\ell}$ or $\gamma_{\ell}\le \gamma_{k}$,
 if the bounded complement of $\gamma_{\ell}$
 is a subset of that of $\gamma_{k}$.
If $\gamma_k$ includes $\gamma_{\ell}$ and
 $\gamma_k\ne \gamma_{\ell}$,
 we write $\gamma_{k}> \gamma_{\ell}$
 or $\gamma_{\ell}< \gamma_{k}$.
In a partially ordered set (poset) ${\mathcal J}$ of Jordan curves
 with inclusion as the partial order, 
we say that 
 $\gamma_k$ \emph{cover}s $\gamma_{\ell}$ in ${\mathcal J}$
 and write `$\gamma_{k} \succ \gamma_{\ell}$'
 or `$\gamma_{\ell}\prec \gamma_{k}$'
 if $\gamma_{\ell}< \gamma_{k}$ and no element $\gamma\in {\mathcal J}$
 satisfies $\gamma_{\ell}< \gamma<\gamma_{k}$.

In Definition \ref{def:YinSet}, 
a regular open set instead of a regular closed set is employed
because the former can be \emph{uniquely} represented
by its boundary Jordan curves
while the latter cannot \cite[Fig. 5]{zhang2020:_boolean}.
 
\begin{thm}[Global topology and boundary representation 
  of connected Yin sets \cite{zhang2020:_boolean}]
  \label{thm:uniqueCases}
  The boundary of any connected Yin set 
  \mbox{${\cal Y}\ne \emptyset, \mathbb{R}^2$}
  can be uniquely partitioned into 
  a finite set of pairwise almost disjoint Jordan curves, 
  which can be uniquely oriented to yield
  a unique representation of ${\cal Y}$ as 
  ${\cal Y} = \bigcap_{\gamma_j\in {\cal J}_{\partial {\cal Y}}}
  \Int(\gamma_j)$
  where
  ${\mathcal J}_{\partial {\mathcal Y}}$,
  the set of oriented boundary Jordan curves of ${\mathcal Y}$,
  must be one of the two types,
  \begin{equation}
    \label{eq:decomTypes}
    \renewcommand{\arraystretch}{1.2}
    \left\{
      \begin{array}{ll}
        {\mathcal J}^-
        =\{\gamma^-_1, \gamma^-_2, \ldots, \gamma^-_{n_-}\}
        & \text{where}\ n_-\ge 1,
        \\
        {\mathcal J}^+
        =\{\gamma^+,\gamma^-_1, \gamma^-_2, \ldots, \gamma^-_{n_-}\}
        & \text{where}\ n_-\ge 0,
      \end{array}
    \right.
  \end{equation}
  and all $\gamma^-_j$'s are negatively oriented
  and mutually incomparable with respect to inclusion.
  In the case of ${\mathcal J}^+$,
  $\gamma^+$ covers $\gamma^-_j$,
  i.e., $\gamma_j^- \prec \gamma^+$ holds 
  for each $j=1,2,\ldots,n_-$.
\end{thm}

A form ${\cal J}^-$ or ${\cal J}^+$ implies
 that the connected Yin set ${\cal Y}$ is unbounded or bounded, 
 respectively.
In Fig.~\ref{fig:multiphase-representation}(a), 
 ${\cal M}_6$ is unbounded while
 all other connected Yin sets are bounded; 
 $n_{-}=2$ for ${\cal M}_{4,1}$ and ${\cal M}_{5}$
 and $n_{-}=0$ for ${\cal M}_{1}$,
 ${\cal M}_{2}$, ${\cal M}_{3}$, and ${\cal M}_{4,2}$.

\begin{thm}[Boolean algebra on the Yin space \cite{zhang2020:_boolean}]
  \label{thm:YinSetsFormABooleanAlgebra}
  The universal algebra 
  $\mathbf{Y}:=
  ({\mathbb Y},\ \cup^{\perp\perp},\ \cap,\ \, ^{\perp},\
  \emptyset,\ \mathbb{R}^2)$
  is a Boolean algebra,
  where the regularized union
  is given by ${\mathcal Y}\cup^{\perp\perp}{\mathcal Z}
  := ({\mathcal Y}\cup {\mathcal Z})^{\perp\perp}$
  for all ${\mathcal Y},{\mathcal Z}\in{\mathbb Y}$.
\end{thm}

The uniqueness of the boundary representation of Yin sets
 in Theorem \ref{thm:uniqueCases}
 implies that $\mathbb{Y}$ and $\mathbb{J}$ are \emph{isomorphic},
 written $\mathbb{Y}\cong \mathbb{J}$,
 where $\mathbb{J}$ is the \emph{Jordan space}
 of posets of oriented Jordan curves.
This isomorphism is exploited in \cite{zhang2020:_boolean}
 to reduce the above Boolean algebra
 to calculating intersections of boundary Jordan curves.

\revise{
On top of its topological and algebraic structures, 
 the Yin space can be further augmented to a metric space. 
To this end, the volume of a Yin set $\mathcal{Y}$ is given by
\begin{equation}
  \label{eq:int}
  \begin{array}{l}
    \Vert \mathcal{Y} \Vert :=
    \left|\int_{\mathcal{Y}}d \mathbf{x}\right|,
  \end{array}
\end{equation}
where the integral can be interpreted as a Riemann integral
since $\mathcal{Y}$ is semianalytic.
The regularized symmetric difference
 $\oplus:\mathbb{Y}\times\mathbb{Y}\rightarrow\mathbb{Y}$
 is defined as 
\begin{equation}
  \label{eq:oplus}
  \mathcal{P} \oplus \mathcal{Q} := \left( \mathcal{P}\,\backslash\,
    \mathcal{Q} \right) \,
  \mathop{\cup}\nolimits^{\perp\perp} \,
  \left( \mathcal{Q} \,\backslash \,\mathcal{P} \right),
\end{equation}
which satisfies
\begin{align}
  \label{eq:oplusconst}
  \forall \mathcal{Y} \in \mathbb{Y},&\quad \mathcal{Y}\oplus \mathcal{Y} =
  \emptyset, \quad\emptyset \oplus \mathcal{Y} = \mathcal{Y}; 
  \\
  \label{eq:triinequ}
  \forall \mathcal{P}, \mathcal{Q} \in \mathbb{Y},&\quad \Vert \mathcal{P}
  \oplus \mathcal{Q} \Vert \le \Vert \mathcal{P}\Vert + \Vert \mathcal{Q}\Vert.
\end{align}

It follows from \eqref{eq:int}, \eqref{eq:oplus},
\eqref{eq:oplusconst}, and \eqref{eq:triinequ}
that $(\mathbb{Y}, d)$ forms a metric space
where the metric $d \,: \, \mathbb{Y} \times \mathbb{Y}
\rightarrow [0, +\infty)$ is 
\begin{equation}
  \label{eq:metric}
  \forall \mathcal{P}, \mathcal{Q} \in \mathbb{Y},\quad
  d(\mathcal{P}, \, \mathcal{Q})\, := \, \Vert \mathcal{P} \oplus
  \mathcal{Q} \Vert.
\end{equation}
}




\section{The multiphase IT problem}
\label{sec:generalIT}
%
%
For any given initial time $t_0$
 and initial position $p_0\in \mathbb{R}^{\Dim}$, 
 the ODE 
\begin{equation}
  \label{eq:velODE}
  \frac{\dif\, \mathbf{x}}{\dif\, t} = \mathbf{u}(\mathbf{x},t)
\end{equation}
 admits a unique solution
 if the time-dependent velocity field $\mathbf{u}(\mathbf{x},t)$
 is continuous in time and Lipschitz continuous in space.
This uniqueness gives rise to a flow map
 $\phi: \mathbb{R}^{\Dim}\times \mathbb{R}\times
 \mathbb{R}\rightarrow \mathbb{R}^{\Dim}$
 that takes the initial position $p_0$
 of a Lagrangian particle $p$,
 the initial time $t_0$, 
 and the time increment $\tau$
 and returns $p(t_0+\tau)$,
 the position of $p$ at time $t_0+\tau$:
\begin{equation} 
  \label{eq:tracing}
  \phi_{t_0}^{\tau}(p) := p(t_0+\tau) 
  = p(t_0) + \int_{t_0}^{t_0+\tau} \mathbf{u}(p(t),t)\,\dif\, t.
\end{equation}
The flow map also generalizes to arbitrary point sets
 in a straightforward way, 
 \begin{equation}
   \label{eq:flowMapActsOnYinSet}
   \phi_{t_0}^{\tau}({\mathcal Y})  
   = \{\phi_{t_0}^{\tau}(p) : p\in {\mathcal Y}\}.
\end{equation}
If we further restrict the above point set to a Yin set,
 then the flow map $\phi$ for given $t_0$ and $\tau$
 can be considered as a unitary operation
 $\phi_{t_0}^{\tau} : \mathbb{Y} \rightarrow \mathbb{Y}$.
It is not difficult \cite[p. 6]{arnold73:_ordin_differ_equat}
 to show
 
\begin{lem}
  \label{lem:diffeomorphism}
  For fixed $t_0$ and $\tau$,
  the flow map $\phi_{t_0}^{\tau}$ 
  in (\ref{eq:flowMapActsOnYinSet}) is a diffeomorphism,
  i.e., a ${\cal C}^1$ bijection
  whose inverse is also ${\cal C}^1$. 
\end{lem}

In the IT problem,
 we are usually given \emph{a priori}
 a velocity field $\mathbf{u}(\mathbf{x},t)$,
 by which each fluid phase is passively advected.
It is via this action of flow maps upon the Yin space
that we formulate

\revise{
\begin{defn}[Multiphase IT]
  \label{def:interfaceTracking}
  For $N_p$ pairwise disjoint moving Yin sets 
  ${\cal M}(t):=\left\{\mathcal{M}_i(t): i=1,\ldots, N_p\right\}$, 
  the \emph{multiphase IT problem}
  is to determine ${\cal M}(T)$ at $T > t_0$
  from ${\cal M}(t_0)$ and a family of diffeomorphic flow maps
  \begin{equation}
    \label{eq:groupOfDiffeomorphisms}
    \Phi(t_0,T):=\{\phi_{t_s}^{\tau}:\
    t_s\in[t_0,T);\, t_s+\tau\in (t_0,T]\}, 
  \end{equation}
  each of which acts on ${\cal M}(t)$
  by (\ref{eq:flowMapActsOnYinSet}). 
\end{defn}
}

Definition \ref{def:interfaceTracking}
 extends the IT problem for a single phase
 in \cite[Def. 3.1]{zhang2018fourth}. 
This extension is \emph{theoretically} trivial
 in that the exact flow map can be applied
 to the Yin sets in any order
 to produce $\mathcal{M}(T)$.
However, challenges of the multiphase IT problem
 mostly lie in the \emph{computational} aspects such as 
 the simultaneous preservation of high-order accuracy, 
 phase adjacency, topological structures, and geometric features.

The setup of the multiphase IT problem in Definition \ref{def:interfaceTracking}
 does not allow topological changes, 
 since they are precluded by
 the diffeomorphic flow maps in Definition \ref{def:interfaceTracking}.
Although in this work we confine ourselves
 to diffeomorphic movements for each phase, 
 the static representation of multiple Yin sets
 in Sec.~\ref{sec:representation} 
 and the ARMS strategy for moving splines
 in Sec.~\ref{sec:algorithms}
 lay a solid ground for future treatments 
 of topological changes.
 
 


\section{Boundary representation
  of static Yin sets}
\label{sec:representation}
The $N_p$ phases are identified with
 a set ${\cal M}:=\left\{\mathcal{M}_i: i=1,\ldots, N_p\right\}$ 
 of pairwise disjoint Yin sets, 
 each of which is 
 $\mathcal{M}_i := \bigcup^{\perp\perp} _{j} \mathcal{M}_{i,j}$ 
 where $\mathcal{M}_{i,j}$ is
 the $j$th connected component of the $i$th phase ${\cal M}_i$. 
If ${\cal M}_i$ is connected,
 we simply write ${\cal M}_i$ for ${\cal M}_{i,1}$. 
\revise{
  Thanks to Theorem \ref{thm:uniqueCases},
  it suffices to represent each ${\cal M}_{i,j}(t)$
  by its boundary $\partial{\cal M}_{i,j}(t)$. 
 }

\begin{ntn}
  \label{ntn:B-REP}
  The $(i,j)$th \emph{poset of oriented Jordan curves}
  of ${\cal M}$
  is denoted by $\Gamma_{i,j}:=\{\gamma_{i,j}^k\}$
  such that
  $\mathcal{M}_{i,j} =
  \cap _{\gamma_{i,j}^k \in \Gamma_{i,j}}
  \Int(\gamma_{i,j}^k)$.
  Denote by $N_{{\cal M}_i}$ 
  the number of connected components of ${\cal M}_i$
  and we write
  \begin{equation}
    \label{eq:B-REP}
    \begin{array}{l}
      \Gamma_i:=\{\Gamma_{i,j}: j=1,\ldots, N_{{\cal M}_i}\},\quad
      \Gamma:= \{\Gamma_i: i=1,\ldots,N_p\}; 
      \\
      \chi(\Gamma_{i,j}):= \bigcup_{\gamma_{i,j}^k \in \Gamma_{i,j}}
      \gamma_{i,j}^k ,\quad 
      \chi(\Gamma_{i}) := \bigcup_{j=1}^{N_{{\cal M}_i}}
      \chi(\Gamma_{i,j}),\quad
      \chi(\Gamma) := \bigcup_{i=1}^{N_p} \chi(\Gamma_{i}), 
    \end{array}
  \end{equation}
  where $\chi(\Gamma_{i,j})$ is a subset of $\mathbb{R}^2$,
  so are $\chi(\Gamma_{i})$
  and the interface $\chi(\Gamma)$.
\end{ntn}
\revise{It is emphasized that
  each $\Gamma_{i,j}$ is a poset of \emph{oriented} Jordan curves
  whereas each $\chi(\Gamma_{i,j})=\partial{\cal M}_{i,j}$
  is an \emph{unoriented} point set in $\mathbb{R}^2$. 
  The poset $\Gamma$ 
  contains both the topology and geometry of all Yin sets.
}

\begin{figure}
  \centering
  \subfloat[Six Yin sets $\{{\cal M}_i: i=1,\ldots,6\}$ 
  and their interface graph
  $G_{\Gamma}=(V_{\Gamma}, E_{\Gamma}, \psi_{\Gamma})$
  ]{
    \includegraphics[width=0.46\linewidth]{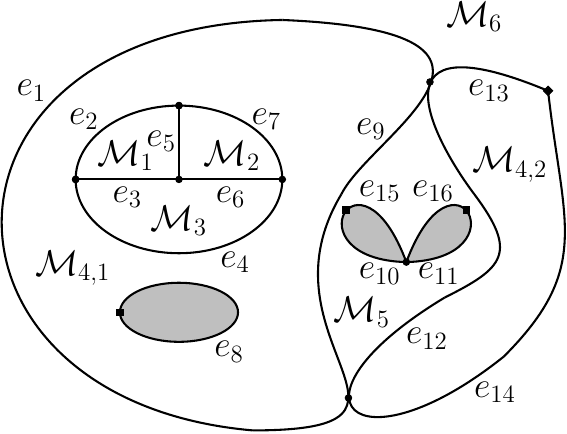}
  }
  \hfill
  \subfloat[
  The directed edges
  that constitute oriented boundary Jordan curves in $\Gamma$]{
    \includegraphics[width=0.48\linewidth]{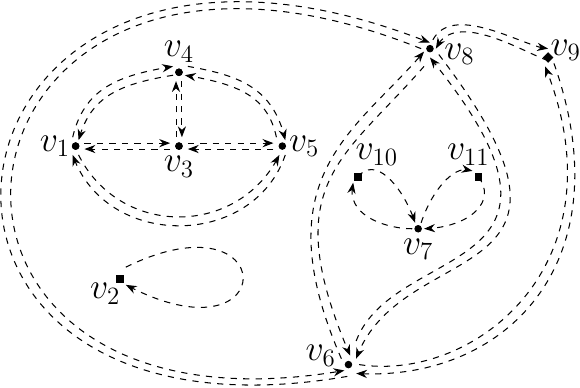}
  }

  \subfloat[Representing oriented boundary Jordan curves
  of the seven connected Yin sets. 
  The set $C_S$ of circuits and the set $T_S$ of trails
  facilitate spline fitting. 
  A smooth Jordan curve such as $\gamma_{4,1}^k$
  corresponds to a circuit in $C_S$. 
  In the last column,
  a trail $\mathbf{e}$
  corresponds to a smooth curve segment $\gamma$
  that is approximated by a not-a-knot spline;  
  $\gamma$ may or may not be closed.
  ]{\small
    \renewcommand{\arraystretch}{1.2}
    \input{{\TABDIR}figure4-1c}
  }
  \caption{The boundary representation
    of six pairwise disjoint Yin sets
    whose regularized union covers the plane
    except the shaded regions. 
    All Yin sets are connected except
    that ${\cal M}_4$ has two components:
    ${\cal M}_4={\cal M}_{4,1}\cup^{\perp\perp} {\cal M}_{4,2}$.
    The interface $\chi(\Gamma)$
    is represented by the graph
    \mbox{$G_{\Gamma} = (V_{\Gamma}, E_{\Gamma}, \psi_{\Gamma})$}
    in Definition \ref{def:topoGraph}
    with $E_{\Gamma}=\{e_i: i=1,\ldots,16\}$ shown in (a)
    and $V_{\Gamma}=\{v_i: i=1,\ldots,11\}$ \mbox{in (b)}.
    Solid dots, diamonds, and solid squares
    respectively denote junctions in Definition \ref{def:junctionPoint}, 
    non-smooth points in Definition \ref{def:kink},
    and basepoints 
    in Definition \ref{def:topoGraph}.
    In (c), 
    oriented boundary Jordan curves
    in Notation \ref{ntn:B-REP}
    and corresponding directed cycles in Notation \ref{ntn:cycleSet}
    are enumerated for each of the seven connected Yin sets.
    As an edge partition of $E_{\Gamma}$,
    the circuits in $C_S$ and the trails in $T_S$
    correspond to smooth closed curves
    and smooth curve segments,
    respectively approximated by periodic splines
    and not-a-knot splines in $S_{CT}$, 
    cf. Definition \ref{def:setOfFittedSplines}
    and Algorithm \ref{alg:partitionGraph}. 
  }
  \label{fig:multiphase-representation}
\end{figure}

\revise{The rule of thumb of our boundary representation of Yin sets
  is to \emph{separate their topology from their geometry}: 
  the former is represented by the incidence function
  of the interface graph in Definition \ref{def:topoGraph}
  and the cycle sets in Notation \ref{ntn:cycleSet}
  while the latter approximated by the cubic splines
  in Sec.~\ref{sec:geometryOfCommonBoundary}. 
As a major advantage of this design, 
 for a homeomorphic flow map
 only the \emph{geometry} of $\Gamma(t)$
 needs to be updated at each time step 
 while the topology of $\Gamma(t)$ can be determined
 from the initial condition $\Gamma(t_0)$ once and for all.
We illustrate key points of this design  
 by a comprehensive example 
 in Fig.~\ref{fig:multiphase-representation}.
}

First, 
 the topology of $\Gamma$ in the case of $N_p>2$
 is fundamentally more complicated than that of $N_p=2$, 
 because the common boundary of any two connected Yin sets
 might not be a Jordan curve, 
 due to the potential presence of T junctions
 such as $v_1,v_3,v_4,v_5$
 in Fig.~\ref{fig:multiphase-representation}(b).
Although X junctions such as $v_7$
 in Fig.~\ref{fig:multiphase-representation}(b)
 may also show up in two-phase flows,
 they tend to appear more frequently in three or more phases,
 cf. $v_{6},v_8$ in Fig.~\ref{fig:multiphase-representation}(b).
In both cases,
 the degree of a junction 
 can be any positive integer greater than two. 
These complications are handled
in Sec.~\ref{sec:topologyOfCommonBoundary}.

Second,
 boundary Jordan curves of adjacent Yin sets
 may have distinct geometric features.
In Fig.~\ref{fig:multiphase-representation},
 $\gamma_{4,1}^{1-}\cap\gamma_1^{+}=e_2$, 
 but $\gamma_{4,1}^{1-}$ is smooth
 while $\gamma_1^{+}$ is only ${\cal C}^0$
 due to the junctions $v_1,v_3,v_4$.
As shown in Fig.~\ref{fig:independent-mars},
 separate approximations of Jordan curves
 may lead to overlaps and/or vacuums of adjacent phases.
This problem can be solved by
 approximating each common boundary \emph{only once}
 with an appropriate spline type;
 see Sec.~\ref{sec:geometryOfCommonBoundary}.
For example,
 we can approximate $\gamma_{4,1}^{1-}$
 with a periodic cubic spline,
 cut the spline into three pieces at $v_1,v_4,v_5$,
 and reuse them in assembling other Jordan curves
 that share common boundaries with $\gamma_{4,1}^{1-}$.

Lastly,
 we combine topological and geometric data structures 
 into an approximation of $\Gamma$
 in Notation \ref{ntn:splineApproxOfGamma};
 see Fig.~\ref{fig:pipeline}
 and the last paragraph of Sec.~\ref{sec:extractingSinglePhase}
 for a summary of our design
 of the discrete boundary representation of multiple Yin sets.

\subsection{Representing the topology 
  of $\Gamma$}
\label{sec:topologyOfCommonBoundary}

At the center of representing
 the \emph{unoriented} point set $\chi(\Gamma)$ 
 and \emph{oriented} boundary curves of the Yin sets is 

\begin{defn}
  \label{def:graph}
  A \emph{graph} is a triple $G=(V,E,\psi)$
  where $V$ is a set of \emph{vertices},
  $E$ is a set of \emph{edges},
  and \mbox{$\psi: E\rightarrow V\times V$} is 
  the \emph{incidence function} given by \mbox{$\psi(e)=(v_s,v_t)$};
  the vertices $v_s$ and $v_t$
  are respectively the \emph{source} and \emph{target} of the edge $e$.
  $G$ is \emph{undirected}
  if we do not distinguish the source and target for any edge;
  $G$ is \emph{directed} if we do for all edges.
  The \emph{set of edges incident to $v\in V$} is 
  \begin{equation}
    \label{eq:adjacentEdgesAtv}
    E_v := \{e\in E: \psi(e)=(v,\cdot)  \text{ or } (\cdot,v)\}.
  \end{equation}
\end{defn}

An edge $e\in E$ is a \emph{self-loop}
 if $\psi(e)=(v,v)$ for some vertex $v$.
The \emph{degree or valence of a vertex}, written $\# E_v$, 
 is the number of edges incident to $v$,
 with each self-loop counted twice.
A \emph{subgraph} of $G=(V,E,\psi)$ 
 is a graph $G'=(V',E',\psi')$ such that
 $V'\subseteq V$, \mbox{$E'\subseteq E$}, and
 $\psi'=\left.\psi\right|_{E'}$. 
 
\begin{defn}[Types of subgraphs]
  \label{def:subgraphs}
  A \emph{walk} is a sequence of edges joining a sequence of vertices.
  A \emph{trail} is a walk where all edges are distinct.
  A \emph{circuit} is a non-empty trail
  where the first and last vertices coincide.
  A \emph{cycle} is a circuit
  where only the first and last vertices coincide.
\end{defn}

 \begin{defn}
   \label{def:planarGraph}
  A \emph{planar graph} is a graph $G=(V,E,\psi)$ satisfying
  \begin{enumerate}[label=(\alph*)]\itemsep0em
  \item each vertex in $V$ is a point in $\mathbb{R}^2$,
  \item each edge in $E$ is a curve $\gamma: [0,1] \to \mathbb{R}^2$
    whose endpoints are in $V$,
  \item two different edges/curves in $E$ do not intersect
    except at vertices in $V$, 
  \item the incidence function is given by
    $\psi(\gamma) := (\gamma(0), \gamma(1))$
    for any $\gamma\in E$.
  \end{enumerate}
\end{defn}

Any planar graph admits a dual graph, 
 \revise{which promptly yields the adjacency relation
 of any given pair of Yin sets \cite[Sec. 4.6]{diestel17:_graph_theor}.
This feature is helpful in coupling
 the proposed IT method with a main flow solver.
}

\begin{defn}
  \label{def:junctionPoint}
  A \emph{junction of the interface} $\chi(\Gamma)$
  is a point $p\in\chi(\Gamma)$
  such that, for any $\epsilon>0$,
  the intersection of $\chi(\Gamma)$
  with the $\epsilon$-open ball centered at $p$ 
  is \emph{not} homeomorphic to the interval $(0,1)$.
  The set of all junctions of $\chi(\Gamma)$
  is denoted by $J_{\Gamma}$. 
\end{defn} 

In particular,
 the T and Y junctions shown in Fig.~\ref{fig:independent-mars}
 are junctions of degree 3 
 and an X junction is a junction of degree 4.
Since we approximate $\chi(\Gamma)$ with cubic splines, 
 a curve is said to be \emph{smooth} 
 if it is ${\cal C}^4$.
If quintic splines were employed,
 it would be appropriate to define 
 a smooth curve as ${\cal C}^6$.
 
\begin{defn}
  \label{def:kink}
  A \emph{non-smooth point of the interface} $\chi(\Gamma)$
  is a point
  $p \in \left(\chi(\Gamma)
  \setminus J_{\Gamma}\right)$
  such that $\chi(\Gamma)$ is not smooth at $p$.
  The set of all non-smooth points of $\chi(\Gamma)$
  is denoted by $K_{\Gamma}$. 
\end{defn} 

\revise{A kink introduced in Sec.~\ref{sec:intro}
 is a special non-smooth point; 
 the latter is more general
 and facilitates the accuracy analysis
 of boundary representations 
 of arbitrary complex geometry.}

By Sec.~\ref{sec:yinsets},
 a curve segment and its corresponding open curve are different; 
 this is exploited in the representation of $\Gamma$. 
\begin{defn}[Interface graph]
  \label{def:topoGraph}
  The \emph{interface graph} of $N_p$ pairwise disjoint Yin sets
  is an undirected planar graph
  $G_{\Gamma}=(V_{\Gamma}, E_{\Gamma}, \psi_{\Gamma})$
  constructed as follows.
  \begin{enumerate}[label=(\alph*)]\itemsep0em
  \item Initialize $V_{\Gamma}\leftarrow J_{\Gamma}\cup K_{\Gamma}$
    and $E_{\Gamma}\leftarrow \emptyset$;
  \item Any curve $\gamma\in\Gamma_E:=\chi(\Gamma)
      \setminus (J_{\Gamma}\cup K_{\Gamma})$
      is one of the three types: 
      (i) a Jordan curve, 
      (ii) an open curve whose corresponding curve segment is not closed,
      or (iii) an open curve whose corresponding curve segment
      is Jordan.      
      \begin{itemize}
      \item For $\gamma$ of type (i),
        add $\gamma(\frac{1}{2})$ into $V_{\Gamma}$
        and add $\gamma$ as a self-loop into $E_{\Gamma}$.
      \item For $\gamma$ of type (ii),
        add into $E_{\Gamma}$ its corresponding curve segment; 
      \item For $\gamma$ of type (iii),
        add $\gamma$ as a self-loop into $E_{\Gamma}$
        if $\gamma(0)=\gamma(1)$ is a non-smooth point;
        otherwise
        add $\gamma(\frac{1}{2})$ into $V_{\Gamma}$
        and add into $E_{\Gamma}$ the two curve segments
        $\gamma([0,\frac{1}{2}])$
        and $\gamma([\frac{1}{2},1])$. 
      \end{itemize}
  \item Deduce the incidence function $\psi_{\Gamma}$ of $G_{\Gamma}$
    from (d) of Definition \ref{def:planarGraph}.
  \end{enumerate}
  The point $v_{\gamma}:=\gamma(\frac{1}{2})$ of types (i, iii) 
  is the \emph{basepoint of the Jordan curve} $\gamma$. 
\end{defn}




See Fig.~\ref{fig:multiphase-representation}
 for an illustration of the construction steps
 in Definition \ref{def:topoGraph}.
For $\gamma$ of type (i) in (b), 
 we add into $V_{\Gamma}$
 the basepoint $v_{\gamma}$ of the Jordan curve 
 so that $\psi_{\Gamma}(\gamma)=(v_{\gamma},v_{\gamma})$. 
For type (iii) where multiple Jordan curves
 intersect at a single junction, 
 it is necessary to add the basepoint of each Jordan curve
 into $V_{\Gamma}$;
 otherwise it would be difficult to enforce
 the smoothness of a trail that spans
 multiple Jordan curves,
 cf. the trail $(e_{15},e_{10},e_{11},e_{16})$
 in Fig.~\ref{fig:multiphase-representation}.
 
By Theorem \ref{thm:uniqueCases}, 
 each boundary Jordan curve $\gamma_{i,j}^k$ of a Yin set
 induces a directed cycle $C_{i,j}^k$,
 whose constituting edges
 come from the interface graph $G_{\Gamma}$
 and inherit their orientations
 from that of $\gamma_{i,j}^k$.
 
\begin{ntn}
  \label{ntn:cycleSet}
  Denote by $C_{i,j}^k$
  the \emph{$(i,j,k)$th directed cycle} 
  of the oriented boundary Jordan curve $\gamma_{i,j}^k$.
  Analogous to Notation \ref{ntn:B-REP},
  the \emph{(directed) cycle sets} of $\Gamma$ are denoted by
  $C_{i,j}:=\{C_{i,j}^k\}$,
  $C_{i}:=\{C_{i,j}: j=1,\ldots, N_{{\cal M}_i}\}$,
  and \mbox{$C:=\{C_{i}: i=1,\ldots,N_p\}$}.
\end{ntn}

\revise{
See Fig.~\ref{fig:multiphase-representation}(b) 
 and the third column of Fig.~\ref{fig:multiphase-representation}(c)
 for all directed cycles of the Yin sets
 in Fig.~\ref{fig:multiphase-representation}(a), 
 where the shaded region 
 does not belong to any of these Yin sets. 
Consequently, 
 the counter-clockwise self-loop with basepoint $v_2$
 is not in any directed cycle,
 nor is any of the two counter-clockwise cycles
 adjacent at the junction $v_7$.
}

\subsection{Approximating the geometry of
  $\chi(\Gamma)$}
\label{sec:geometryOfCommonBoundary}

The interface topology is captured
 in $\psi_{\Gamma}$ and $C$
 while its geometry in $E_{\Gamma}$. 
 
\begin{defn}
  \label{def:splineBoundarySegs}
  The \emph{spline edge set} $S_E$ of an interface $\chi(\Gamma)$
  is a set of splines that approximate curves in $E_{\Gamma}$. 
  The \emph{set of marker sequences} of $\chi(\Gamma)$ is 
  \begin{equation}
    \label{eq:setOfMarkSequences}
    E_X:=\{(v_i,X_{1},\ldots, X_{N_{\gamma}-1},v_j):
    \gamma\in E_{\Gamma},\ 
    \psi_{\Gamma}(\gamma)=(v_i,v_j)\},
  \end{equation}
  where $X_{1},\ldots, X_{N_{\gamma}-1}$
  are points on $\gamma$ selected as its \emph{interior markers}.
\end{defn}

Besides the one-to-one correspondence
 between elements in $S_E$ and $E_{\Gamma}$,
 any two corresponding elements in $S_E$ and $E_{\Gamma}$
 are made homeomorphic by using a sufficient number of interior markers
 for the sequence in $E_X$.
In Subsections \ref{sec:periodic-cubic-splines-with}
 and \ref{sec:NAK-cubic-splines-with}, 
 we discuss two types of splines employed
 for generating $S_E$.

\subsubsection{Cubic spline functions}
\label{sec:cubic-splines}

The \emph{arc length} of a continuously differentiable curve
 $\gamma:[0,1]\rightarrow\mathbb{R}^2$
 is a continuous function \mbox{$s_{\gamma}: [0,1]\rightarrow[0,L_{\gamma}]$}
 where $s_{\gamma}(l) := \int_{0}^{l}
 \sqrt{x'_{\gamma}(\tau)^2+y'_{\gamma}(\tau)^2} \dif \tau$
 and  $L_{\gamma}$ is the total length of $\gamma$.
Reparametrize $\gamma$
 as $[0,L_{\gamma}]\rightarrow \mathbb{R}^2$, 
 consider $\gamma$ as two coordinate functions
 $x_{\gamma},y_{\gamma}: [0,L_{\gamma}]\rightarrow \mathbb{R}$
 with the same domain, 
 and we approximate $x_{\gamma}$ and $y_{\gamma}$ separately via 

\begin{defn}[Space of spline functions]
  \label{def:splineFuncs}
  Over a \emph{breakpoint sequence} \mbox{$X_b:=(l_i)_{i=0}^{N}$}
  that partitions $[a,b]$ as
  \begin{equation}
    \label{eq:subdivision}
    a=l_0< l_1 < \cdots < l_{N} =b; \quad N\ge 3, 
  \end{equation}
  the space of \emph{spline functions of degree $m\in \mathbb{N}$
    and smoothness class $j\in\mathbb{N}$} is
  \begin{equation}
    \label{eq:splineFunctions}
    \mathbb{S}_m^j(X_b) := \bigl\{
    s\in {\cal C}^j[a,b]: \ \forall i=0,\ldots,N-1, \ 
    \left. s \right |_{[l_i, l_{i+1}]}\in \mathbb{P}_m
    \bigr\},
  \end{equation}
  where $\mathbb{P}_m$ is the space of polynomials 
  whose degrees are no more than $m$
  and each $l_i$ is called a \emph{breakpoint} of $s$. 
\end{defn}

\begin{figure}
  \centering
  \subfloat[A spline function in $\mathbb{S}_3^2(X_b)$; 
  see (\ref{eq:splineFunctions})]{
    \includegraphics[width=0.47\linewidth]{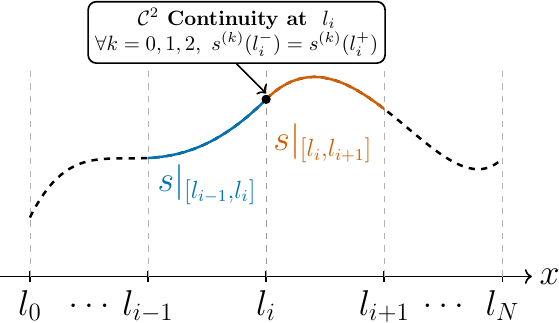}
  }
  \hfill
  \subfloat[The cumulative chordal length in
  (\ref{eq:cumulativeChordalLength})]{
    \includegraphics[width=0.45\linewidth]{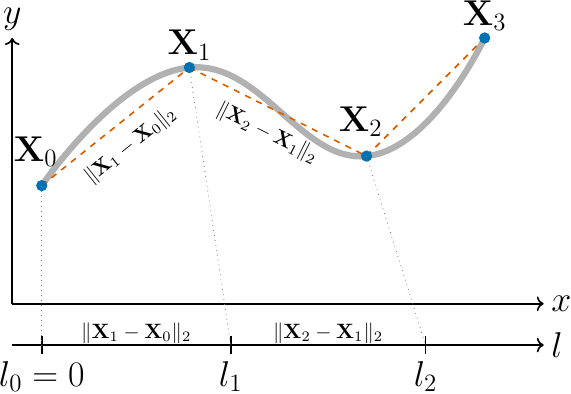}
  }
  \caption{\revise{Spline curve fitting
      based on the cumulative chordal lengths.
      In subplot (b), the thick solid curve
      represents an interface to be approximated
      and the dots interface markers.}}
  \label{fig:spline}
\end{figure}

$\mathbb{S}_3^2$ is probably the most popular class of spline functions. 
As implied by (\ref{eq:splineFunctions})
 \revise{and shown in Fig.~\ref{fig:spline}(a)},
the restriction of any $s\in \mathbb{S}_3^2$ on a subinterval 
is a cubic polynomial 
and thus $4N$ coefficients need to be determined.
In interpolating a function $f:[a,b]\rightarrow\mathbb{R}$
 by $s\in\mathbb{S}_3^2$,
 the number of equations given
 by the interpolation conditions at all breakpoints
 and by the continuity requirements at interior breakpoints
 is respectively $N+1$ and $3(N-1)$,
 leading to a total of $4N-2$ equations.
The last two equations come from 

\begin{defn}[Types of cubic spline functions]
  \label{def:splineTypes}
  A \emph{periodic cubic spline function}
  is a spline function $s\in \mathbb{S}_3^2$ 
  satisfying $s(0)=s(1)$, $s'(0)=s'(1)$, and $s''(0)=s''(1)$.
  A \emph{not-a-knot cubic spline function}
  \cite{deBoor1966method}
  is a spline function $s\in \mathbb{S}_3^2$
  such that $s'''(l)$ exists at $l=l_1$ and $l=l_{N-1}$,
  cf. Definition \ref{def:splineFuncs}. 
\end{defn}

A \emph{spline} is a pair of spline functions, 
 each of which is interpreted as a coordinate function
 of the approximated curve. 
A smooth closed curve 
 is approximated by a periodic cubic spline, 
 although it is not necessarily Jordan
 due to potential self-intersections.
In contrast, a smooth curve segment 
 can be approximated by a not-a-knot cubic spline.
In both cases,
 the \emph{cumulative chordal length}
 is a discrete counterpart
 of the arc length of $\gamma\in\Gamma_E$
 from a sequence of distinct markers
 $(\mathbf{X}_i)_{i=0}^{N}$ on $\gamma$: 
\begin{equation}
  \label{eq:cumulativeChordalLength}
  l_0=0; \quad
  \forall i=1,\ldots,N, \ \ 
  l_i=l_{i-1} + \|\mathbf{X}_{i}-\mathbf{X}_{i-1}\|_2, 
\end{equation}
where $\|\cdot\|_2$ denotes the Euclidean-norm;
\revise{see Fig. \ref{fig:spline}(b)}. 
Having the same breakpoint sequence $X_b = (l_i)_{i=0}^N$, 
 the two coordinate spline functions
 are determined \emph{separately}
 and then combined as the interpolatory spline
 of $\gamma$.
 

\subsubsection{Periodic cubic splines}
\label{sec:periodic-cubic-splines-with}

The following conclusions
 on the unique existence and interpolation accuracy 
 of periodic cubic splines
 have been proved in \cite[Sec.~2.2.1]{hu2025:_arms}.

\begin{thm}
  \label{thm:periodicSplineErrorEstimates}
  A periodic spline $p\in \mathbb{S}^2_3(X_b)$
  that interpolates a periodic function
  $f\in {\cal C}^2([a,b])\cap {\cal C}^4([a,b]\setminus X_b)$
  at $X_b$ 
  can be uniquely determined and satisfies 
  \begin{equation}
    \label{eq:periodicSplineErrorEstimates}
    \forall l \in[a,b], \ \forall j= 0, 1, 2, \
    \left|p^{(j)}(l)-f^{(j)}(l)\right| \le c_j h^{4-j}
    \max_{\xi\in[a,b]\setminus X_b}\left|f^{(4)}(\xi)\right|,
  \end{equation}
  where $h:=\max_{i=1}^{N}|l_{i}-l_{i-1}|$ 
  and the constants are 
  $c_0=\frac{1}{16}$ and $c_1=c_2=\frac{1}{2}$.
\end{thm}

\revise{
  A periodic spline is homeomorphic to a cycle
  while a not-a-knot spline to an interval.
  This fundamental difference
  leads to substantial difference
  in the analysis of their accuracy and stability.
}




 
\subsubsection{Not-a-knot cubic splines} 
\label{sec:NAK-cubic-splines-with}

\revise{
  In this section, we analyze not-a-knot splines 
  to prove their unique existence and their fourth-order accuracy
  in interpolating smooth curves.
 }

\begin{lem}
  \label{lem:3moments}
  Let $p \in \mathbb{S}_3^2(X_b)$ be a not-a-knot spline
  interpolating a function $f$ at the nodes $X_b$.
  Then the second derivatives $M_i = p''(l_i)$ satisfy:
  \begin{eqnarray}
    \label{eq:3Moments}
    \forall i= 1,\ldots, N-1, \quad
    \mu_i M_{i-1} + 2M_i + \lambda_i M_{i+1}
    = 6 f[l_{i-1},l_i,l_{i+1}],  \\
    \label{eq:3Moments2}
    \lambda_1 M_0 - M_1 + \mu_1 M_2 = 0, \quad 
    \lambda_{N-1} M_{N-2} - M_{N-1} + \mu_{N-1} M_{N} = 0, 
  \end{eqnarray}
  where
  \begin{equation}
    \label{eq:mu/lambda}
    \begin{array}{l}
      \forall i=1,\ldots,N-1, \qquad 
      \mu_i := \frac{l_i-l_{i-1}}{l_{i+1}-l_{i-1}},\quad
      \lambda_i := \frac{l_{i+1}-l_{i}}{l_{i+1}-l_{i-1}},
    \end{array}
  \end{equation}
  and the \emph{divided difference} is recursively given by
  \begin{equation}
    \label{eq:dividedDiff}
    \begin{array}{l}
      f[l]:=f(l); \quad 
      f[l_{0},l_1,\ldots, l_{j}]:=
      \frac{f[l_1,\ldots, l_{j}]-f[l_0,\ldots, l_{j-1}]}{l_j-l_0}. 
    \end{array}
  \end{equation}
\end{lem}
\begin{proof}
  Taylor expansion of $p(l)$ at $l_i$ yields
  \begin{equation}
    \label{eq:proofLem3Moments1}
    \begin{array}{l}
      p(l) = f_i + p'(l_i)(l-l_i)
      + \frac{M_i}{2}(l-l_i)^2 + \frac{p'''(l_i)}{6}(l-l_i)^3.
    \end{array}
  \end{equation}
  Differentiate (\ref{eq:proofLem3Moments1}) twice, 
  set $l = l_{i+1}$,
  and we have $p'''(l_i) = \frac{M_{i+1}-M_i}{l_{i+1} - l_i}$,
  the substitution of which back into (\ref{eq:proofLem3Moments1}) yields
  \begin{equation}
    \label{eq:proofLem3Moments3}
    \begin{array}{l}
      p'(l_i) = f[l_i,l_{i+1}] - \frac{1}{6}(M_{i+1}+2M_i)(l_{i+1}-l_i).
    \end{array}
  \end{equation}
  Similarly,
  for $l\in [l_{i-1},l_i)$,
  differentiate (\ref{eq:proofLem3Moments1}) twice,
  set $l = l_{i-1}$, and we have
  $p'''(l_i) = \frac{M_{i-1}-M_i}{l_{i-1} - l_i}$ and 
  \begin{equation}
    \label{eq:proofLem3Moments4.5}
    \begin{array}{l}
      p'(l_i) = f[l_{i-1}, l_i] - \frac{1}{6}(M_{i-1}+2M_i)(l_{i-1} - l_i).
    \end{array}
  \end{equation}
  Then (\ref{eq:3Moments}) follows from 
  subtracting (\ref{eq:proofLem3Moments3}) 
  from (\ref{eq:proofLem3Moments4.5})
  and applying (\ref{eq:dividedDiff}).

  By Definition \ref{def:splineTypes}, 
  the not-a-knot boundary condition requires continuity 
  of $p'''(l)$ at both $l_1$ and $l_{N-1}$, i.e.,
  $\forall i\in \{1,N-1\}, 
    \frac{M_{i+1}-M_i}{l_{i+1} - l_i} 
    = \frac{M_{i-1}-M_i}{l_{i-1} - l_i}$,
  which, when multiplied by
  $\frac{(l_{i+1} - l_i)(l_{i} - l_{i-1})}{l_{i+1} - l_{i-1}}$,
  yields (\ref{eq:3Moments2}). \qed
\end{proof}

\revise{
\begin{lem}
  \label{lem:uniquenessNAK}
  For any $f \in \mathcal{C}^2([a,b])$
  and a breakpoint sequence $X_b$ over $[a,b]$, 
  there exists a unique not-a-knot spline
  $p \in \mathbb{S}^2_3(X_b)$ that interpolates $f$.
\end{lem}
\begin{proof}
By Lemma~\ref{lem:3moments}, we have 
\begin{equation}
\label{eq:miLinearSystem2}
\begin{array}{l}
  A\mathbf{M}
  :=\begin{bmatrix}
      \lambda_1 & -1    & \mu_1     &           &               &               &               \\
      \mu_1     & 2     & \lambda_1 &           &               &               &               \\
                & \mu_2 & 2         & \lambda_2 &               &               &               \\
                &       & \ddots    & \ddots    & \ddots        &               &               \\
                &       &           & \mu_{N-2} & 2             & \lambda_{N-2} &               \\
                &       &           &           & \mu_{N-1}     & 2             & \lambda_{N-1} \\
                &       &           &           & \lambda_{N-1} & -1            & \mu_{N-1}
    \end{bmatrix}
  \begin{bmatrix}
    M_0 \\ M_1 \\ M_2 \\ \vdots \\M_{N-2} \\ M_{N-1} \\ M_N
  \end{bmatrix}
  =\mathbf{b},
\end{array}
\end{equation}
where $\mu_i, \lambda_i$ are defined in (\ref{eq:mu/lambda}),
$b_0 = 0$, $b_N = 0$, and 
\begin{displaymath}
  \begin{array}{l}
    \forall i = 1,\ldots,N-1,\quad 
    b_i = 6f[l_{i-1},l_i,l_{i+1}].
  \end{array}
\end{displaymath}

Define $T_1,T_2,T_3\in \mathbb{R}^{(N + 1)\times(N + 1)}$ as
$T_3=\text{diag}(1, \frac{1}{1 - \mu_1}, 1, \ldots, 1,
\frac{1}{1 - \lambda_{N-1}}, 1)^{\top}$, 
\begin{equation}
  \label{eq:TDef}
  T_1 := 
  \begin{bmatrix}
    1 & 1 &        &   &   \\
    0  & 1 & 0       &   &   \\
      & \ddots  & \ddots & \ddots  &   \\
      &   &      0 & 1 & 0 \\
      &   &        & 1 & 1 \\
  \end{bmatrix},
  \ 
  T_2 := 
  \begin{bmatrix}
    1      & 0 &        &        &        &   &                \\
    -\mu_1 & 1 & 0      &        &        &   &                \\
           & 0 & 1      & 0      &        &   &                \\
           &   & \ddots & \ddots & \ddots &   &                \\
           &   &        & 0      & 1      & 0 &                \\
           &   &        &        & 0      & 1 & -\lambda_{N-1} \\
           &   &        &        &        & 0 & 1
  \end{bmatrix}.
\end{equation}
Then $\mu_1, \lambda_1, \mu_{N-1}, \lambda_{N-1} 
 \in (0, 1)$, $\mu_i + \lambda_i = 1$, and $N \ge 3$ yield
\begin{equation}
  \label{eq:TransADef}
  T_3 T_2 T_1 A =
  \begin{bmatrix}
  1 & 1    & 1     &           &               &               &               \\
  & \tfrac{2 - \mu_1}{1 - \mu_1} & \tfrac{1 - 2\mu_1}{1 - \mu_1}  \\
  & \mu_2 & 2         & \lambda_2 &               &               &               \\
  &       & \ddots    & \ddots    & \ddots        &               &               \\
  &       &           & \mu_{N-2} & 2             & \lambda_{N-2} &               \\
  &       &           &           & \tfrac{1 - 2\lambda_{N-1}}{1 - \lambda_{N-1}} & \tfrac{2 - \lambda_{N-1}}{1 - \lambda_{N-1}}              & \\
  &       &           &           & 1 & 1            & 1    \\
  \end{bmatrix}
\end{equation}
and $T_3 T_2 T_1 \mathbf{b} = 
(b_1, b_1, b_2, \ldots, b_{N-2}, b_{N-1}, b_{N-1})^{\top}$. 
Consequently, this transformed linear system 
\begin{equation}
  \label{eq:transEqs}
  \begin{array}{l}
  T_3 T_2 T_1 A \mathbf{M} = T_3 T_2 T_1 \mathbf{b}
  \end{array}
\end{equation}
 can be written as
\begin{align}
  \label{eq:transMid}
  &\forall j = 2, \cdots, N-2, \quad
  \mu_j M_{j-1} + 2 M_j + \lambda_{j} M_{j+1} =
  6 f[l_{j-1}, l_j, l_{j+1}], \\
  \label{eq:transEnd1}
  &\begin{cases}
    \tfrac{2-\mu_1}{1 - \mu_1} M_1 + \tfrac{1 - 2\mu_1 }{1 - \mu_1 } M_2
  = 6 f[l_{0}, l_1, l_{2}], \\
    \tfrac{1 - 2\lambda_{N-1} }{1 - \lambda_{N-1} } M_{N-2} + \tfrac{2-\lambda_{N-1}}{1 - \lambda_{N-1}} M_{N-1} = 6 f[l_{N-2}, l_{N-1}, l_{N}],
  \end{cases}
  \\
  \label{eq:transEnd0}
  &\begin{cases}
  M_0 + M_1 + M_2 = 6 f[l_{0}, l_1, l_{2}], \\
  M_{N-2} + M_{N-1} + M_N = 6 f[l_{N-2}, l_{N-1}, l_{N}].
  \end{cases}
\end{align}
The equations corresponding to \eqref{eq:transMid} and \eqref{eq:transEnd1} 
 do not involve the variables $M_0$ and $M_N$; 
 we denote this linear subsystem as
\begin{equation}
  \label{eq:subEqus}
  \begin{array}{l}
    A_s \mathbf{M}_s = \mathbf{b}_s.  
  \end{array}
\end{equation}
It follows from \eqref{eq:mu/lambda}
that $A_s$ is strictly row diagonally dominant with
\begin{equation}
  \label{eq:AsRowdiagdominant}
  \begin{array}{ll}
   &\forall j = 2, \cdots, N-2, \quad \left| 2 \right| - \left|\mu_j\right| - \left|\lambda_j\right| = 1, \\
   &\left|\frac{2-\mu_1}{1 - \mu_1}\right|
  - \left| \frac{1 - 2\mu_1 }{1 - \mu_1 }\right|
  \ge \min\left( \frac{2-\mu_1}{1 - \mu_1} - \frac{1 - 2\mu_1 }{1 - \mu_1 },
  \frac{2-\mu_1}{1 - \mu_1} + \frac{1 - 2\mu_1 }{1 - \mu_1 }\right)
  = \min\left(\frac{1 + \mu_1 }{1 - \mu_1 }, 3\right)
  > 1,
  \\  
  &\left|\frac{2-\lambda_{N-1}}{1 - \lambda_{N-1}}\right|
  - \left| \frac{1 - 2\lambda_{N-1} }{1 - \lambda_{N-1} }\right|
  \ge \min\left(\frac{1 + \lambda_{N-1} }{1 - \lambda_{N-1} }, 3\right)
  > 1,
  \end{array}
\end{equation}
 which implies that the matrix $A_s$ is invertible 
 and the linear subsystem \eqref{eq:subEqus} 
 admits a unique solution $\mathbf{M}_s$.
Together with \eqref{eq:transEnd0}, it follows 
 that the full system \eqref{eq:transEqs} yields 
 a unique solution $\mathbf{M}$,   
 thereby uniquely determining the not-a-knot spline 
 $p \in \mathbb{S}^2_3(X_b)$. \qed
\end{proof}

A matrix $A\in \mathbb{R}^{n\times n}$ is
\emph{strictly diagonally dominant} (by rows) if
we have $|a_{ii}| > \sum_{j\ne i, j= 1}^n |a_{ij}|$
for each $i=1, 2, \ldots, n$.

\begin{thm}
  \label{thm:inftyNormEstiFromDiagDom}
  Suppose that $A$ is strictly diagonally dominant by rows
  and write $\alpha:= \min_i (|a_{ii}|-\sum_{j\ne i} |a_{ij}|)$.
  Then $\left\|A^{-1}\right\|_{\infty} \le \frac{1}{\alpha}$.
\end{thm}
\begin{proof}
  The proof is due to Varah \cite{varah1975}
  and we rephrase it here for completeness.

  Write $A\mathbf{x}=\mathbf{y}$
  and we have
  $\|A^{-1}\|_{\infty}
      = \sup_{\mathbf{y}\ne 0} \frac{\|A^{-1}\mathbf{y}\|_{\infty}}
      {\|\mathbf{y}\|_{\infty}}
      = \sup_{\mathbf{x}\ne 0} \frac{\|\mathbf{x}\|_{\infty}}
      {\|\mathbf{A\mathbf{x}}\|_{\infty}}$.
  Thus the reciprocal of $\|A^{-1}\|_{\infty}$
  is the minimum contraction
  that $A$ applies to any vector, i.e.,
  $\|A^{-1}\|_{\infty}^{-1}
  = \inf_{\mathbf{x}\ne 0} \frac{\|A\mathbf{x}\|_{\infty}}
      {\|\mathbf{\mathbf{x}}\|_{\infty}}$.
  Hence it suffices to show that 
  $\alpha \|\mathbf{x}\|_{\infty}\le \|A\mathbf{x}\|_{\infty}$
  holds for an arbitrary $\mathbf{x}$.
  Suppose $|x_k|=\|\mathbf{x}\|_{\infty}$. 
  Then the definition of $\alpha$ gives 
  $0<\alpha\le |a_{kk}|-\sum_{j\ne k} |a_{kj}|$
  and 
  \begin{displaymath}
    \begin{array}{rl}
      0<\alpha |x_k|
      &\le |a_{kk}||x_k|-\sum_{j\ne k} |a_{kj}||x_j|
      \le |a_{kk}x_k|-\sum_{j\ne k} |a_{kj}x_j|
      \\
      &\le \sum_{j} |a_{kj}||x_j|
      \le \max_k |\sum_{j} a_{kj}x_j| = \|A\mathbf{x}\|_{\infty}.\qed
    \end{array}
  \end{displaymath}
\end{proof}

\begin{lem}
  \label{lem:2ndDerivativeBound}
  For any given $f \in \mathcal{C}^2([a,b])$, 
  the not-a-knot spline $p \in \mathbb{S}^2_3(X_b)$
  that uniquely interpolates $f$ satisfies 
  \begin{equation}
    \label{eq:2ndDerivativeBound}
    \begin{array}{l}
      \forall l\in [a,b], \qquad
      |p''(l)| < 9 \max_{\xi\in [a,b]}\left|f''(\xi)\right|.
    \end{array}
  \end{equation}
\end{lem}
\begin{proof}
  Since $p''$ is piecewise linear on $[l_i, l_{i+1}]$,
  its maximum absolute value must occur at some breakpoint $l_{j}$. 
  By \eqref{eq:TransADef}, \eqref{eq:AsRowdiagdominant} 
   and Theorem~\ref{thm:inftyNormEstiFromDiagDom},
   the linear subsystem 
   $A_s \mathbf{M}_s = \mathbf{b}_s$ 
   in \eqref{eq:subEqus}
   satisfies $\|A_s^{-1}\|_\infty \le 1$ and 
  \begin{equation*}
    \begin{array}{l}
      \|\mathbf{b}_s\|_\infty = \|\mathbf{b}\|_\infty 
      = \max\nolimits_{1 \le j \le N-1} |6 f[l_{j-1}, l_j, l_{j+1}]|
      \le 3\max\nolimits_{\xi\in [a,b]}\left|f''(\xi)\right|.
    \end{array}
  \end{equation*} 
  Consequently,
  \begin{equation*}
    \begin{array}{ll}
      \max\nolimits_{1 \le j \le N-1} |M_j| 
      &= \|\mathbf{M}_s\|_\infty 
      = \|A_s^{-1} \mathbf{b}_s\|_\infty  
      \le \|A_s^{-1}\|_\infty \| \mathbf{b}_s\|_\infty 
      \le  \| \mathbf{b}\|_\infty \\
      &\le 3\max\nolimits_{\xi\in [a,b]}\left|f''(\xi)\right|.
    \end{array}
  \end{equation*} 
  Finally, by combining \eqref{eq:transEnd1} and \eqref{eq:transEnd0}
   to eliminate $M_1$, we obtain
  \begin{align*}
                      &
      (2 - \mu_1)(M_0 + M_1 + M_2) -
      (1 - \mu_1)\left(\tfrac{2-\mu_1}{1 - \mu_1} M_1 + \tfrac{1 - 2\mu_1 }{1 - \mu_1 } M_2\right)
      = 6f[l_{0}, l_1, l_{2}]  \notag     \\
      \implies &
      (2 - \mu_1)M_0 + (1 + \mu_1) M_2 = 6 f[l_{0}, l_1, l_{2}]  \notag  \\
      \implies &
      |M_0| \le \tfrac{1}{2 - \mu_1} |6f[l_{0}, l_1, l_{2}]|
      + \tfrac{1 + \mu_1}{2 - \mu_1} |M_2|
      \le \tfrac{2 + \mu_1}{2 - \mu_1} \| \mathbf{b}\|_\infty 
      < 3 \| \mathbf{b}\|_\infty \label{eq:M0Bound}  \\
      \implies &
      |M_0| < 9 \max\nolimits_{\xi\in [a,b]}\left|f''(\xi)\right|,
      \notag
  \end{align*}
   where the third inequality follows from
   $\mu_1 = \frac{l_1 - l_0}{l_2 - l_0} < 1$.
  The conclusion for $M_N$ follows from similar arguments. \qed
\end{proof}
 
\begin{thm}
  \label{thm:notaKnotSplineErrorEstimates}
  For any given $f \in \mathcal{C}^2([a,b])
  \cap \mathcal{C}^4([a,b]\setminus X_b)$, 
  the not-a-knot spline $p \in \mathbb{S}^2_3(X_b)$
  that interpolates $f$ satisfies 
  \begin{equation}
    \label{eq:notaKnotSplineErrorEstimates}
    \forall l\in[a,b], \ \forall j= 0, 1, 2, \
    \left|p^{(j)}(l)-f^{(j)}(l)\right| < c_j h^{4-j}
    \max_{\xi\in[a,b]\setminus X_b}\left|f^{(4)}(\xi)\right|,
  \end{equation}
  where the constants are given 
  by $c_0=\frac{5}{32}$, 
  $c_1=c_2=\frac{5}{4}$.
\end{thm}
\begin{proof}
  By Lemma~\ref{lem:uniquenessNAK}, 
  the not-a-knot spline exists and is unique.

  For $j=2$, 
  we interpolate $f''(l)$ with some $\tilde{p}\in \mathbb{S}_1^0$
  and integrate $\tilde{p}$ twice to get $\hat{p}\in \mathbb{S}^2_3$
  so that $\hat{p}''$ interpolates $f''$ over $X_b$. 
  The Cauchy remainder theorem of polynomial interpolation yields 
  \begin{displaymath}
    \begin{array}{l}
      \exists  \xi_i \in (l_i,l_{i+1})  
      \text{ s.t. } \forall l\in [l_i,l_{i+1}],\, 
      |f''(l)-\tilde{p}(l)|\le \frac{1}{2}\left|f^{(4)}(\xi_i)\right|
      |(l-l_i)(l-l_{i+1})|
    \end{array}
  \end{displaymath}
  and hence we have
  \begin{equation}
    \label{eq:cubicSplineErrorEstimatesProof1}
    \begin{array}{l}
      |f''(l)-\hat{p}''(l)|\le \frac{h^2}{8}\max_{l\in[a,b]\setminus X_b}|f^{(4)}(l)|.
    \end{array}
  \end{equation}
  
  It follows from $\hat{p}(l)\in \mathbb{S}_3^2$ that
  $p(l) - \hat{p}(l)$ interpolates $f(l) - \hat{p}(l)$.
  Then Lemma \ref{lem:2ndDerivativeBound} yields
  $|p''(l)-\hat{p}''(l)| < 9
  \max_{\xi\in[a,b]}|f''(\xi)-\hat{p}''(\xi)|$
  for any $l \in[a,b]$, 
  which, together with the triangular inequality, gives
  \begin{equation}
    \label{eq:cubicSplineErrorEstimatesProof2}
    \begin{array}{rl}
      \forall l \in [a, b],\quad
      |f''(l)-p''(l)| 
      &< 10 \max_{\xi\in[a,b]}|f''(\xi)-\hat{p}''(\xi)| \\
      &\le \frac{5}{4}
      h^2\max_{\xi\in[a,b]\setminus X_b}\left|f^{(4)}(\xi)\right|, 
    \end{array}
  \end{equation}
  where the second step follows from (\ref{eq:cubicSplineErrorEstimatesProof1}).

  Finally, for $j=1$,
  the interpolating conditions give $f(l) - p(l) = 0$
  for $l = l_i, l_{i+1}$ 
  and Rolle's theorem implies
  $f'(\xi_i)-p'(\xi_i)=0$ for some $\xi_i \in (l_i, l_{i+1})$.
  Then the second fundamental theorem of calculus
  yields
  \begin{displaymath}
    \begin{array}{rl}
      \forall l\in [l_i, l_{i+1}],\quad
      f'(l) - p'(l) = 
        \int_{\xi_i}^l\left(f''(t)-p''(t)\right)\dif t,
    \end{array} 
  \end{displaymath}
  which, together with the integral mean value theorem
  and (\ref{eq:cubicSplineErrorEstimatesProof2}), gives
  \begin{displaymath}
    \begin{array}{rl}
      \forall l \in [a, b],\quad
      |f'(l) - p'(l)|_{l\in[l_i, l_{i+1}]} 
       &= \left| l - \xi_i \right|
      \left|f''(\eta_i) - p''(\eta_i)\right| \\
      &< \frac{5}{4}
      h^3\max_{\xi\in[a,b]\setminus X_b}\left|f^{(4)}(\xi)\right|.
    \end{array}
  \end{displaymath}
  
  For $j=0$, the interpolation of $f(l) - p(l)$
  with some $\bar{p}\in \mathbb{S}_1^0$ dictates
  $\bar{p}(l) \equiv 0$ for any $l \in [a, b]$. 
  Hence, for any $l \in [l_i, l_{i+1}]$, we have 
  \begin{displaymath}
    \begin{array}{rl}
      |f(l) - p(l)|
      = |f(l) - p(l) -\bar{p}|
      &\le \frac{1}{8}(l_{i+1}-l_i)^2\max_{\xi\in(l_i,l_{i+1})}|f''(\xi)-p''(\xi)| \\
        & < \frac{5}{32} h^4\max_{\xi\in[a,b]\setminus X_b}|f^{(4)}(\xi)|,
    \end{array}
  \end{displaymath}
  where the first inequality follows from
  the Cauchy remainder theorem 
  and the second 
  from (\ref{eq:cubicSplineErrorEstimatesProof2}).  \qed
\end{proof}

\subsubsection{The execellent conditioning
  of cubic periodic and not-a-knot splines}
\label{sec:conditioning}

In this subsection,
we prove that upper bounds of the condition number for fitting cubic
periodic and not-a-knot splines are 3 and 54, respectively. 

\begin{thm}
  \label{thm:periodicSplineConditionNumber}
  The max-norm-based condition number
  of fitting the unique periodic cubic spline 
  is bounded by 3.
\end{thm}
\begin{proof}
  By \cite[Lemma 2.7]{hu2025:_arms},
  the linear system for determining $M_i=p''(l_i)$'s
  is
  \begin{equation}
    \label{eq:miLinearSystem2Periodic}
    \begin{array}{l}
      A\mathbf{M}
      :=\begin{bmatrix}
        2 & \lambda_0    &      &           &               &           \mu_0          \\
        \mu_1     & 2     & \lambda_1 &           &               &                \\
        & \mu_2 & 2         & \lambda_2 &               &               \\
        &       & \ddots    & \ddots    & \ddots        &               \\
        &       &           & \mu_{N-2} & 2             & \lambda_{N-2} \\
        \lambda_{N-1}&       &           &           & \mu_{N-1}     & 2            
      \end{bmatrix}
      \begin{bmatrix}
        M_0 \\ M_1 \\ M_2 \\ \vdots \\M_{N-2} \\ M_{N-1} 
      \end{bmatrix}
      =\mathbf{b},
    \end{array}
  \end{equation}
  where $\mu_i, \lambda_i$ are defined in (\ref{eq:mu/lambda})
  for 
  $i=0,1,\ldots,N-1$ with $l_{-1} = l_{N-1}$.
  Hence $\|A\|_{\infty}=3$.
  By Theorem \ref{thm:inftyNormEstiFromDiagDom}, 
  $\alpha=1$ and
  $\|A^{-1}\|_{\infty}\le 1$.
  \qed
\end{proof}

For not-a-knot splines,
 the matrix in the above linear system
 is not diagonally dominant. 
So transformations are needed
 to decouple the first and the last lines 
 from other lines
 that are strictly diagonally dominant.
 
\begin{lem}
  \label{lem:inverseA-inftyNorm}
  The matrix $A$ in \eqref{eq:miLinearSystem2} 
  with $N \ge 5$
  satisfies $\|A^{-1}\|_{\infty} < 18$. 
\end{lem}
\begin{proof}
  Define an elementary transformation $T_4\in\mathbb{R}^{(N + 1) \times (N + 1)}$ as
  \begin{equation}
    \label{eq:T4Def}
    T_4 :=
    \begin{bmatrix}
      1 & 0 &          &        &          &   &   \\
      0 & 1 & -\beta_1                             \\
        & 0 & 1        & 0      &          &   &   \\
        &   & \ddots   & \ddots & \ddots   &   &   \\
        &   &          & 0      & 1        & 0 &   \\
        &   &          &        & -\beta_2 & 1 & 0 \\
        &   &          &        &          & 0 & 1 \\
    \end{bmatrix},
  \end{equation}
   where $\beta_1 := \tfrac{1 - 2\mu_1}{2(2 - \mu_1)}$ and 
   $\beta_2 := \tfrac{1 - 2\lambda_{N-1}}{2(2 - \lambda_{N-1})}$ 
   satisfy
   $\beta_1, \beta_2 \in \left(-\frac{1}{2}, \frac{1}{4}\right)$.
   Then
   \begin{equation*}
     T_2 T_1 A T_4 =
     \begin{bmatrix}
       1 & 1           & 1 -\beta_1            &           &                               &                     &   \\
        & {2 - \mu_1} & \tfrac{1 - 2\mu_1}{2}                                                                       \\
        & \mu_2       & 2 -\beta_1 \mu_2      & \lambda_2 &                               &                     &   \\
        &             & \ddots                & \ddots    & \ddots                        &                     &   \\
        &             &                       & \mu_{N-2} & 2  - \beta_2 \lambda_{N-2}    & \lambda_{N-2}       &   \\
        &             &                       &           & \tfrac{1 - 2\lambda_{N-1}}{2} & {2 - \lambda_{N-1}} &   \\
        &             &                       &           & 1  - \beta_2                  & 1                   & 1 \\
    \end{bmatrix}.
  \end{equation*}
  Consequently, the linear system 
   $T_2 T_1 A T_4 \mathbf{x} = \mathbf{y}$ can be expanded as
  \begin{align}
    \label{eq:transATMid}
    & \forall j = 3, \cdots, N-3, \quad
    \mu_j x_{j-1} + 2 x_j + \lambda_{j} x_{j+1} =
    y_j,
    \\ \label{eq:transATEnd2}
    &
    \begin{cases}
      \mu_2 x_1 + (2 - \beta_1 \mu_2) x_2 + \lambda_2 x_3
      = y_2, \\
      \mu_{N-2} x_{N-3} + (2 - \beta_2 \lambda_{N-2}) x_{N-2}
      + \lambda_{N-2} x_{N-1}
      = y_{N-2},
    \end{cases}
    \\
    \label{eq:transATEnd1}
    & (2 - \mu_1)x_1 + \tfrac{1-2\mu_1}{2}x_2 = y_1, \quad
    \tfrac{1-2\lambda_1}{2}x_{N-2} + (2 - \lambda_{N-1})x_{N-1}  = y_{N-1}
    \\
    \label{eq:transATEnd0}
    & \begin{cases}
        x_0 + x_1 + (1 - \beta_1)x_2 = y_0, \\
        (1 - \beta_2)x_{N-2} + x_{N-1} + x_N = y_N.
      \end{cases}
  \end{align}
  Note that \eqref{eq:transATMid}, \eqref{eq:transATEnd2}, 
   and \eqref{eq:transATEnd1} constitute a linear subsystem 
   $A_t \mathbf{x}_t = \mathbf{y}_t$. 
  Due to \eqref{eq:mu/lambda} and 
  $\beta_1, \beta_2 \in (-\frac{1}{2}, \frac{1}{4})$,
  $A_t$ is strictly diagonally dominant by rows: 
  \begin{equation*}
    \begin{array}{ll}
      & \forall j = 3, \cdots, N-3, \quad \left| 2 \right| - \left|\mu_j\right| - \left|\lambda_j\right| = 1 > \frac{1}{2}, \\
      & \left|{2-\mu_1}\right|
      - \left| \frac{1 - 2\mu_1 }{2}\right|
      \ge \frac{1}{2}, \quad
      \left|{2-\lambda_{N-1}}\right|
      - \left| \frac{1 - 2\lambda_{N-1} }{2}\right|
      \ge \frac{1}{2},
      \\
      & \left|2 - \beta_1 \mu_2\right|
      - \left| \mu_2\right| -  \left| \lambda_2\right|
      \ge \frac{1}{2}, \quad
      \left|2 - \beta_2 \lambda_{N-2}\right|
      - \left|\mu_{N-2}\right| -  \left|  \lambda_{N-2}\right|
      \ge \frac{1}{2},
    \end{array}
  \end{equation*}
   which, together with Theorem~\ref{thm:inftyNormEstiFromDiagDom},
   gives that $\|A_t^{-1}\|_\infty \le 2$. 
  Hence, we have
  \begin{equation}
    \label{eq:Atboundxi}
    \begin{array}{ll}
      \max\nolimits_{1 \le j \le N-1} |x_j|
      = \|\mathbf{x}_t\|_\infty
      \le \|A_t^{-1}\|_\infty \| \mathbf{y}_t\|_\infty
      \le  2\| \mathbf{y}\|_\infty.
    \end{array}
  \end{equation}
  \eqref{eq:transATEnd1} and \eqref{eq:transATEnd0} yield
  \begin{equation*}
    \label{eq:Atboundx0}
    \begin{array}{ll}
      |x_0| & = |y_0 - x_1 - (1-\beta_1)x_2| 
      = \left|y_0 - \left(\frac{y_1}{2-\mu_1} - \beta_1 x_2\right)
      - (1-\beta_1)x_2\right| \\
      & \le |y_0| + \frac{1}{2-\mu_1}|y_1| + (1 - 2\beta_1)|x_2|
      < 6 \|y\|_\infty,
    \end{array}
  \end{equation*}
  where the last inequality follows from 
  $\mu_1 \in (0,1), \beta_1 \in (-\frac{1}{2}, \frac{1}{4})$ and \eqref{eq:Atboundxi}.
  Similarly we have $|x_N| < 6 \|y\|_\infty$.
  Thus, the solution to $T_2 T_1 A T_4 \mathbf{x} = \mathbf{y}$
   satisfies
  \begin{equation*}
    \begin{array}{l}
      \|(T_2 T_1 A T_4)^{-1} \mathbf{y}\|_\infty
      = \|\mathbf{x}\|_\infty
      < 6 \|\mathbf{y}\|_\infty,
    \end{array}
  \end{equation*}
   indicating that $\|(T_2 T_1 A T_4)^{-1}\|_\infty \le 6$.
  Finally, we have
  \begin{displaymath}
    \|A^{-1}\|_\infty
    \le \|(T_2 T_1 A T_4)^{-1}\|_\infty 
    \| T_4\|_\infty  \| T_2T_1\|_\infty      
     \le 6  (1 + \max(|\beta_1|, |\beta_{2}|)) 2
      < 18,
  \end{displaymath}
   where the second inequality follows from \eqref{eq:TDef}
   and \eqref{eq:T4Def}. \qed
\end{proof}

\begin{thm}
  \label{thm:condA-inftyNorm}
  The matrix $A$ in \eqref{eq:miLinearSystem2}
  with $N \ge 5$
  satisfies $\cond_\infty(A) < 54$.
\end{thm}
\begin{proof}
  \eqref{eq:miLinearSystem2} implies 
   $\|A\|_\infty = \max\nolimits_{1 \le i \le N+1} 
    \sum_{j = 1}^{N+1} |a_{ij}| = 3$. 
  Then Lemma~\ref{lem:inverseA-inftyNorm} gives 
   $\cond_\infty(A) = \|A\|_\infty \|A^{-1}\|_{\infty} < 54$.  
   \qed
\end{proof}

}

\revise{
\subsubsection{The $(r, h)$-regularity
  of the breakpoint sequence of a spline}
\label{sec:r-h-regularity}

As the main conclusion of the previous subsection,
the excellent conditioning of periodic and not-a-knot splines
does not depend on the variation of chordal lengths
in the breakpoint sequence.
Nonetheless,
 for the purpose of good efficiency of marker usage, 
 some regularity on these chordal lengths
 is still desirable.

\begin{defn}[The $(r, h)$-regularity]
  \label{def:regularityCondition}
  A spline or its breakpoint sequence is 
  \emph{$(r, h)$-regular}
  if the maximum and minimal distances
  between all pairs of adjacent breakpoints
  are $h$ and $rh$, respectively. 
\end{defn}

The $(r, h)$-regularity in Definition \ref{def:regularityCondition}
 should be imposed on practical problems
 by specifying certain values of $r$ and $h$.
The upper bound $h$ can be calculated 
 according to the desired accuracy of the problem at hand
 via Theorems \ref{thm:periodicSplineErrorEstimates}
 and \ref{thm:notaKnotSplineErrorEstimates}. 
A value of $r$ close to 1 ensures
 a good efficiency of all breakpoints.
We recommend using \mbox{$r\in [\frac{1}{10},\frac{1}{3}]$}
 for most applications,
 with the upper bound $\frac{1}{3}$
 derived in the first paragraph
 in the proof of Lemma \ref{lem:rhRegularity}. 
 
However,
 when the variation of interface curvature 
 spans many orders of magnitude,
 one might want to accordingly select a small value of $r$ 
 so that high-curvature arcs
 \emph{can} be resolved by more densely populated markers. 
See Sec.~\ref{sec:curvatureBasedARMS}
 for related discussions. 
}

\subsection{Combining topology and geometry}
\label{sec:extractingSinglePhase}

While we definitely approximate a smooth closed curve 
 by a periodic cubic spline, 
 it might not be correct
 to approximate a curve segment $\gamma\in E_{\Gamma}$
 by a not-a-knot spline
 fitted through markers on $\gamma$.
In Fig.~\ref{fig:multiphase-representation},
 the smoothness of $\gamma_{4,1}^{1-}$ is lost
 if the edges $e_2,e_7,e_4$
 are approximated separately by not-a-knot splines.
Our solution to this problem starts with 

\begin{defn}
  \label{def:edgePairing}
  An \emph{edge pairing of an undirected graph} $G=(V,E,\psi)$
  is a set {$R^{\text{EP}}:=\{R^{\text{EP}}_v: v\in V\}$}, 
  where $R^{\text{EP}}_v \subseteq E_v\times E_v$,
  cf. (\ref{eq:adjacentEdgesAtv}), is
  a \emph{set of pairs of adjacent edges at $v$}
  such that each self-loop $e\in E_v$ appears and only appears
  in the pair $(e,e)$ 
  while any other edge in $E_v$
  appears at most once
  across all pairs in $R^{\text{EP}}_v$. 
\end{defn}


\begin{defn}
  \label{def:smoothnessIndicator}
  The \emph{smoothness indicator of an interface graph $G_{\Gamma}$}
  is an edge pairing of $G_{\Gamma}$ such that, 
  for each $v\in V_{\Gamma}$, each $(e_l, e_r)\in R^{\text{EP}}_v$
  indicates that $e_l$ and $e_r$ connect smoothly at $v$;
  in particular, $e_l=e_r$ corresponds to a smooth self-loop. 
\end{defn}

For the interface graph in Fig.~\ref{fig:multiphase-representation}(a),
 the smoothness indicator is given by
\begin{equation}
  \label{eq:RepVals}
  \small
  \begin{array}{l}
    R^{\text{EP}}_{v_1}=\left\{(e_2,e_4)\right\},\ 
    R^{\text{EP}}_{v_3}=\left\{(e_3,e_6)\right\},\ 
    R^{\text{EP}}_{v_4}=\left\{(e_2,e_7)\right\},\ 
    R^{\text{EP}}_{v_5}=\left\{(e_4,e_7)\right\}; 
    \\
    R^{\text{EP}}_{v_6}=\left\{(e_1,e_9),(e_{12},e_{14})\right\},\ 
    R^{\text{EP}}_{v_{7}}=\left\{(e_{10},e_{11})\right\},\
    R^{\text{EP}}_{v_{8}}=\left\{(e_1,e_9),(e_{12},e_{13}) \right\}; 
    \\
    R^{\text{EP}}_{v_{2}}=\left\{(e_8,e_8)\right\},
    R^{\text{EP}}_{v_{10}}=\left\{(e_{10},e_{15})\right\},
    R^{\text{EP}}_{v_{11}}=\left\{(e_{11},e_{16})\right\};
    R^{\text{EP}}_{v_{9}}=\emptyset;
  \end{array}
\end{equation}
where the first two lines correspond
 to T and X junctions, respectively.
If there exist multiple edges
 that connect smoothly to a given $e_l$, 
 we try to select an edge $e_r$ such that
 $e_r$ and $e_l$ belong to the same cycle $C_{i,j}^k$.
 
\begin{algorithm}
  \caption{: $(C_S, T_S)$ = \revise{\texttt{partitionEdgeSet}}
    ($G$,\ $R^{\text{EP}}$)}
  \begin{algorithmic}[1]
    \REQUIRE An undirected graph $G=(V,E,\psi)$
    in Definition \ref{def:graph}; \\
    \hspace{.6cm} the edge pairing $R^{\text{EP}}$ of $G$
    in Definition \ref{def:edgePairing}.
    \ENSURE A set $C_S$ of circuits and a set $T_S$ of trails
    as in Definition \ref{def:subgraphs}.
    \PostConditions
    (a) The edges in $C_S$ and $T_S$ partition $E$; \\
    \hspace{2.05cm}
    (b) $\forall \mathbf{e}\in (C_S\cup T_S)$,
    $\forall (e_l,e_r)$ adjacent in $\mathbf{e}$, 
    $\exists u\in V$ s.t. $(e_l,e_r)\in R^{\text{EP}}_u$.
    \STATE Initialize $C_S \leftarrow \emptyset$;
    $T_S \leftarrow \emptyset$ 
    \WHILE{$E \neq \emptyset$}
    \STATE Initialize a trail $\mathbf{e}\leftarrow (e)$
    with an arbitrary edge $e \in E$
    \STATE Initialize $(u_l,u_r)\leftarrow\psi(e)$
    and remove $e$ from $E$
    \STATE Initialize $e_1\leftarrow e$;\  $e_m\leftarrow e$;\ 
    $u_0\leftarrow u_l$;\  $u_m\leftarrow u_r$
    \hfill // $\mathbf{e}$ contains $m$ edges 
    \WHILE{$\exists (e_m,e_r)\in R_{u_m}^{\text{EP}}$
      \textbf{and} $e_r\in E$}
    \STATE $\mathbf{e}\leftarrow (\mathbf{e}, e_r)$;\ 
    $e_m\leftarrow e_r$;\  $u_m\leftarrow u_r$ where
    $(u_m,u_r)=\psi(e_r)$
    \STATE Remove $e_r$ from $E$
    \hfill // grow the right end of $\mathbf{e}$
    \ENDWHILE
    \WHILE{$\exists (e_l,e_1)\in R_{u_0}^{\text{EP}}$
      \textbf{and} $e_l\in E$}
    \STATE $\mathbf{e}\leftarrow (e_l, \mathbf{e})$;\ 
    $e_1\leftarrow e_l$;\  $u_0\leftarrow u_l$ where
    $(u_l,u_0)=\psi(e_l)$
    \STATE Remove $e_l$ from $E$
    \hfill // grow the left end of $\mathbf{e}$
    \ENDWHILE
    \IF{$u_{m} = u_{0}$ \textbf{and}
      $(e_{1}, e_{m})\in R^{\text{EP}}_{u_{0}}$}
    \STATE Add $\mathbf{e}$ to $C_S$
    \ELSE
    \STATE Add $\mathbf{e}$ to $T_S$
    \ENDIF
    \ENDWHILE
  \end{algorithmic}
  \label{alg:partitionGraph}
\end{algorithm}

\revise{
Given an undirected graph $G$
 and its edge pairing $R^{\text{EP}}$,
 Algorithm \ref{alg:partitionGraph} 
 constructs a set $C_S$ of circuits and a set $T_S$ of trails, 
 the union of which contains exactly the same edges
 as those in $E_G$.}
 
\begin{lem}
  \label{lem:algorithmIsCorrect}
  Algorithm \ref{alg:partitionGraph} stops and
  its post-conditions hold. 
\end{lem}
\begin{proof}
  By Sard's theorem and Definition \ref{def:YinSet},
  the total number of junctions and non-smooth points is finite
  and thus $E$ is also finite.
  Inside the three while loops,
  any edge added to the trail is immediately deleted
  and thus $\# E$
  decreases strictly monotonically
  for each while loop.
  Therefore, eventually we have $E=\emptyset$
  and the algorithm stops.

  The trail $\mathbf{e}$ is initialized at line 3 with a single edge.
  We grow $\mathbf{e}$ by appending edges in $R^{\text{EP}}$
  to its left and right ends
  until there are no edges in $E$
  to be paired with these ends.
  During this process,
  the removal of $e_r$ and $e_l$ at lines 8 and 12 implies 
  the distinctness of edges in $\mathbf{e}$,
  ensuring that $\mathbf{e}$ is indeed a trail,
  cf. Definition~\ref{def:subgraphs}.
  Also by Definition~\ref{def:subgraphs}, 
  the trail in line 15 is a circuit.
  Then post-condition~(a) follows from the classification in lines 14--18
  and the fact that, inside the outermost while loop, 
  all edges in $E$ have been visited. 
  Post-condition (b) follows from
  lines 6, 10, 14 as each edge added to the trail or circuit
  satisfies the pairing condition. \qed
\end{proof}

For the interface graph in Fig.~\ref{fig:multiphase-representation}(a)
 and the edge pairing in (\ref{eq:RepVals}),
 the output of Algorithm \ref{alg:partitionGraph}
 is shown in the last two columns
 of Fig.~\ref{fig:multiphase-representation}(c).
For the initial trail $(e_8)$,
 the two while loops in lines 6--13 are skipped
 and the condition at line 14 holds with 
 $u_0=u_m=v_2$ and $(e_8,e_8)\in R^{\text{EP}}_{v_2}$; 
 thus $(e_8)$ is added into $C_S$.
For the initial trails $(e_{10})$ or $(e_{11})$,
 the two while loops in lines 6--13 extend the trail to
 $\mathbf{e} = (e_{15}, e_{10}, e_{11}, e_{16})$
 due to (\ref{eq:RepVals}); 
 however, the condition $(e_{15},e_{16})\in R^{\text{EP}}_{v_7}$
 at line 14 does not hold,
 so the trail $\mathbf{e}$ is added into $T_S$.
Similarly,
 $(e_{13},e_{12},e_{14})$ in Fig.~\ref{fig:multiphase-representation}
 violates the pairing condition at line 14
 and is also added into $T_S$.
 
The input parameter $G$ of Algorithm \ref{alg:partitionGraph}
 is not required
 to have the structure of an interface graph or even a planar graph; 
 similarly, $R^{\text{EP}}$ is not the smoothness indicator
 in Definition \ref{def:smoothnessIndicator}
 but the edge pairing in Definition \ref{def:edgePairing}. 
Nonetheless,  $R^{\text{EP}}$ is \emph{interpreted}
 as the smoothness indicator of the interface graph
 so that the output $C_S$ and $T_S$ correspond respectively
 to smooth closed curves approximated by periodic splines
 and to smooth curve segments
 approximated by not-a-knot splines. 
$C_S\cup T_S$ is not isomorphic to $E_{\Gamma}$ 
 and neither is $C_S$ to $C$: 
 circuits preserve smoothness in fitting splines
 while cycles represent the topology of each Yin set,
 cf. Definition \ref{def:subgraphs}.
This discussion suggests the need
 for some set of splines isomorphic to $C_S\cup T_S$.

\begin{defn}
  \label{def:setOfFittedSplines}
    The \emph{set $S_{CT}$ of fitted splines} 
    is built by first concatenating marker sequences in $E_X$
    in (\ref{eq:setOfMarkSequences})
    according to circuits or trails in $C_S\cup T_S$
    and then interpolating each concatenated breakpoint sequence:
    a periodic spline for a circuit
    and a not-a-knot spline for a trail (that is not a circuit).
\end{defn}

Cutting the fitted splines in $S_{CT}$
 at junctions and non-smooth points,  
 we obtain the spline edge set $S_E$
 in Definition \ref{def:splineBoundarySegs}.

\revise{
\begin{ntn}
  \label{ntn:splineApproxOfGamma}
  Denote by $\widetilde{\Gamma}:= (C, S_E)$
  a \emph{spline approximation of} the interface $\chi(\Gamma)$,
  where
  the cycle set $C$ in Notation \ref{ntn:cycleSet}
  records the topology of $\chi(\Gamma)$
  and the spline edge set $S_E$
  in Definition \ref{def:splineBoundarySegs} 
  approximate the geometry of $\chi(\Gamma)$.
\end{ntn}

}


\begin{figure}
  \centering
  \includegraphics[width=0.95\textwidth]{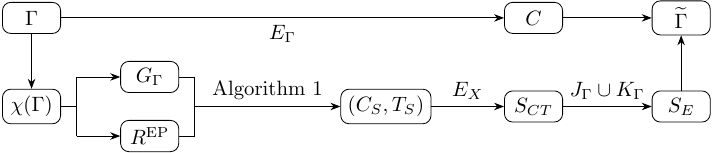}
  \caption{The pipeline of constructing $\widetilde{\Gamma}$
    as a spline approximation of 
    the poset $\Gamma$ that contains the topology and geometry
    of all Yin sets in Notation \ref{ntn:B-REP}.
    Recall that $\chi(\Gamma)$
    denotes the interface in Notation \ref{ntn:B-REP},
    $G_{\Gamma}=(V_\Gamma, E_\Gamma, \psi_\Gamma)$
    the interface graph in Definition \ref{def:topoGraph}, 
    $R^{\text{EP}}$ the smoothness indicator of $G_{\Gamma}$
    in Definition \ref{def:smoothnessIndicator}, 
    $C_S$ and $T_S$ the sets of circuits and trails
    generated by Algorithm \ref{alg:partitionGraph}, 
    $E_X$ the set of marker sequences
    in Definition \ref{def:splineBoundarySegs}, 
    $S_{CT}$ the set of fitted splines
    in Definition \ref{def:setOfFittedSplines}, 
    and $J_{\Gamma}$ and $K_{\Gamma}$
    the set of junctions in Definition \ref{def:junctionPoint}
    and non-smooth points in Definition \ref{def:kink}, respectively.
    The cycle set $C$ and the spline edge set $S_E$
    constitute $\widetilde{\Gamma}$
    in Notation \ref{ntn:splineApproxOfGamma}.
  }
  \label{fig:pipeline}
\end{figure}

Our design of the boundary representation of multiple Yin sets
 is concisely summarized in Fig.~\ref{fig:pipeline}.
By Theorem \ref{thm:uniqueCases},
 each of the $N_p$ Yin sets
 in ${\cal M}$ is uniquely represented 
 by a set $\Gamma_{i}=\{\Gamma_{i,j}\}$
 of posets of oriented Jordan curves. 
\revise{The topology of $\chi(\Gamma)$ is captured 
 in the interface graph
 $G_{\Gamma}=(V_{\Gamma}, E_{\Gamma},\psi_{\Gamma})$
 as well as by the cycle set $C$. 
Given $G_{\Gamma}$ and
 the smoothness indicator $R^{\text{EP}}$, 
 Algorithm \ref{alg:partitionGraph} partitions $E_{\Gamma}$
 into $(C_S, T_S)$,
 which, together with the set $E_X$ of marker sequences, 
 yields the set $S_{CT}$ of fitted splines.
Elements in $C_S\cup T_S$ have a one-to-one
 correspondence to the splines in $S_{CT}$. 
Since the edges in $C_S\cup T_S$
 form a partition of all edges in $E_{\Gamma}$, 
 $S_{CT}$ covers all arcs in $\chi(\Gamma)$, 
 so do $S_E$ since it is obtained 
 by cutting splines in $S_{CT}$
 at junctions and non-smooth points.
In our implementation,
 $S_E$ is simply an array of splines
 and each directed cycle $C_{i,j}^k$ in $C$
 a sequence of indices of the splines in $S_E$, 
 cf. Fig.~\ref{fig:multiphase-representation}(c).
}
 


\section{Evolving the interface
 with a general ARMS strategy}
\label{sec:algorithms}
\revise{At the initial time $t_0$,
  we represent the initial Yin sets ${\cal M}(t_0)$
  by the poset $\Gamma(t_0)$ in Notation \ref{ntn:B-REP} 
  and employ the pipeline in Fig.~\ref{fig:pipeline}
  to produce the data
  $C$, $G_{\Gamma}$, $R^{\text EP}$, $S_{CT}^0$, and $S_E^0$, 
  as explained in the ending paragraph of the previous section. 
  $\widetilde{\Gamma}(t_0):=(C,S^0_E)$
  amounts to an approximation of the initial condition $\Gamma(t_0)$.
  Then we numerically solve the multiphase IT problem
  in Definition \ref{def:interfaceTracking}
 by evolving $S_{CT}^0$ from $t=t_0$ to $t=T$
 over a finite number of time steps with uniform size $k$.
During this process,
$C$, $G_{\Gamma}$, and $R^{\text EP}$ remain constant.
We emphasize that it is $S_{CT}^0$ instead of $S_{E}^0$
 that is evolved because
 $S_{CT}^0$ already contains the topological and geometric information
 in $C_S$, $T_S$, and $R^{\text EP}$. 
 
Hereafter, a computed value at time $t_n=t_0+nk$
 is indicated by a superscript $^n$, 
 e.g., $S_{CT}^n$ denotes the set of fitted splines at $t_n$
 and ${\cal M}^n$ the IT result
 given by $C$ and $S_E^n$.
}

\begin{defn}
  \label{def:MARS}
  A \emph{MARS method} is an IT method of the form
  \begin{equation}
    \label{eq:mars}
    \mathcal{M}^{n+1} :=
    (\chi_{n+1} \circ \varphi_{t_n}^k \circ \psi_n) \mathcal{M}^n,
  \end{equation}
  where $\mathcal{M}^n \in \mathbb{Y}$ is an approximation of
  $\mathcal{M}(t_n) \in \mathbb{Y}$,
  $\varphi_{t_n}^k: \mathbb{Y} \rightarrow \mathbb{Y}$
  a fully discrete mapping operation
  that approximates the exact flow map 
  in (\ref{eq:flowMapActsOnYinSet}), 
  $\psi_n : \mathbb{Y} \rightarrow \mathbb{Y}$
  an augmentation operation at $t_n$
  to prepare $\mathcal{M}^n$ for $\varphi_{t_n}^k$,
  and $\chi_{n+1}: \mathbb{Y} \rightarrow \mathbb{Y}$
  an adjustment operation after the mapping $\varphi_{t_n}^k$.
\end{defn}

\revise{
  Previously, we have developed two cubic MARS methods,
  one \cite{zhang2018fourth} for
  the set $\Phi$ of flow maps in \eqref{eq:groupOfDiffeomorphisms}
  coming from the ODE \eqref{eq:velODE}
  while the other \cite{hu2025:_arms}
  for the scenario where $\Phi$ is implicitly specified
  by a curve shortening flow.
  The success of these MARS methods
  hinges on the ARMS strategy that
  maintains the $(r_{\Tiny},h_L)$-regularity
  by adding and removing interface markers,
  which correspond to
  $\psi_n$ and $\chi_{n+1}$ in (\ref{eq:mars}), respectively. 
  However, both methods are limited to two-phase flows. 
  In addition,
  since the interface curves
  are approximated only by periodic splines,
  these methods lack an intrinsic mechanism
  to handle non-smooth points accurately and efficiently. 

 In light of 
  Sec.~\ref{sec:representation}, 
  the key to generalizing the aforementioned cubic MARS methods
  to multiphase IT
  is how to maintain the $(r_{\Tiny},h_L)$-regularity
  of the splines in $S_{CT}^n$ under the action of flow maps in $\Phi$. 
 The answer lies in two main aspects.
 In Sec.~\ref{sec:arms-strategy},
  we propose a more general ARMS strategy for evolving $S_{CT}^n$, 
  with particular attention 
  to retaining junctions and non-smooth points
  in the marker sequences.
 In Sec.~\ref{sec:curvatureBasedARMS}, 
  we vary the value of $h_L$ according to the local curvature
  so that the IT error,
  rather than the marker distribution, 
  tends to be uniform along the interface.
 Finally in Sec.~\ref{sec:cubicMARS},
  the two components 
  assemble into the multiphase cubic MARS method.
}


\subsection{The ARMS strategy that retains characteristic markers}
\label{sec:arms-strategy}

\begin{figure}
  \centering
  \subfloat[trace forward in time the breakpoints (interface markers) of $\partial{\cal M}^n$]{
    \includegraphics[width=0.47\linewidth]{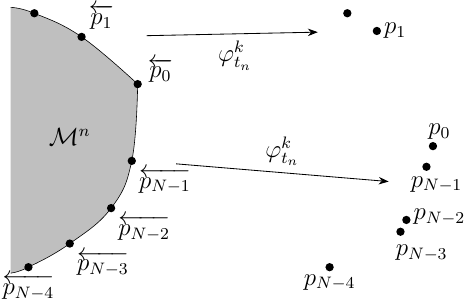}
  }
  \hfill
  \subfloat[enforce $h_L^*=(1-2r_{\Tiny})h_L$ 
  as the upper bound of chordal lengths by adding new marker $q_0$]{
    \includegraphics[width=0.47\linewidth]{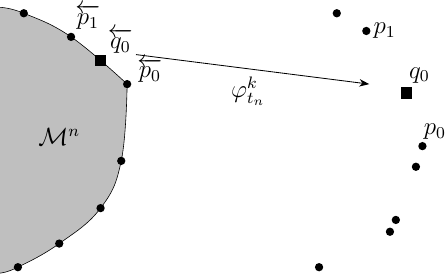}
  }

  \subfloat[enforce $r_{\Tiny}h_L$ as the lower bound of
  chordal lengths %
  by removing the marker $p_{N-1}$ near the kink $p_0$
  and the marker $p_{N-2}$ near $p_{N-3}$]{
    \includegraphics[width=0.47\linewidth]{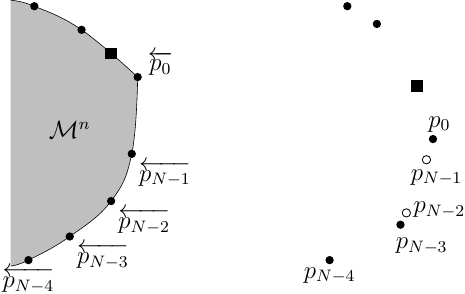}
  }
  \hfill
  \subfloat[fit a new spline 
  as (part of) the boundary $\partial{\cal M}^{n+1}$]{
    \includegraphics[width=0.47\linewidth]{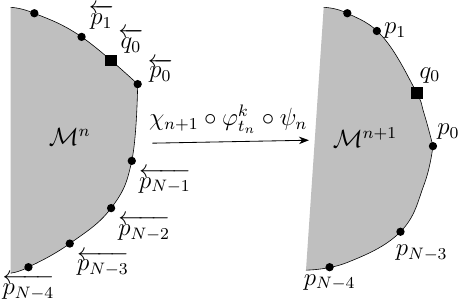}
  }
  \caption{
    The ARMS strategy.
    In subplot (a), 
    interface markers on $\partial{\cal M}^n$ are mapped 
    to their images
    by the discrete flow map $\varphi_{t_n}^{k}$.
    In subplot (b), 
    the distance between $p_0$ and $p_1$ 
    is larger than $h_L^*$ and thus 
    a new marker $\protect\overleftarrow{q_0}$
    (the solid square) is added
    and the new preimage sequence
    is mapped to time $t_n+k$ 
    so that distances between both $p_0, q_0$ and $q_0, p_1$ 
    are smaller than $h_L^*$.
    \revise{In subplot (c), 
    the distance between $p_0$ and $p_{N-1}$ 
    is found to be smaller than the lower bound $r_{\Tiny}h_L$
    the marker $p_{N-1}$ (the hollow circle) is removed
    because the kink $p_0$ must remain in the sequence. 
    Also, the marker $p_{N-2}$ is removed
    since the distance between $p_{N-2}$ and $p_{N-3}$ 
    is smaller than $r_{\Tiny}h_L$.
    In subplot (d),
    we obtain (part of) $\partial{\cal M}^{n+1}$
    by fitting a new not-a-knot spline 
    through the updated chain of markers 
    ``$p_0 \rightarrow q_0 \rightarrow p_1 \rightarrow \cdots
    \rightarrow p_{N-4} \rightarrow p_{N-3} 
    \rightarrow p_0$.''}     
  }
  \label{fig:ARMS}
\end{figure}

\revise{
The original ARMS strategy \cite[Def. 3.5]{hu2025:_arms}
 acts on periodic splines only 
 and all breakpoints are treated equally 
 in terms of being added or removed.
When generalizing ARMS to not-a-knot splines, 
 we encounter the essential difficulty
 that not all breakpoints are created equal: 
 vertices of the interface graph $G_{\Gamma}$ are special
 because they carry topological information
 such as the basepoint of a cycle
 and the two ends of a trail; 
 hence they should never be removed 
 during a homeomorphic IT
 so that the topological data in $G_{\Gamma}$ and $C$
 \emph{indeed} remain constants.

In Definition \ref{def:ARMS},  
 we formalize the generic ARMS strategy
 that applies not only to periodic splines
 but also to not-a-knot splines. 
The input set $\mathbf{z}^n$
 of characteristic points of $\mathbf{s}^n$
 is supposed to contain all vertices of the interface graph.
As shown in (ARMS-3) and Fig.~\ref{fig:ARMS} (c),
 markers in $\mathbf{z}^n$ are retained
 when enforcing the lower bound of chordal lengths. 
In Algorithm \ref{alg:ARMS},
 we elaborate steps of the generic ARMS strategy 
 for better reproducibility. 
} 

\begin{defn}[General ARMS strategy for a single spline]
  \label{def:ARMS}
  Given
  \begin{itemize}
  \item a discrete flow map
    $\varphi_{t_n}^k: \mathbb{Y}\rightarrow \mathbb{Y}$, 
  \item a periodic or not-a-knot cubic spline $\mathbf{s}^n$ 
    whose breakpoint sequence $(X_i)_{i=0}^{N^n}$
    is $(r_{\Tiny},h_L)$-regular in the sense of
    Definition \ref{def:regularityCondition}, 
  \item a subset $\mathbf{z}^n\subset (X_i)_{i=0}^{N^n}$
  that characterizes $\mathbf{s}^n$, 
  \end{itemize}
  the \emph{ARMS strategy} generates
  from $(\varphi_{t_n}^k, \mathbf{s}^n, \mathbf{z}^n, r_{\Tiny}, h_L)$
  a pair $(\mathbf{s}^{n+1}, \mathbf{z}^{n+1})$, 
  where $\mathbf{z}^{n+1}$ is the set
  of characterizing breakpoints of $\mathbf{s}^{n+1}$,
  by steps as follows. 
  \begin{enumerate}[label=(ARMS-\arabic*), leftmargin=5em]\itemsep0em
  \item Initialize $(p_i)_{i=0}^{N^{n+1}}$
    with $N^{n+1} \gets N^{n}$ and 
    $p_i \leftarrow\varphi_{t_n}^k(X_i)$; 
    also set \mbox{$\mathbf{z}^{n+1}= \varphi_{t_n}^k(\mathbf{z}^n)$}.
    \label{en:arms-init}
  \item For a chordal length $\|p_j-p_{j+1}\|_2$ 
    greater than $h_L^*:=(1-2r_{\Tiny})h_L$, 
    \begin{enumerate}
    \item locate $X_j=(x(l_j), y(l_j))$ and 
      $X_{j+1}=(x(l_{j+1}), y(l_{j+1}))$
      on $\mathbf{s}^n$
      as preimages of $p_j$ and $p_{j+1}$,
    \item divide the parametrization interval $[l_j, l_{j+1}]$ 
      into $\left\lceil\frac{\|p_j-p_{j+1}\|_2}{h_L^*}\right\rceil$ 
      equidistant subintervals, 
      compute the corresponding new markers on $\mathbf{s}^n(l)$,
      insert them between $X_j$ and $X_{j+1}$, and
    \item insert the images of new markers under $\varphi_{t_n}^k$ 
      into the new sequence between $p_j$ and $p_{j+1}$.
    \end{enumerate}
    Repeat (a, b, c) until no chordal length 
    is greater than $h_L^*$.
    \label{en:arms-add}
  \item Remove chords of negligible lengths
    from the sequence $(p_i)_{i=0}^{N^{n+1}}$: 
    \begin{enumerate}
    \item for each $p_j\in \mathbf{z}^{n+1}$,
      keep removing $p_{j+1}$ from the breakpoint sequence
      until $\|p_j-p_{j+1}\|_2\ge r_{\Tiny}h_L$ holds
      and keep removing $p_{j-1}$ from the sequence
      until $\|p_j-p_{j-1}\|_2\ge r_{\Tiny}h_L$ holds, 
    \item locate $p_{\ell}\in (X_i)$ satisfying
      $\|p_{\ell}-p_{\ell+1}\|_2< r_{\Tiny}h_L$ and set $j=\ell$, 
    \item if $\|p_j-p_{j+1}\|_2< r_{\Tiny}h_L$, 
      keep removing $p_{j+1}$ from the point sequence until 
      $\|p_j-p_{j+1}\|_2\ge r_{\Tiny}h_L$ holds for the new $p_{j+1}$,
    \item increment $j$ by 1 and repeat (b, c)
      until all chords are checked.
    \end{enumerate}
    \label{en:arms-remove}
  \item Construct spline $\mathbf{s}^{n+1}$
    from the new breakpoint sequence $(p_i)_{i=0}^{N^{n+1}}$. 
    \label{en:arms-construct}
  \end{enumerate}
  We have assumed cyclic indexing for periodic splines.
\end{defn}

\begin{algorithm}
  \caption{: 
    \revise{$(\mathbf{s}^{n+1}, \mathbf{z}^{n+1})$ = \texttt{ARMS}
      ($\varphi_{t_n}^k$,\ $\mathbf{s}^n$,\ $\mathbf{z}^n$,\
      $r_{\Tiny}$, $h_L$)}}
  \begin{algorithmic}[1]
    \REQUIRE A discrete flow map
    $\varphi_{t_n}^k: \mathbb{Y}\rightarrow \mathbb{Y}$;
    \\ \hspace{.6cm}
    a cubic spline $\mathbf{s}^n$
    and a subset $\mathbf{z}^n$ of
    the breakpoint sequence of $\mathbf{s}^n$; 
    \\ \hspace{.6cm}
    the regularity parameters $(r_{\Tiny}, h_L)$.
    \ENSURE A pair $(\mathbf{s}^{n+1}, \mathbf{z}^{n+1})$
    where $\mathbf{z}^{n+1}$ is a subset of
    the breakpoint sequence of $\mathbf{s}^{n+1}$.
    
    \PreConditions \hspace{0.3mm}
    $r_{\Tiny}\in (0,\frac{1}{3})$ and 
    $\mathbf{s}^n$ is $(r_{\Tiny},h_L)$-regular
    \PostConditions The breakpoint sequence $(X_i)_{i=0}^{N^{n+1}}$ 
     of $\mathbf{s}^{n+1}$ is $(r_{\Tiny},h_L)$-regular.
    
    \Comment  // trace forward in time
    \STATE Initialize $N^{n+1} \leftarrow N^{n}$ and
    $\mathbf{z}^{n+1} \leftarrow \varphi_{t_n}^k(\mathbf{z}^n)$
    \STATE Initialize a sequence $\mathbf{p}:=(p_i)_{i=1}^{N^{n+1}}$
    with $p_i \leftarrow \varphi_{t_n}^k(X_i)$

    \Comment // enforce $h_L^{\ast}$ as the upper bound of chordal lengths
    \STATE Set $h_L^* \leftarrow (1-2r_{\Tiny})h_L$

    \WHILE{$\exists j$ such that $\|p_j-p_{j+1}\|_2 > h_L^*$}
      \STATE Locate $X_j, X_{j+1}$ on $\mathbf{s}^n$ as preimages of $p_j, p_{j+1}$
      \STATE Calculate $M \leftarrow
      \left\lceil \frac{1}{h_L^*}\|p_j-p_{j+1}\|_2  \right\rceil$
      where $\lceil\cdot\rceil$ is the ceiling function
      \STATE Locate $M$ equidistant markers on $\mathbf{s}^n$
      between $X_j$ and $X_{j+1}$
      \STATE Insert images of these $M$ markers
      under $\varphi_{t_n}^k$ in between $p_j$  and $p_{j+1}$
    \ENDWHILE

    \Comment // enforce $r_{\Tiny}h_L$ as the lower bound
    of chordal lengths near the markers in $\mathbf{z}^{n+1}$
    \FORALL{$p_j \in \mathbf{z}^{n+1}$}
      \WHILE{$\|p_j-p_{j+1}\|_2 < r_{\Tiny}h_L$}
        \STATE Remove $p_{j+1}$ from the sequence $\mathbf{p}$
      \ENDWHILE
      \WHILE{$\|p_j-p_{j-1}\|_2 < r_{\Tiny}h_L$}
        \STATE Remove $p_{j-1}$ from the sequence $\mathbf{p}$
      \ENDWHILE
    \ENDFOR

    \Comment // enforce $r_{\Tiny}h_L$ as the lower bound
    of chordal lengths for other markers
    \STATE Initialize index $j \leftarrow 0$
    \WHILE{$j < \text{length}(\mathbf{p})-1$}
    \WHILE{$\|p_j-p_{j+1}\|_2 < r_{\Tiny}h_L$}
          \STATE Remove $p_{j+1}$ from the sequence $\mathbf{p}$
        \ENDWHILE
        \STATE $j \leftarrow j + 1$
    \ENDWHILE
    \STATE Construct a periodic or not-a-knot spline $\mathbf{s}^{n+1}$
    from the sequence $\mathbf{p}$ 
  \end{algorithmic}
  \label{alg:ARMS}
\end{algorithm} 

\revise{
Instead of $h_L^* =(1 - r_{\Tiny})h_L$
 in the original ARMS strategy \cite[Def. 3.5]{hu2025:_arms}, 
 $h_L^* = (1 - 2r_{\Tiny})h_L$ is used in (ARMS-2)
 so that the $(r_{\Tiny},h_L)$-regularity can be strictly maintained
 even at the presence of special markers in $\mathbf{z}^n$.
 }

\begin{lem}
  \label{lem:rhRegularity}
  Provided that its preconditions hold,
  the ARMS strategy in Definition \ref{def:ARMS}
  generates a periodic or not-a-knot spline
  that is $(r_{\Tiny},h_L)$-regular.
\end{lem}
\begin{proof}
  If $r_{\Tiny}\ge \frac{1}{2}$,
  we have $h_L^* = (1 - 2r_{\Tiny})h_L\le 0$
  and the ARMS strategy cannot continue.
  The reasonable choice of $h_L^*\ge r_{\Tiny}h_L$
  yields $r_{\Tiny}\le \frac{1}{3}$.
  
  Although \ref{en:arms-add} ensures
  that no chordal length is greater than 
  $h_L^{\ast}$,
  the removal of certain markers in enforcing 
  the lower bound of $r_{\Tiny}h_L$ on chordal lengths
  during \ref{en:arms-remove} 
  may increase the maximum chordal length.

  As shown in Fig.~\ref{fig:ARMS}(c), 
   neighboring markers of characterizing breakpoints in 
   $\mathbf{z}^{n+1}$ may be removed in (ARMS-3a), 
   resulting in the maximum chordal length being increased to 
   $h_L^{\ast} + r_{\Tiny}h_L$.
  Also, the removal of ordinary markers in (ARMS-3c)
   may further increase the chordal length to
   $h_L^{\ast} + r_{\Tiny}h_L + r_{\Tiny}h_L$,
   which equals $h_L$
   and is the largest possible value of the maximum chordal length.
  Then the proof is completed by 
  Definition \ref{def:regularityCondition}. \qed
\end{proof}

\revise{
\subsection{The $(r_{\Tiny}, h_L)$-regularity:
 constant $h_L$ vs. curvature-dependent $h_L$}
\label{sec:curvatureBasedARMS}
}

Although the Lagrangian grid of moving markers suffices to
 evolve the interface, 
 an Eulerian grid is needed
 to couple an IT method with a main flow solver.
Assuming for simplicity that
 the Eulerian grid has a uniform size $h$
 along each dimension, 
 we can specify a constant $h_L = r_h h^{\alpha}$
 to relate the length scale of the interface
 to that of the bulk flow.
Since the interface is a set of codimension one
 and an IT has at least $O\left(\frac{1}{h}\right)$ time steps,
 the complexity of a MARS method with $h_L = r_hh^{\alpha}$
 is $O\left(\frac{1}{h^{1+\alpha}}\right)$.
In contrast,
 the optimal complexity of a main flow solver 
 is $O\left(\frac{1}{h^3}\right)$ in two dimensions.
Thus it is fine to choose $\alpha\le 2$ 
 as this does not increase the complexity of the entire solver. 

Appropriate values of $r_h$ and $\alpha$
depend on the ratio of the length scale of the interface
to that of the bulk flow; 
in the usual case of the former
being much smaller than the latter,
we recommend choosing $r_h\ll 1$ with $\alpha=1$
or even $\alpha>1$ with $r_h< 1$
so that the numerical condition $O(\frac{h_L}{h})\ll 1$
reflects the real physics.
This complete decoupling of $h_L$ from $h$
is an advantage of cubic MARS methods; 
it is difficult for VOF methods
to have subgrid resolution due to their interface
reconstruction from volume fractions of the Eulerian grids.

With $h_L= r_h h^{\alpha}$
 and a constant $r_{\Tiny}\approx 0.1$,
 the ARMS strategy 
 maintains a roughly uniform distribution of markers
 along the interface. 
However, 
 for IT problems with orders-of-magnitude variations
 of interface curvature,
 the limited range of $[r_{\Tiny}h_L, h_L]$ 
 might result in very large errors at high-curvature markers
 and/or unnecessarily small errors at low-curvature markers;
 the former deteriorates accuracy
 while the latter deteriorates efficiency.

These adverse effects are 
 caused by the equidistribution of markers
 and can be remedied by varying the marker density
 according to the local curvature,
 which eventually leads to an equidistribution of IT errors
 along the interface.
Intuitively, over the same arc length,
 more markers are needed to resolve a high-curvature arc
 than a low-curvature arc.

Following \cite{hu2025:_arms}, 
 we fix $r_{\Tiny}$ as a constant
 and vary $h_L$ according to the \emph{radius of curvature}
 $\rho := \frac{1}{|\kappa|}$.
If the fitted spline is ${\cal C}^2$ at the $i$th marker $X_i$,
 we calculate $\rho_i:=\rho(X_i)$
 from the piecewise polynomials at $X_i$; 
 otherwise we compute the two radii $\rho_i^{\pm}$
 by one-sided differentiation of polynomials at the two sides of $X_i$.
The variation of $\rho_i$'s
 across all interface markers in all time steps
 is characterized by 
 \begin{equation}
   \label{eq:radiiVariationPars}
   \begin{array}{l}
   \rho_{\min} := \max_n\left(\rho_{\min}^c,  {\min}_i \rho_i^n\right); 
   \quad
   \rho_{\max} := \min_n\left(\rho_{\max}^c, {\max}_i \rho_i^n\right), 
   \end{array}
 \end{equation}
 where $\rho_i^n$ is the curvature of the $i$th marker $X_i$
 at the $n$th time step
 and the interval $[\rho_{\min}^c, \rho_{\max}^c]$
 should cover almost the entire (say, 99\%) range 
 of the radii of curvature
 during the simulation. 
The user-specified constants $\rho_{\min}^c>0$ and $\rho_{\max}^c$
 prevent the distraction from linear segments with zero curvature
 and arcs with exceptionally high curvature,
 respectively. 

In this work, 
 the function $h_L$ is designed to be 
 monotonically increasing, 
 \begin{equation}
   \label{eq:curvature-based-arms}
   h_L(\rho)
   := \begin{cases}
     r_{\min} h_L^c & \text{if} ~ \rho \le \rho_{\min};
     \\
     r_{\min} h_L^c + (1 - r_{\min}) h_L^c \cdot
     \sigma^c \left(
       \frac{\rho^{} - \rho_{\min}^{}}{\rho_{\max}^{} -
         \rho_{\min}^{}} \right)
     & \text{if} ~ \rho_{\min} < \rho < \rho_{\max};
     \\
     h_L^c & \text{if} ~ \rho \ge \rho_{\max},
   \end{cases}
\end{equation}
where $h_L^c$ denotes 
the maximum allowable distance between any two adjacent markers, 
the continuous bijection $\sigma^c: [0, 1] \rightarrow [0, 1]$ 
must satisfy $\sigma^c(0) = 0$ and $\sigma^c(1) = 1$, 
the parameter $r_{\min}$ is defined as
$r_{\min} := \max \left( r_{\min}^c,
  \frac{\rho_{\min}}{\rho_{\max}} \right)$, 
and $r_{\min}^c\in (0,1]$ is a user-specified constant 
for controlling 
$R_{\max}:=\frac{1}{r_{\min}^c r_{\Tiny}}$, 
the \emph{highest possible ratio
  of the longest chordal length over the shortest one}.
If we had omitted $r_{\min}^c$ 
 to define $r_{\min}$ as $\frac{\rho_{\min}}{\rho_{\max}}$,
 then $R_{\max}$ would be
 $\frac{\rho_{\max}^c}{\rho_{\min}^c r_{\Tiny}}$, 
 which might be much larger, 
 e.g., 
 the test case in Sec.~\ref{sec:testDeformation}
 gives $\frac{\rho_{\max}^c}{\rho_{\min}^c r_{\Tiny}}=2\times 10^6$, 
 cf. \eqref{eq:curvatureStrategyConstants}. 
We remark that, even for this large ratio of chordal lengths,
 our cubic spline fitting is still well conditioned;
 see the last paragraph of Sec.~\ref{sec:geometryOfCommonBoundary}.
 
We refer to (\ref{eq:curvature-based-arms}) as
 the \emph{curvature-based formula for the maximum chordal length}.
To demonstrate its adaptive power,  
suppose a line segment has its initial marker density at $\frac{1}{h_L^c}$
 and its curvature is increasing. 
Then (\ref{eq:curvature-based-arms}) dictates
 that more markers will be added to the arc,
 increasing the marker density up to $\frac{1}{r_{\min}^ch_L^c}$.
If the distances between adjacent markers are further reduced
 by the flow map,
 the marker density may continue to increase
 up to  $\frac{1}{r_{\Tiny}r_{\min}^ch_L^c}$.
Therefore, 
 the maximum increase/decrease ratio of the marker density
 is $R_{\max}$.

To unleash the power of \eqref{eq:curvature-based-arms}, 
 one needs to specify, for the problem at hand,
 the form of $\sigma^c$
 and all values of
 $\rho^c_{\min}$, $\rho^c_{\max}$, $h_L^c$, and $r_{\min}^c$.
In light of discussions in this subsection, 
 we summarize, for the most general scenario, 
 how to determine these user-specified constants.
 \begin{enumerate}[label=(USC-\arabic*)]
 \item Select $r_{\Tiny}$ based on how uniform 
   the markers should be on a curve with \emph{constant} curvature;
   for example, $r_{\Tiny}=0.1$ bounds
   the variation of distances by 10
   between adjacent markers on circular arcs or line segments. 
 \item Set $h_L^c=r_h h^{\alpha}$ where
   the values of $r_h$ and $\alpha$
   are determined according to
   the ratio of the interface length scale 
   to that of the bulk flow.
 \item Choose $\sigma^c$ to be $\sigma^c(x)=x$;
   other possibilities are $\sigma^c(x) = x^2$, $x^{\frac{1}{2}}$,
   etc.
 \item Estimate the range of all possible values
   of radii of curvature
   to determine $\rho^c_{\min}$ and $\rho^c_{\max}$; 
   see the paragraph that contains \eqref{eq:radiiVariationPars}.
 \item Pick a value for the maximum allowed ratio $R_{\max}$
   (of the longest chordal length over the shortest one)
   and calculate $r_{\min}^c=\frac{1}{r_{\Tiny}R_{\max}}$.
 \end{enumerate}

\subsection{The multiphase cubic MARS method}
\label{sec:cubicMARS}

We combine the previous subsections to propose

\begin{defn}[The multiphase cubic MARS method]
  \label{def:MARS-for-multiple-phases}
  Given
  \begin{itemize}
  \item a discrete flow map $\varphi_{t_n}^k$ 
    that approximates a homeomorphic flow map $\phi$, 
  \item the cycle set $C$ and the interface graph $G_{\Gamma}$
  of the initial condition ${\cal M}(t_0)$, 
  \item a pair $(S_{CT}^0, Z_{CT}^0)$
  where $S_{CT}^0$ 
  is the set of fitted splines in Definition \ref{def:setOfFittedSplines}
  and the function $Z_{CT}^0: S_{CT}^0 \to V_{\Gamma}$
  is given by
  \mbox{$Z_{CT}^0(\mathbf{s}) = \mathbf{s}\cap V_{\Gamma}$}, 
  \item values of the parameters $r_{\Tiny}$ and $h_L$
    where $h_L$ is given in \eqref{eq:curvature-based-arms}
    with the user-specified constants determined
    by (USC-1--5), 
  \end{itemize}
  \emph{the multiphase cubic MARS method} for the IT problem
  in Definition \ref{def:interfaceTracking}
  advances $(S^n_{CT},Z_{CT}^n)$
  to $(S^{n+1}_{CT},Z_{CT}^{n+1})$ as follows.
  \begin{enumerate}[label=(\alph*)]
  \item For each spline $\mathbf{s}^n\in S^n_{CT}$
    and its characterization set
    $\mathbf{z}^n=Z_{CT}^n(\mathbf{s}^n)$, 
    obtain $\mathbf{s}^{n+1}\in S^{n+1}_{CT}$ and $\mathbf{z}^{n+1}$
    by invoking 
    \texttt{ARMS}
      ($\varphi_{t_n}^k$,\ $\mathbf{s}^n$,\ $\mathbf{z}^n$,\
      $r_{\Tiny}$, $h_L$)
      in Algorithm \ref{alg:ARMS}. 
    All pairs in $\left\{(\mathbf{s}^{n+1},\mathbf{z}^{n+1}):
    \mathbf{s}^{n+1}\in S^{n+1}_{CT}\right\}$
    constitute $Z_{CT}^{n+1}$. 
  \item (optional)
    Assemble $\widetilde{\Gamma}^{n+1}$
    by first converting $S^{n+1}_{CT}$ to $S^{n+1}_E$
    and then connect splines in $S_E^{n+1}$
    according to the cycle set $C$;
    see Fig.~\ref{fig:pipeline}. 
  \end{enumerate}
\end{defn}

Step (b) in Definition \ref{def:MARS-for-multiple-phases}
 is optional since $\widetilde{\Gamma}^{n+1}$
 is not needed in evolving the interface.
Nonetheless, an IT method should 
 readily generate an approximate interface whenever needed, 
 hence step (b) is necessary. 


\revise{ 
\begin{lem}
  \label{lem:complexityOfARMS}
  For the multiphase MARS method,
  the memory usage during the entire simulation
  is $O(\frac{1}{h_L})$, 
  so is its CPU time per time step.
\end{lem}
\begin{proof}
  By Definition \ref{def:YinSet}, 
  any Yin set contains only a finite number
  of junctions and non-smooth points.
  Hence the complexity of the interface graph is $O(1)$.
  By Theorem \ref{thm:uniqueCases},
  the boundary of any Yin set
  contains only a finite number of cycles.
  Hence the memory cost of all topological data
  is $O(1)$
  for any finite number of Yin sets.
  
  By Definition \ref{def:YinSet},
  the boundary of any Yin set is bounded.
  Furthermore, the boundary curves being piecewise analytic
  implies that the total length of the interface is $O(1)$,
  which, together with the $(r_{\Tiny}, h_L)$-regularity
  and $r_{\Tiny}$ being a constant,
  implies that 
  the total number of markers on the interface is $O(\frac{1}{h_L})$.
  Hence the memory cost of the geometric info
  $S_E$ or $S_{CT}$ is $O(\frac{1}{h_L})$.
  Therefore, the memory usage of the multiphase cubic MARS method
  is always $O(\frac{1}{h_L})$ during the entire simulation.

  In Algorithm \ref{alg:ARMS},
  the CPU time for fitting all cubic splines in $S_{CT}^n$
  is $O(\frac{1}{h_L})$ since the tridiagonal linear system
  in \eqref{eq:miLinearSystem2}
  can be solved by the Thomas algorithm
  with the optimal complexity of $O(\frac{1}{h_L})$.
  The cost of other operations
  in Algorithm \ref{alg:ARMS}
  is linear with respect to the number of markers.
  \qed
\end{proof}

Lemma \ref{lem:complexityOfARMS} states that,
 with respect to the number of interface markers, 
 the complexity of the multiphase cubic MARS method
 is linear and thus optimal.
This optimal complexity also holds
 for VOF methods, level-set methods, and front tracking methods. 
As discussed in the opening paragraph
 of Sec.~\ref{sec:curvatureBasedARMS},
 the interface is a set of codimension 1 
 and one can choose the values of $r_h$ and $\alpha$
 in the relation $h_L=r_h h^{\alpha}$
 so that the cost of IT is only a small fraction of
 that of the main flow solver.
}

\revise{
Under the MARS framework,
 the IT error of the multiphase cubic MARS method 
 is shown to be $E_{\mathrm{IT}}(T)  = {O}(h_L^4) + {O}(k^{4\alpha})$
 for any $T>t_0$,
 where 
 $4\alpha$ is the order of accuracy of the employed ODE solver.
In addition, 
 the multiphase cubic MARS method can be extended
 to the case of flow maps 
 in Definition \ref{def:interfaceTracking} 
 being not homeomorphic, 
 via pinpointing new intersections of edges of the interface graph
 and duly updating the topological data of tracked phases.
We defer details of these investigations to future papers.
}




\section{Tests}
\label{sec:tests}
In this section,
 we perform a variety of benchmark tests
 to demonstrate the power of the proposed method
 in accurately and efficiently tracking multiple materials
 with arbitrarily complex topology and geometry. 
To approximate the exact flow map
 with the fourth-, sixth-, and eighth-order accuracy, 
 we employ the classic fourth-order Runge--Kutta method,
 the explicit one-step method by Verner
 \cite{verner78:_explic_runge_kutta_method_estim},
 and that by Dormand and Prince \cite{dormand81:_high_runge_kutta},
 respectively.
These methods are chosen solely based on ease of implementation,
 and the convergence rates of our method
 would be qualitatively the same 
 if another time integrator of the same order were employed.

\revise{In light of \eqref{eq:metric}},
 the IT error of a phase ${\cal M}_i$
 at time $t_n$ is defined as 
 \begin{equation}
   \label{eq:geomError}
   E^{g}_{i}(t_n) =
   \|\mathcal{M}_i(t_n) \oplus \mathcal{M}^{n}_i\|
   = \sum\nolimits_{{\mathcal C}_{\mathbf{j}}\subset \Omega}
   \left\|\left(\mathcal{M}_i(t_n)\cap{\mathcal C}_{\mathbf{j}}\right)
     \oplus
     \left(\mathcal{M}^n_i\cap{\mathcal C}_{\mathbf{j}}\right)
   \right\|,
 \end{equation}
 where
 ${\mathcal M}^{n}_i$,
 the computational result from $\widetilde{\Gamma}^n$, 
 approximates the exact result ${\mathcal M}_i(t_n)$ 
 and ${\cal C}_{\mathbf{j}}$'s are the control volumes 
 that partition the computational domain $\Omega$, 
 i.e., $\bigcup^{\perp\perp}_{\mathbf{j}}{\cal C}_{\mathbf{j}}=\Omega$
 and 
 $\mathbf{i}\ne\mathbf{j}
   \implies {\cal C}_{\mathbf{i}}\cap{\cal C}_{\mathbf{j}}=\emptyset$.

As $h\to 0$, 
 the computation of symmetric differences in (\ref{eq:geomError})
 tends to be more and more ill-conditioned, 
 \revise{since the computed and exact arcs
   become more and more parallel}.
Hence in practice we approximate $E^g_i$ with
\begin{equation}
  \label{eq:errorIT}
  E_i(t_n) := \sum\nolimits_{\mathcal{C}_{\mathbf{j}} \subset \Omega} \big|
  \|\mathcal{M}_i(t_n) \cap \mathcal{C}_{\mathbf{j}}\|
  - \|\mathcal{M}_i^n \cap \mathcal{C}_{\mathbf{j}}\|
  \big|.
\end{equation}
The \emph{total IT error} of $N_p$ phases
 is defined as $\sum_{i=1}^{N_p}E_i$. 

\subsection{Vortex shear of a quartered circular disk}
\label{sec:testVortex}

 \begin{figure}
  \centering
  \subfloat[$t = 0$]{
    \includegraphics[width=0.28\linewidth]{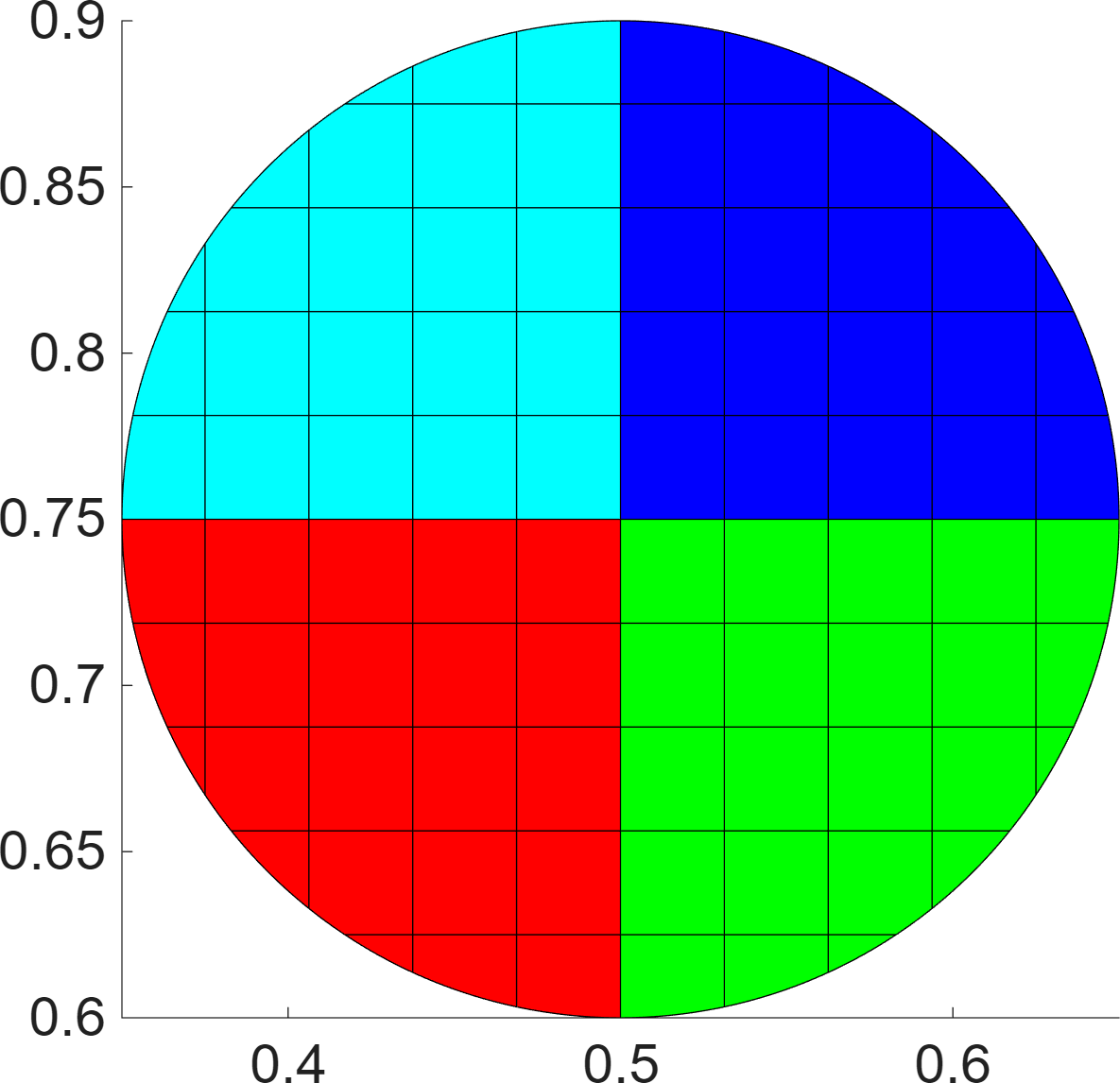}
  }
  \hfill
  \subfloat[$T=4$; $t = \frac{1}{2}T$]{
    \includegraphics[width=0.28\linewidth]{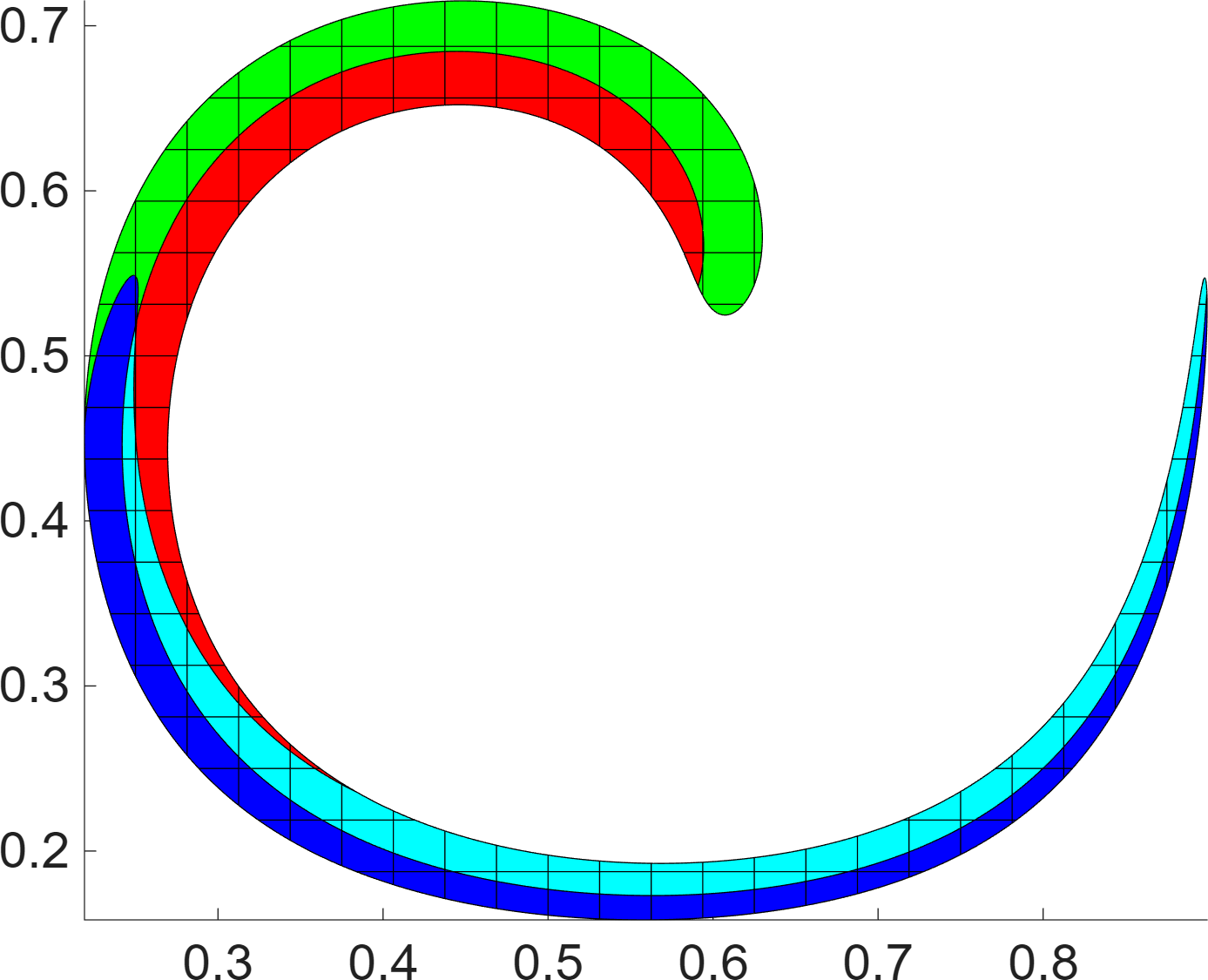}
  }
  \hfill
  \subfloat[$T=4$; $t = T$]{
    \includegraphics[width=0.28\linewidth]{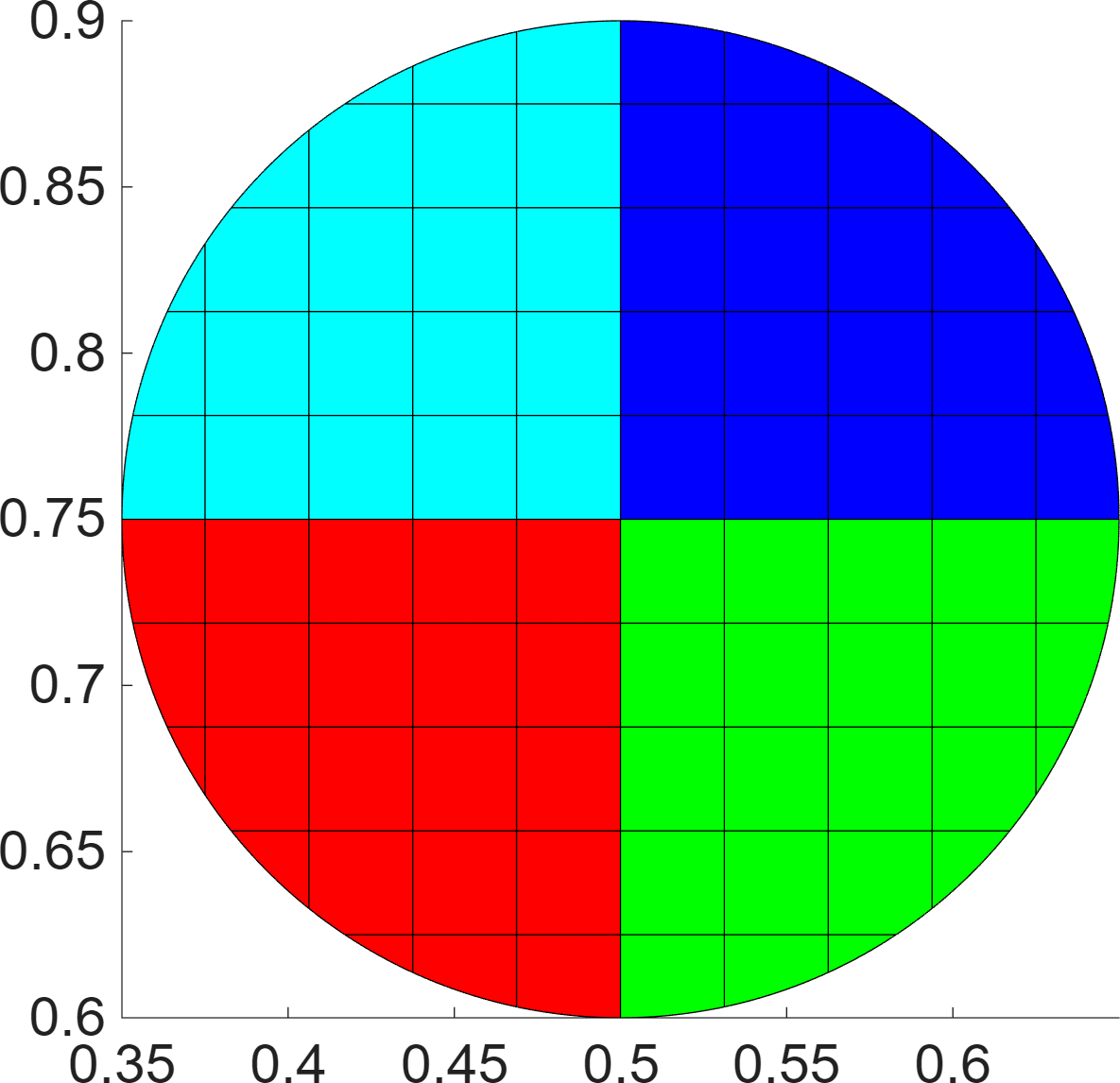}
  }

  \subfloat[$T=8$; $t = \frac{1}{4}T$]{
    \includegraphics[width=0.28\linewidth]{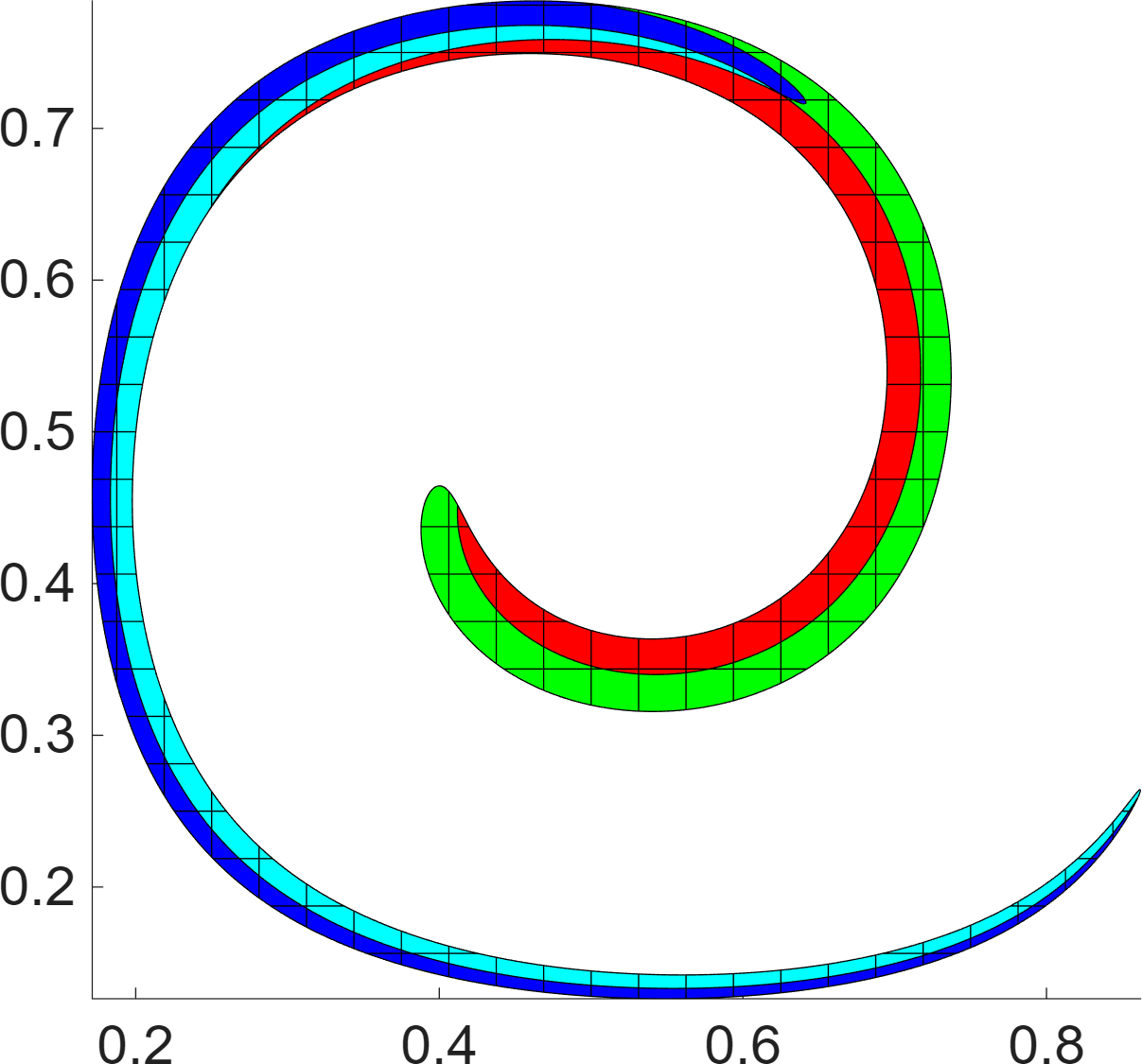}
  }
  \hfill
  \subfloat[$T=8$; $t = \frac{1}{2}T$]{
    \includegraphics[width=0.28\linewidth]{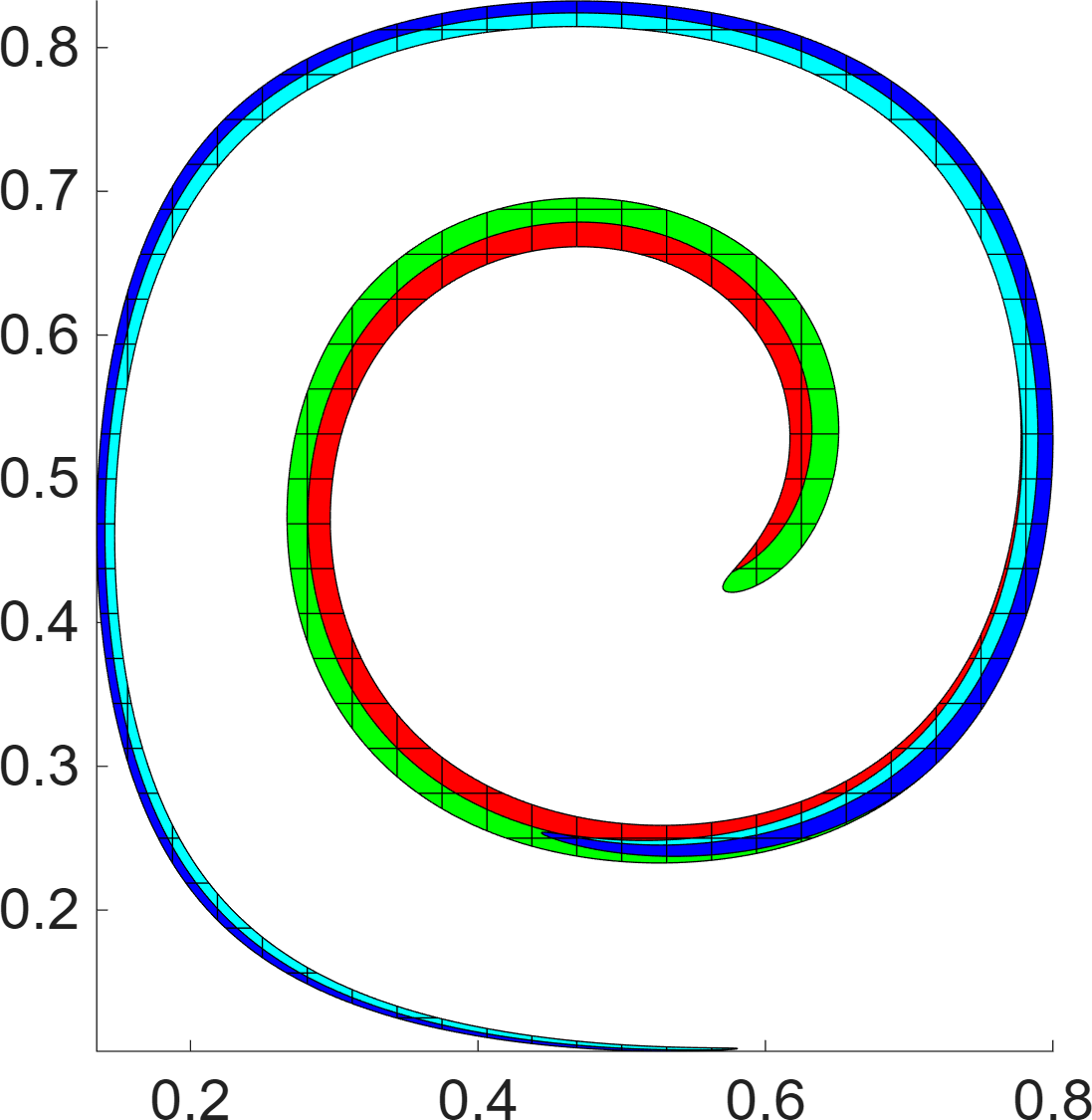}
  }
  \hfill
  \subfloat[$T=8$; $t = T$]{
    \includegraphics[width=0.28\linewidth]{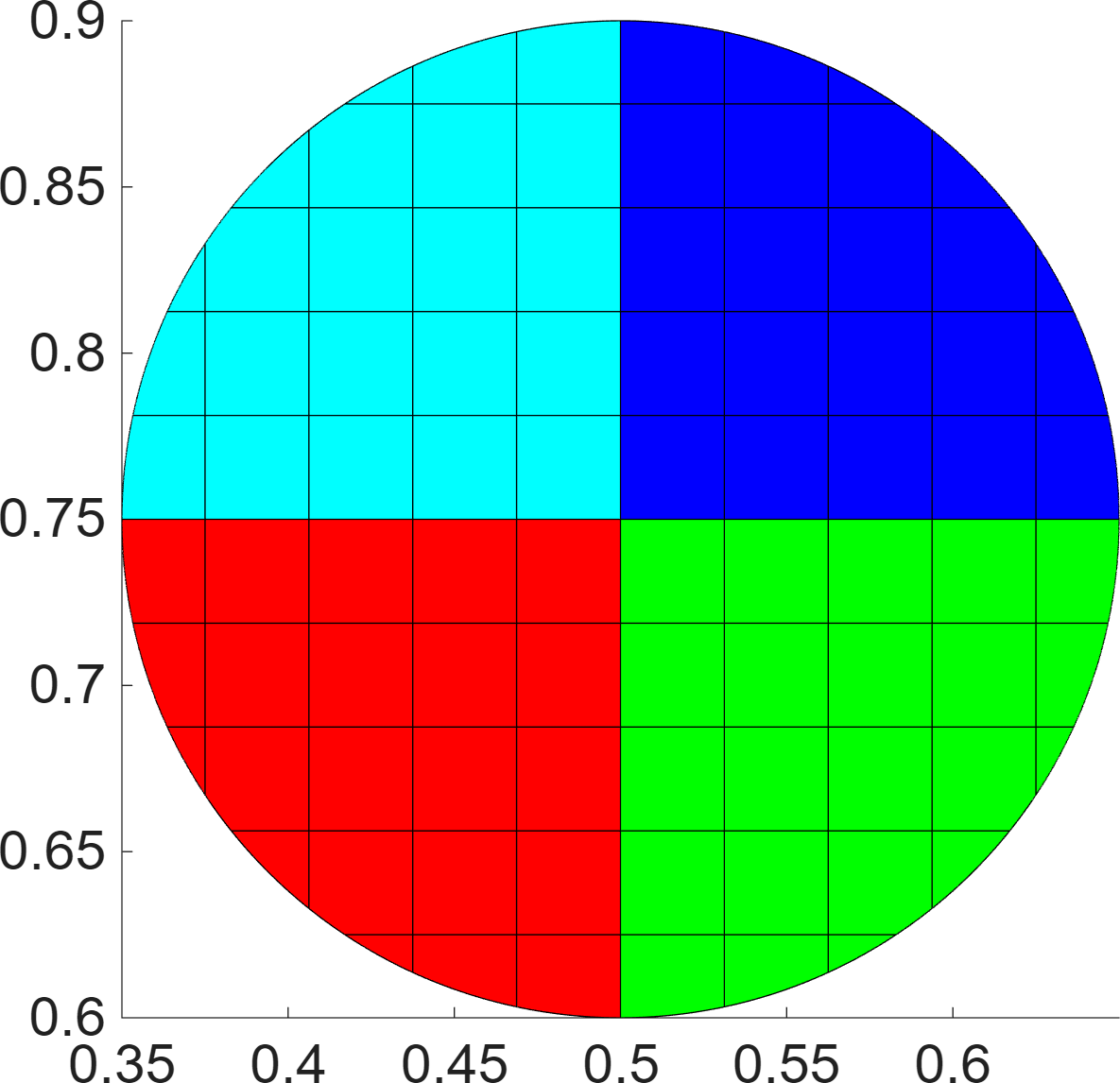}
  }

  \subfloat[$T=12$; $t = \frac{1}{4}T$]{
    \includegraphics[width=0.28\linewidth]{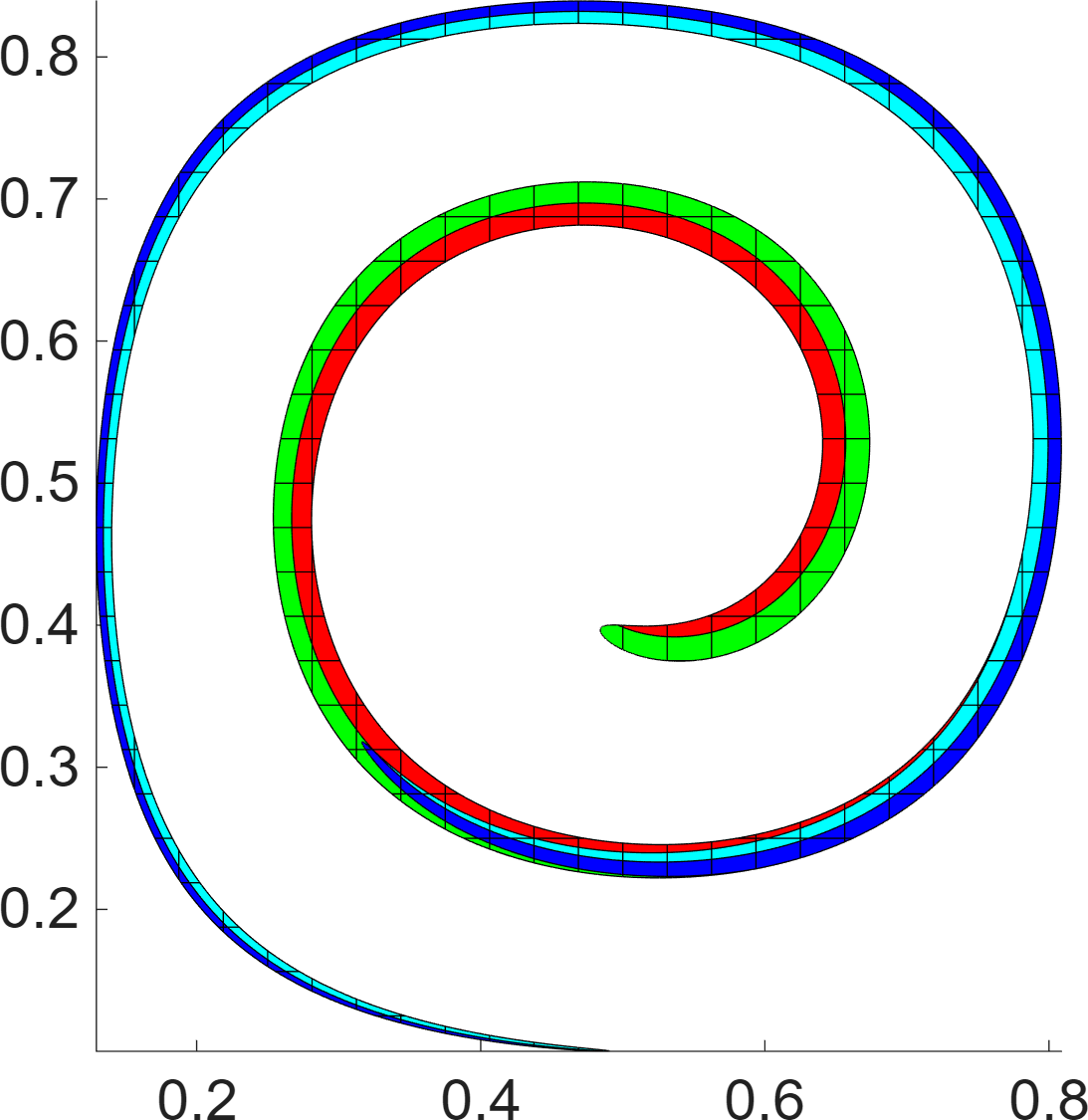}
  }
  \hfill
  \subfloat[$T=12$; $t = \frac{1}{2}T$]{
    \includegraphics[width=0.28\linewidth]{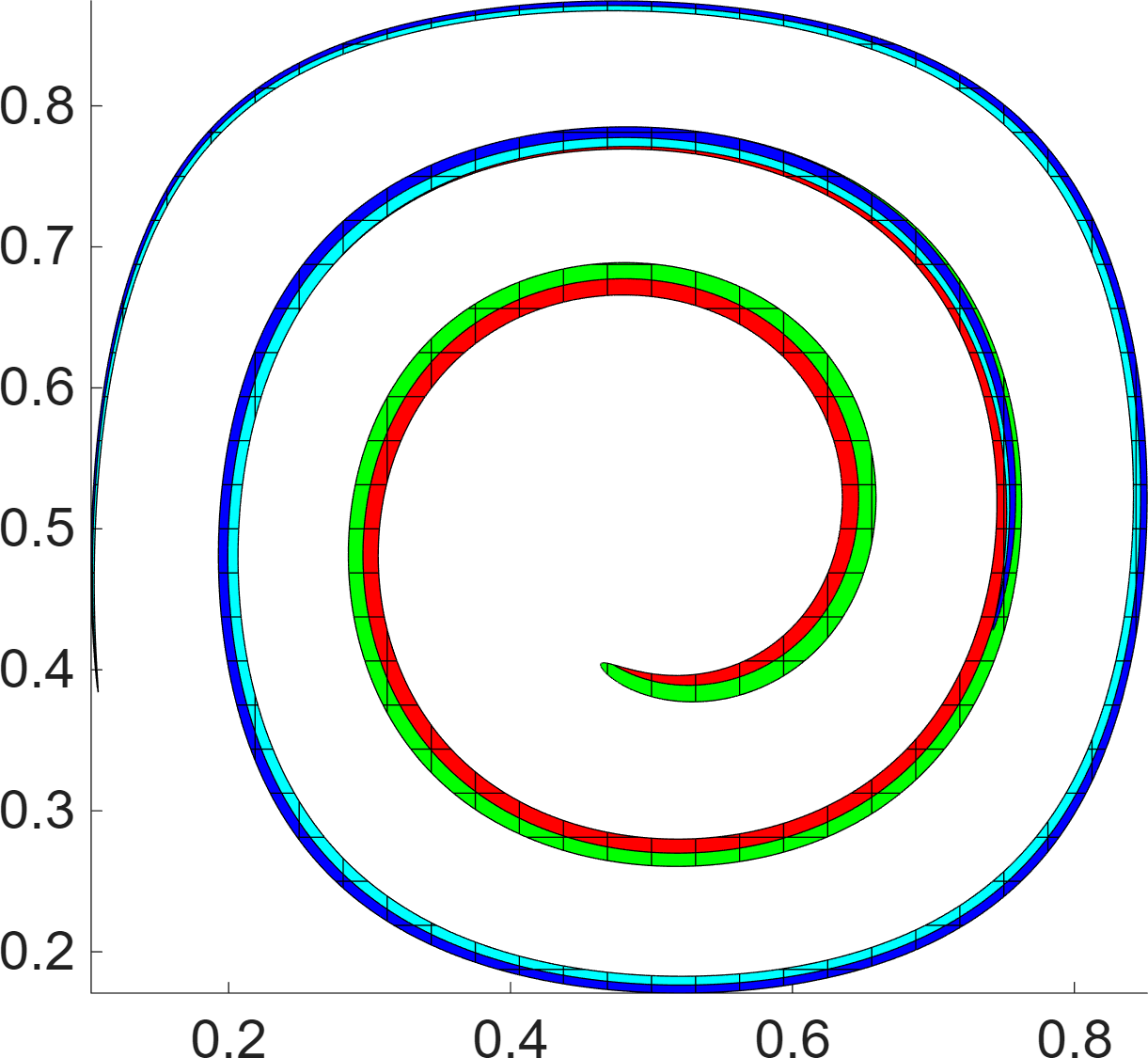}
  }
  \hfill
  \subfloat[$T=12$; $t = T$]{
    \includegraphics[width=0.28\linewidth]{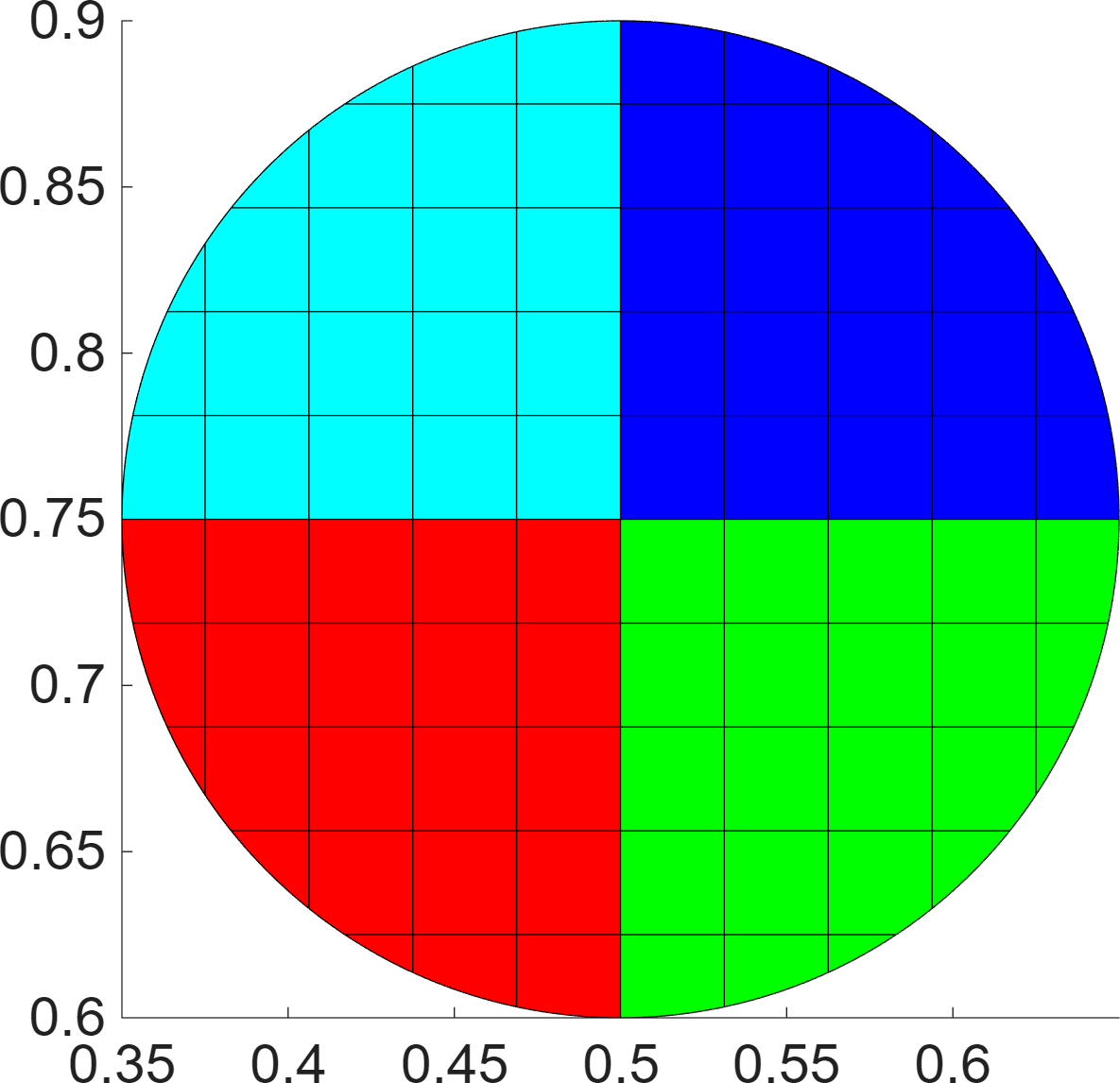}
  }

  \subfloat[$T=16$; $t = \frac{1}{4}T$]{
    \includegraphics[width=0.28\linewidth]{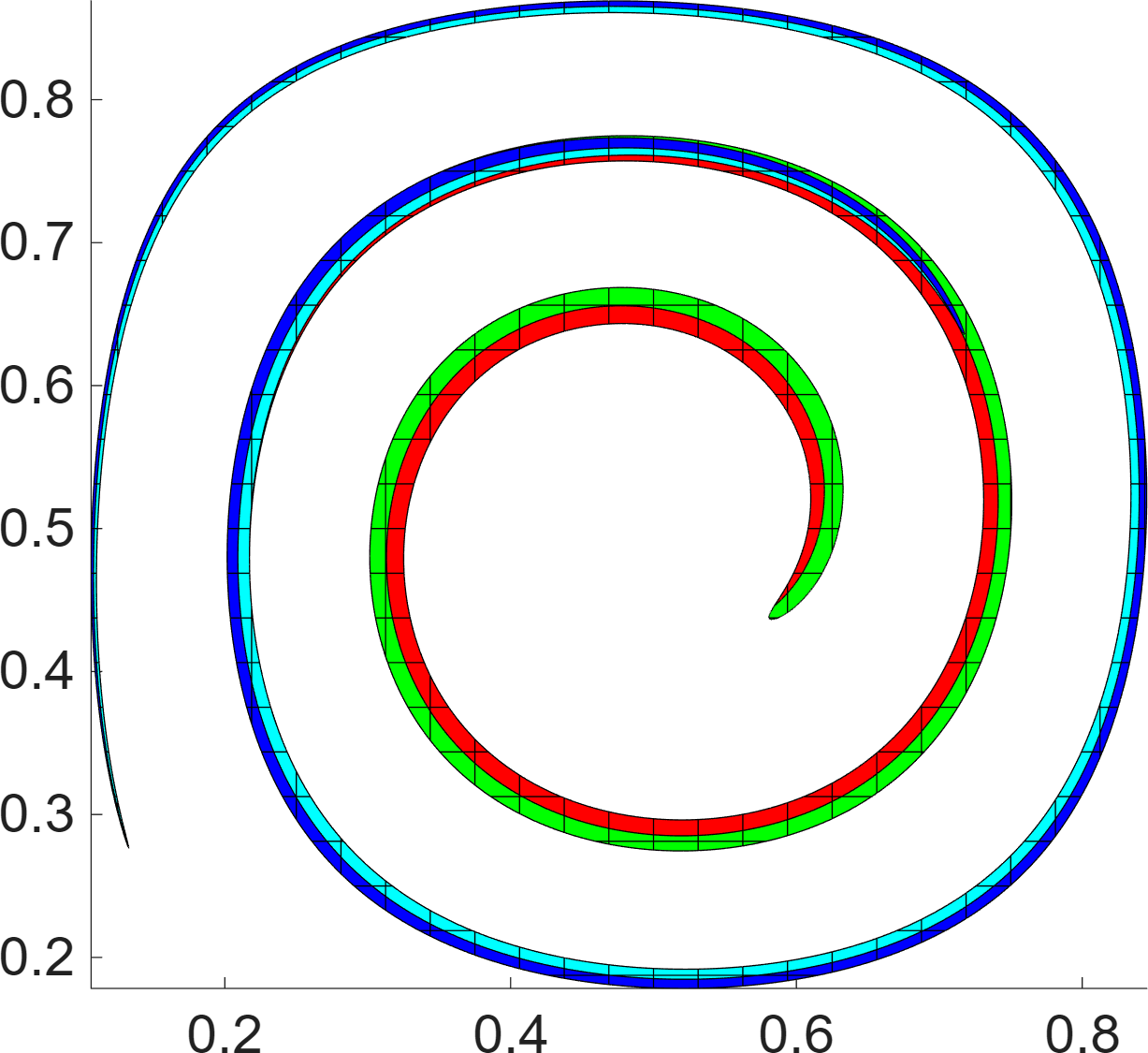}
  }
  \hfill
  \subfloat[$T=16$; $t = \frac{1}{2}T$]{
    \includegraphics[width=0.28\linewidth]{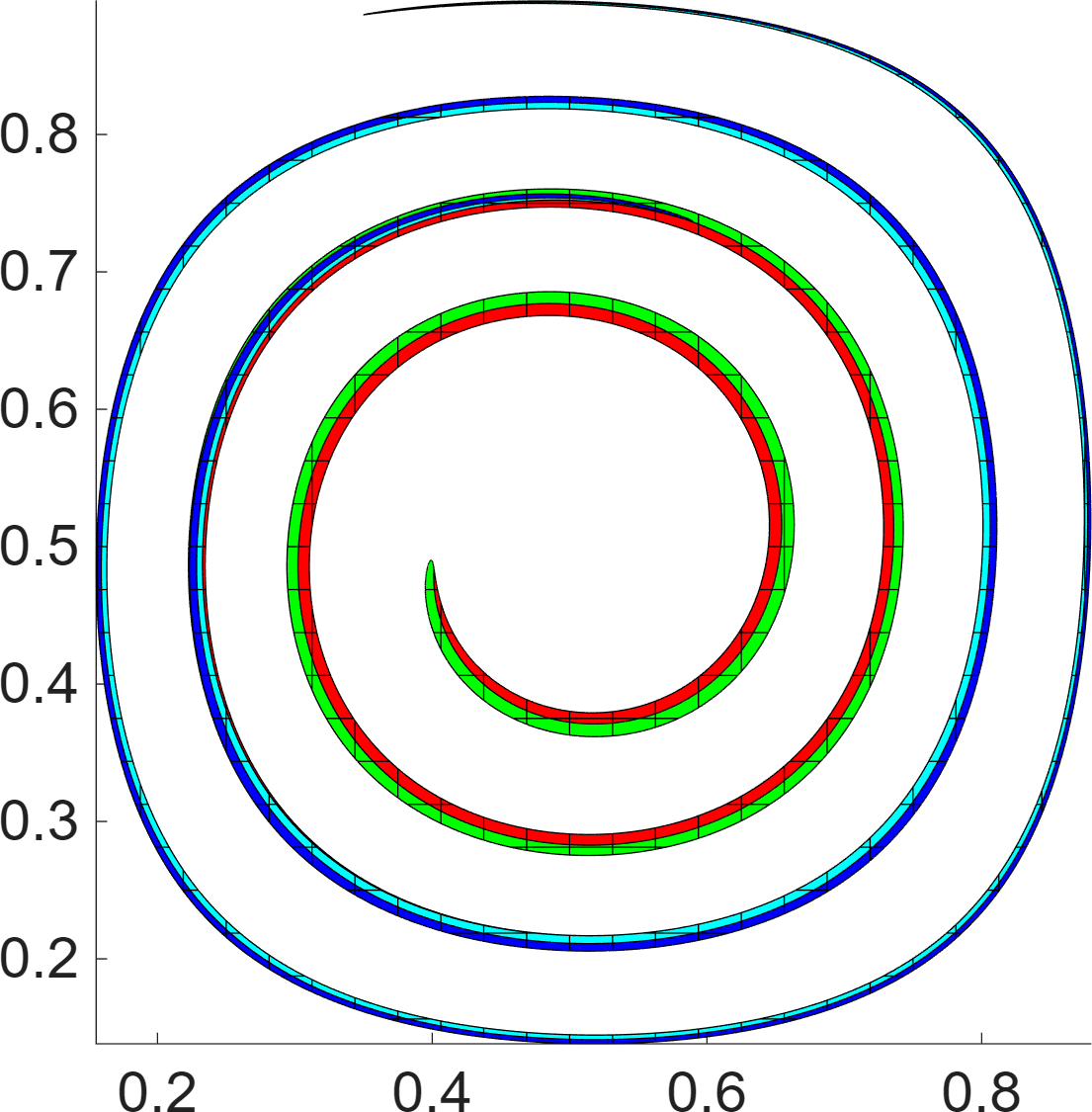}
  }
  \hfill
  \subfloat[$T=16$; $t = T$]{
    \includegraphics[width=0.28\linewidth]{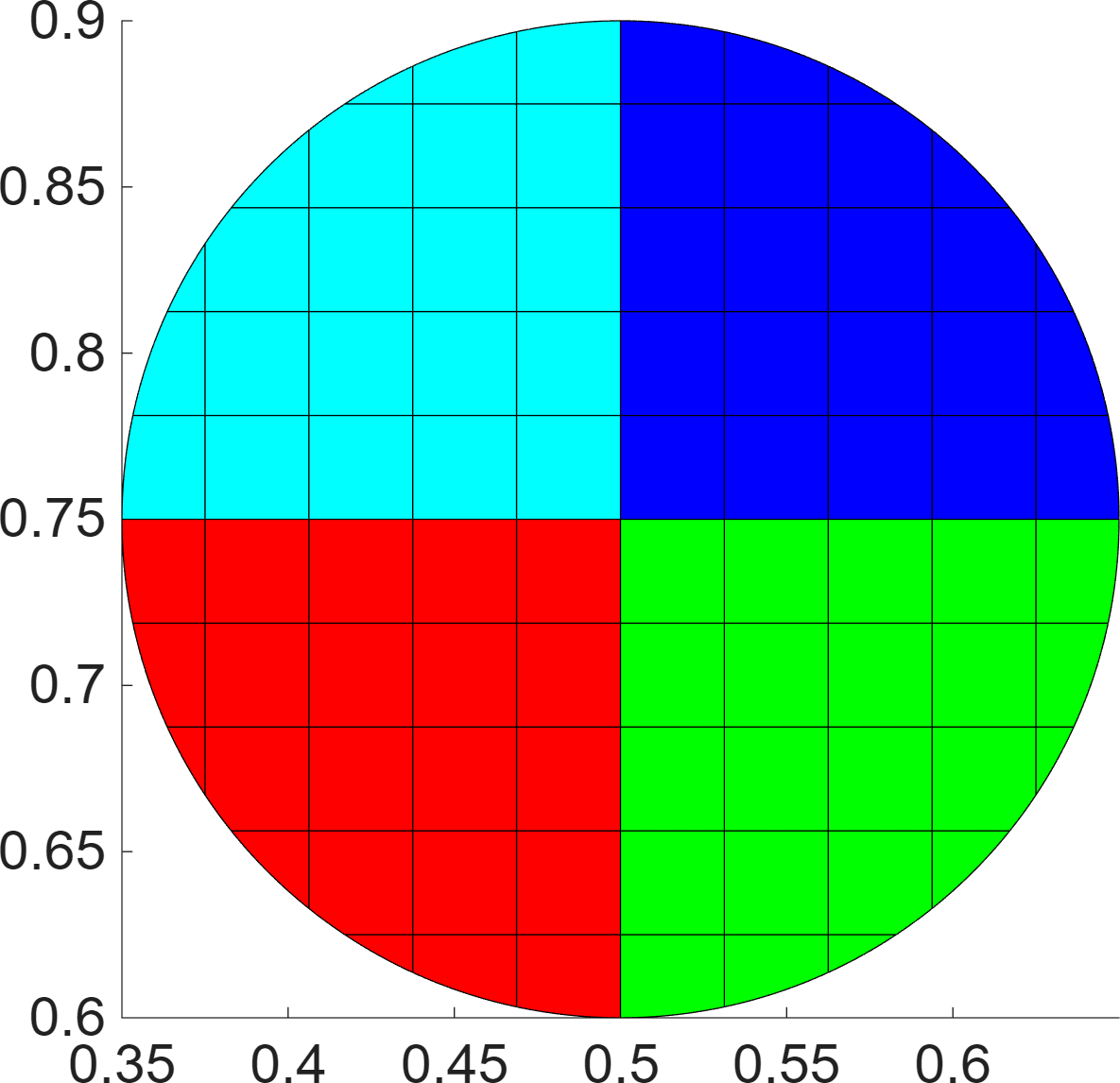}
  }
  \caption{Solutions of the cubic MARS method
    for the vortex shear test with $T = 4, 8, 12, 16$
    on the Eulerian grid of $h = \frac{1}{32}$.
    The initial distances between adjacent markers
    for $T=4,8$ and $T=12,16$
    are respectively set to the uniform constant $0.1 h$
    and the varying value $\frac{1}{2}h_L(\rho)$
    with $h_L$ defined in (\ref{eq:curvature-based-arms}).
    See Table \ref{tab:test-quartered-circular-vortex1}
    for values of other parameters.
  }
  \label{fig:test-quartered-circular-disk-T4}
\end{figure}

\begin{table}
  \centering
    \caption{IT Errors and convergence rates
      of the proposed cubic MARS method with $k = \frac{1}{8}h$
      in solving the vortex shear test of $T = 4, 8, 12, 16$.
      The first part is based on the total IT error 
      while the second on IT errors of individual phases,
      cf (\ref{eq:errorIT}).
    }
    \small
    \renewcommand{\arraystretch}{1.2}
    \input{{\TABDIR}table6.1}
    \label{tab:test-quartered-circular-vortex1}
\end{table}

\begin{table}
\centering
\caption{Total IT errors and convergence rates
  of the multiphase cubic MARS methods
  compared with those of VOF/MOF methods
  in solving three vortex-shear tests.
  For test (a),
  the last three lines are taken from \cite[Tab. 7]{sijoy2010volume},
  where the two-letter acronyms
  LV, NI, MC, and MB
  stand for the LVIRA algorithm
  \cite{puckett1991:_volum_inter_track_algor_applic},
  de Niem's intersection check method \cite{deniem2007volume},
  Mosso and Clancy's method \cite{mosso1995geometrically},
  and a combination of MC \cite{mosso1995geometrically}
  and Benson's method \cite{benson1998eulerian},
  respectively.
  The last two lines for tests (b) and (c)
  are taken from \cite[Tab. 2]{hergibo2023moment}
  and \cite[Tab. 4]{hergibo2023moment}, respectively.
  For all MARS methods,
  we use $h_L=0.2h$ or $h_L^c=0.2h$;
  see Table \ref{tab:test-quartered-circular-vortex1}
  for values of other parameters.
}
\small
\renewcommand{\arraystretch}{1.2} 
\input{{\TABDIR}table6.2}
\label{tab:comparisonVortex}
\end{table}

Referring to Definition \ref{def:interfaceTracking},
 the flow map of this test is that of the ODE
 $\frac{\dif \mathbf{X}}{\dif t}=\mathbf{u}(\mathbf{X},t)$
 with $\mathbf{u}=(\frac{\partial \psi}{\partial y},
 -\frac{\partial \psi}{\partial x})$
 determined from the stream function
 \begin{equation}
   \label{eq:vortex-velocity}
  \begin{array}{l}
   \psi(x,y) = - \frac{1}{\pi} \sin^2(\pi x) \sin^2(\pi y)
   \cos \left( \frac{\pi t}{T} \right),
  \end{array}
\end{equation}
 where the time period $T=4,8,12,16$.
At time $t=\frac{T}{2}$,
 the velocity field is reversed by the cosinusoidal temporal factor
 so that the exact solution $({\cal M}_i(t))_{i=1}^{5}$ at $t=T$
 is the same as the initial condition
 $({\cal M}_i(t_0))_{i=1}^{5}$ at $t_0=0$.
As shown in Fig.~\ref{fig:test-quartered-circular-disk-T4}(a),
 the four colored Yin sets 
 constitute a circular disk
 with its radius as $0.15$ and its center at $[0.5, 0.75]^T$
 while the last Yin set is 
 the unbounded complement of the circle.

For this IT problem of five phases, 
 the cases \mbox{$T=4,8$} are solved
 by the multiphase cubic MARS method
 with the constant ARMS strategy 
 while the cases \mbox{$T=12,16$} by that with
 the curvature-based ARMS strategy (\ref{eq:curvature-based-arms}) and 
 \begin{equation}
   \label{eq:vortex-test-parametersCB}
   \begin{array}{l}
   r_{\Tiny} = 0.1;  \
   (\rho_{\min}^c, \rho_{\max}^c) = (10^{-5},0.2); \
   r_{\min}^c = 0.01, 0.005;\
   \sigma^c(x) = x.
   \end{array}
 \end{equation}

The time step sizes are set to $k = \frac{1}{8}h$
 for $T=4,8,12,16$,
 so that IT errors are dominated
 not by temporal discretizations of flow maps
 but by spatial approximations of the interface;
 otherwise the temporal symmetry in (\ref{eq:vortex-velocity})
 would lead to convergence rates higher than expected,
 such as those in Table \ref{tab:compare-result}. 
In (\ref{eq:vortex-test-parametersCB}),
 we choose $r_{\Tiny} = 0.1$ to limit to one order of magnitude
 the variation of chordal lengths
 caused by the \emph{random} tangential advection of markers. 
The values $(\rho_{\min}^c, \rho_{\max}^c) = (10^{-5},0.2)$
 come from the presence of line segments
 in the initial condition
 and the fact of the radius of the initial circle being 0.15. 
As $T$ increases from 12 to 16,
 we decrease $r_{\min}^c$ from $0.01$ to $0.005$
 to account for the larger deformation.
 

In Fig.~\ref{fig:test-quartered-circular-disk-T4}, 
 we plot our solutions on the Eulerian grid of $h=\frac{1}{32}$
 at key time instances. 
During the entire simulation,
 the interface graph
 $G_{\Gamma}=(V_{\Gamma},E_{\Gamma},\psi_{\Gamma})$
 that represents
 the initial topology of the five phases
 remains the same: 
 $E_{\Gamma}$ always consists of the eight edges
 that connect the five vertices in $V_{\Gamma}=J_{\Gamma}$,
 i.e., the four T junctions on the circle
 and the X junction inside the circle. 
By Algorithm \ref{alg:partitionGraph},
 $C_S$ contains only a single circuit of the four T junctions
 while $T_S$ has two trails
 that correspond to the two disk diameters.
The constancy of these topological data confirms
 the validity and efficiency of separating topology from geometry.
At each time step, the geometry is evolved
 by updating the splines in $S_{CT}^n$. 
Despite the enormous deformations
 and the large size of the Eulerian grid, 
 each phase remains connected without
 generating any flotsam
 and the difference between the final solution
 and the initial condition is indiscernible.

To visually compare results of MARS and VOF methods
 reviewed in Sec.~\ref{sec:intro},
 we note that the vortex shear test with $T=4$
 is the same as that in \cite[Sec. 5.5]{schofield2009second}.
Thus Fig.~\ref{fig:test-quartered-circular-disk-T4}(a,b,c)
 compares directly to \cite[Fig. 17]{schofield2009second},
 where
 both the material-order-dependent Young's method \cite{youngs1982:_time}
 and the material-order-independent
 power diagram method \cite{schofield2009second} %
 generate flotsam, 
 failing to preserve the connectedness of the deforming phases.
Also shown in \cite[Fig. 17]{schofield2009second}
 are the prominently different geometric features
 between the final solutions 
 and the initial conditions of these VOF methods.

In Table \ref{tab:test-quartered-circular-vortex1}
 we present, for all cases of $T=4,8,12,16$,
 total IT errors $\sum_{i=1}^{5} E_i$
 and convergence rates of MARS methods 
 with both constant and curvature-based ARMS strategies.
For $T=16$ and $h_L^c=0.2h$,
we also show, for each individual phase,
results based on $E_i$ in (\ref{eq:errorIT}).
In all cases and for all phases, 
 fourth-, sixth-, and eighth-order convergence rates
 are clearly demonstrated for the choices of $h_L$ or $h_L^c$
 being $O(h)$, $O(h^{\frac{3}{2}})$, and $O(h^2)$,
 respectively. 
The smallest total errors $2.87\times 10^{-13}$ and $2.17\times 10^{-13}$ 
 indicate excellent conditioning of both ARMS strategies.

To quantitatively compare MARS with VOF/MOF methods,
 we first quote from \cite[p. 744]{schofield2009second} that, 
 for the test of $T=4$ with $h=\frac{1}{64}$ and $k=\frac{1}{8}h$, 
 the smallest IT errors of Young's method \cite{youngs1982:_time}
 and the power diagram method \cite{schofield2009second}
 are respectively $1.28\times 10^{-3}$
 and $1.35\times 10^{-4}$,
 which are much larger than $2.53\times 10^{-9}$,
 the total IT error of MARS
 in the case of $h=\frac{1}{32}$, $h_L=0.2h$, and $k=\frac{1}{8}h$ 
 in Table \ref{tab:test-quartered-circular-vortex1}. 
Then in Table \ref{tab:comparisonVortex}
 we compare our cubic MARS method with VOF/MOF methods
 for solving three other vortex-shear tests in the literature. 
The circular disk in test (a)
 consists of three phases with two triple points
 while tests (b) and (c) are the classic two-phase test
 with $T=8$ and 12, respectively.
For all tests,
 the proposed cubic MARS method is more accurate than VOF/MOF methods
 by many orders of magnitude.

\begin{figure}
  \centering
  \subfloat[$t = 0$]{
    \includegraphics[width=0.31\linewidth]{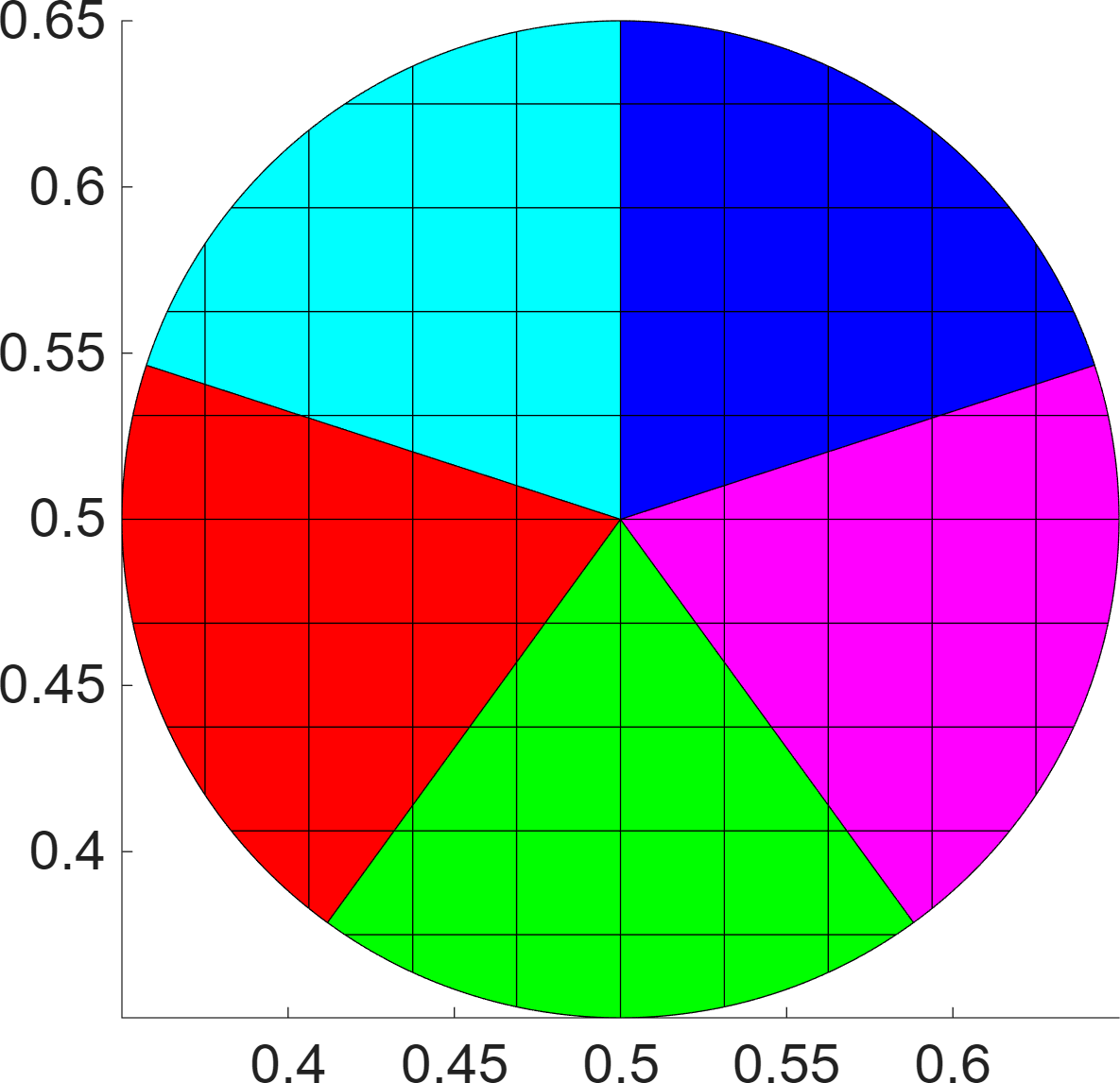}
  }
  \hfill
  \subfloat[$T = 2; t = \frac{1}{2}T$]{
    \includegraphics[width=0.31\linewidth]{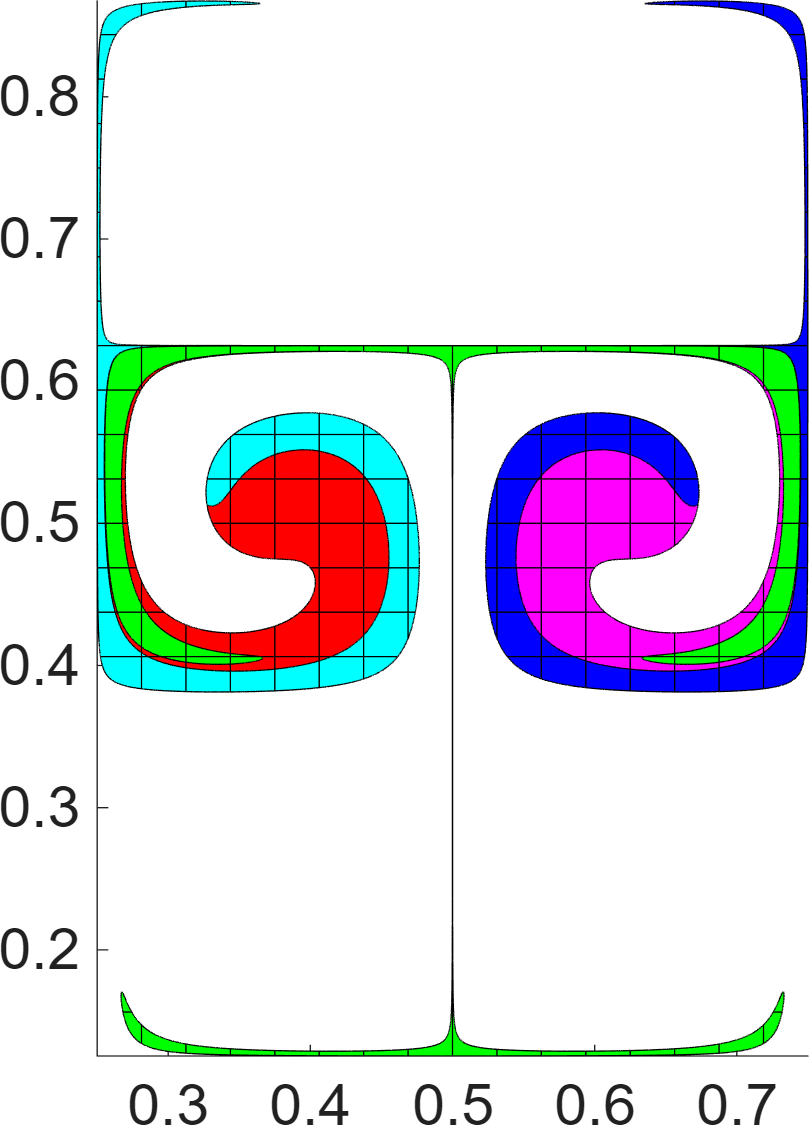}
  }
  \hfill
  \subfloat[$T = 2; t = T$]{
    \includegraphics[width=0.31\linewidth]{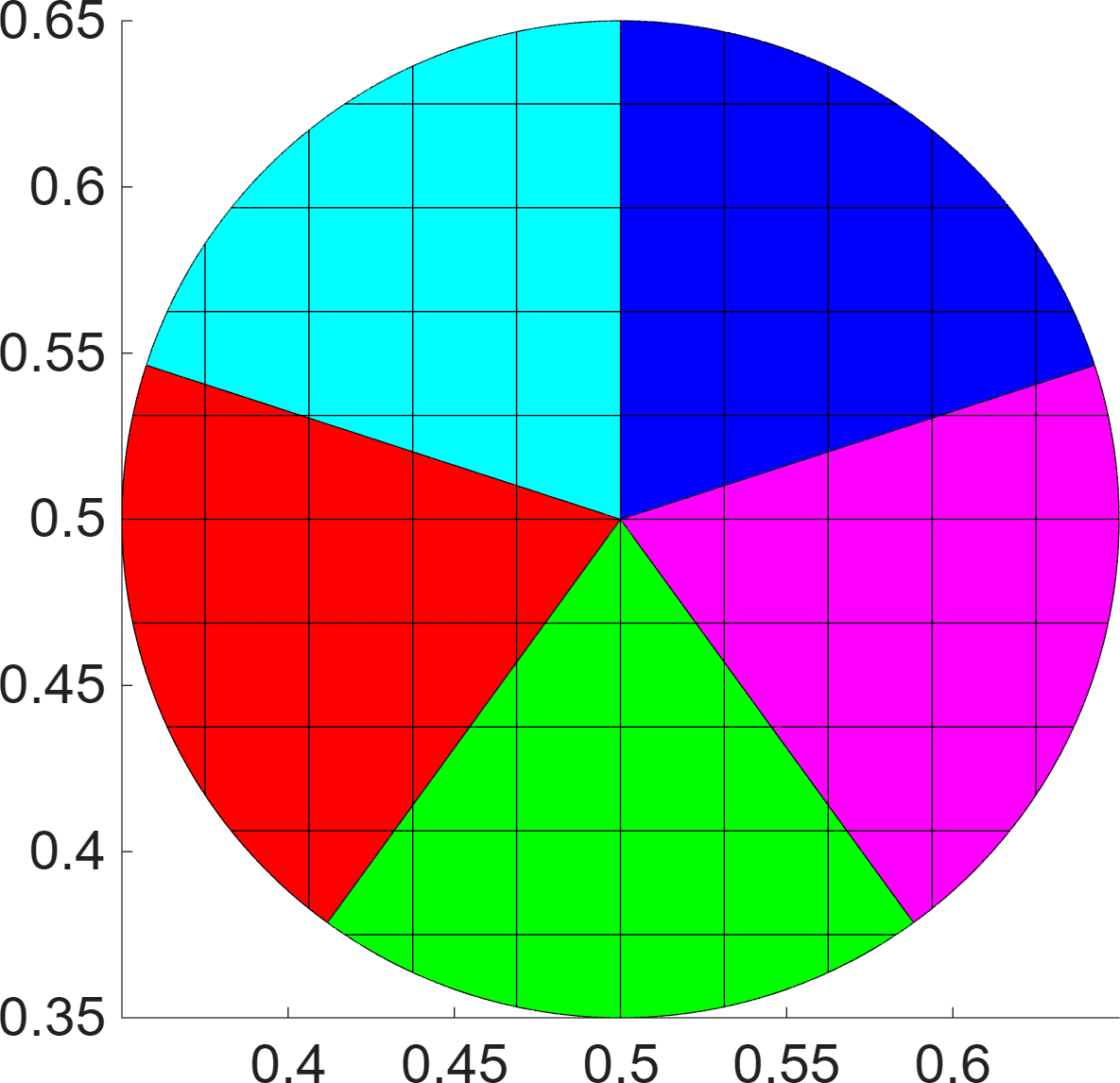}
  }

  \subfloat[$T = 4; t = \frac{1}{4}T$]{
    \includegraphics[width=0.31\linewidth]{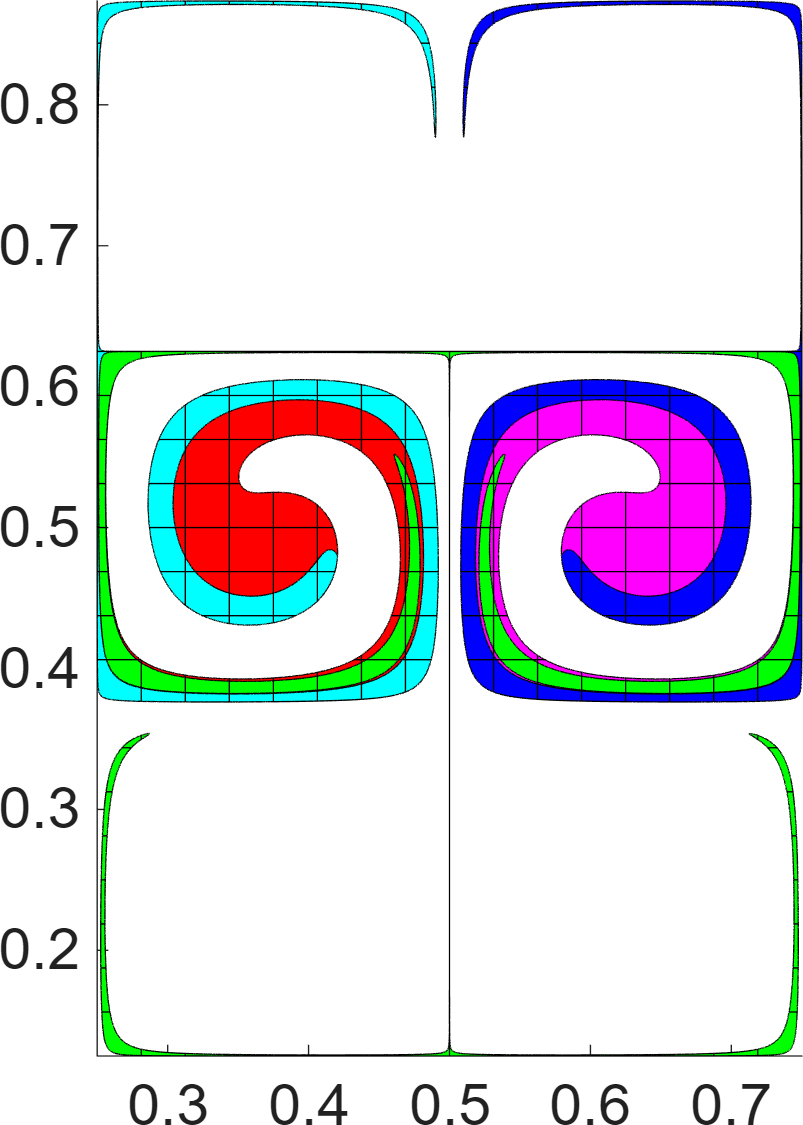}
  }
  \hfill
  \subfloat[$T = 4; t = \frac{1}{2}T$]{
    \includegraphics[width=0.31\linewidth]{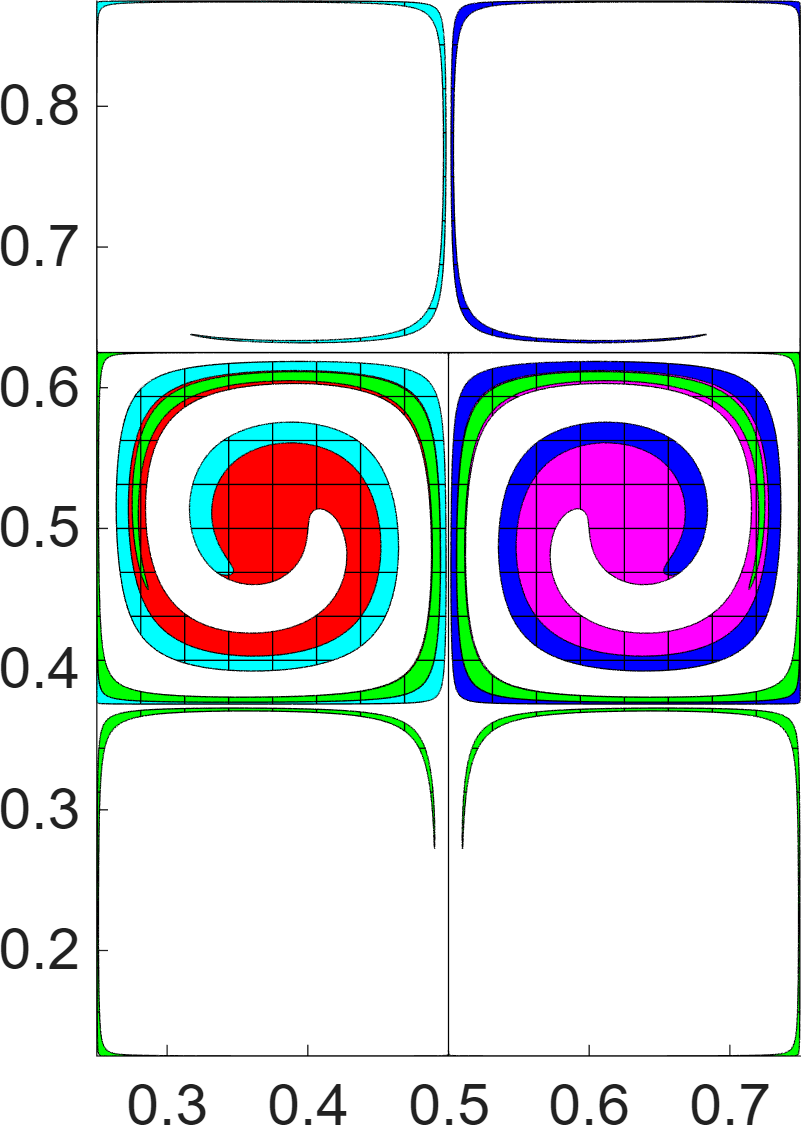}
  }
  \hfill
  \subfloat[$T = 4; t = T$]{
    \includegraphics[width=0.31\linewidth]{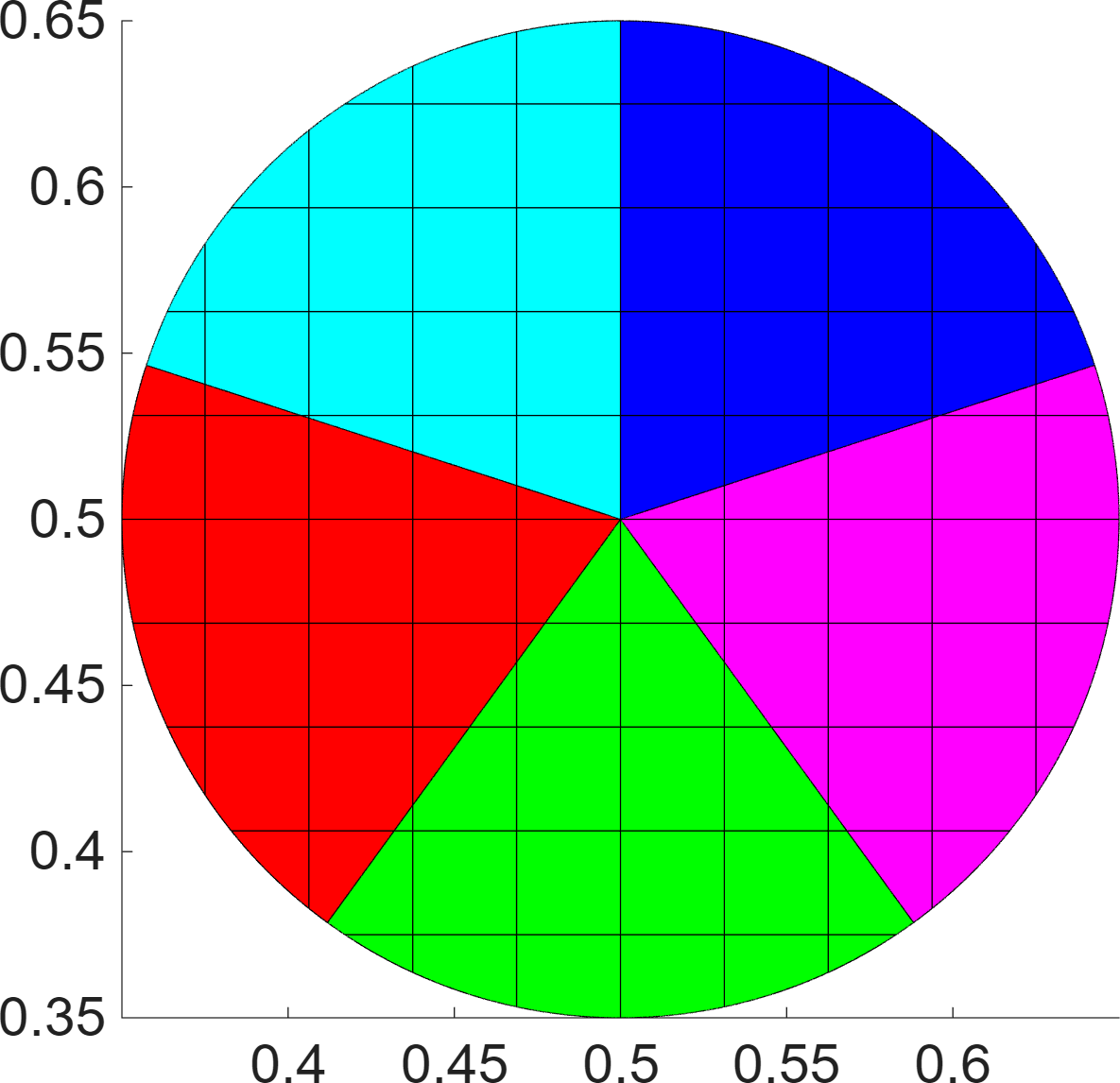}
  }
  
  \subfloat[$T=2$; normalized length of $\partial {\cal M}_i(\frac{t}{T})$
  ]{
    \includegraphics[width=0.222\linewidth]{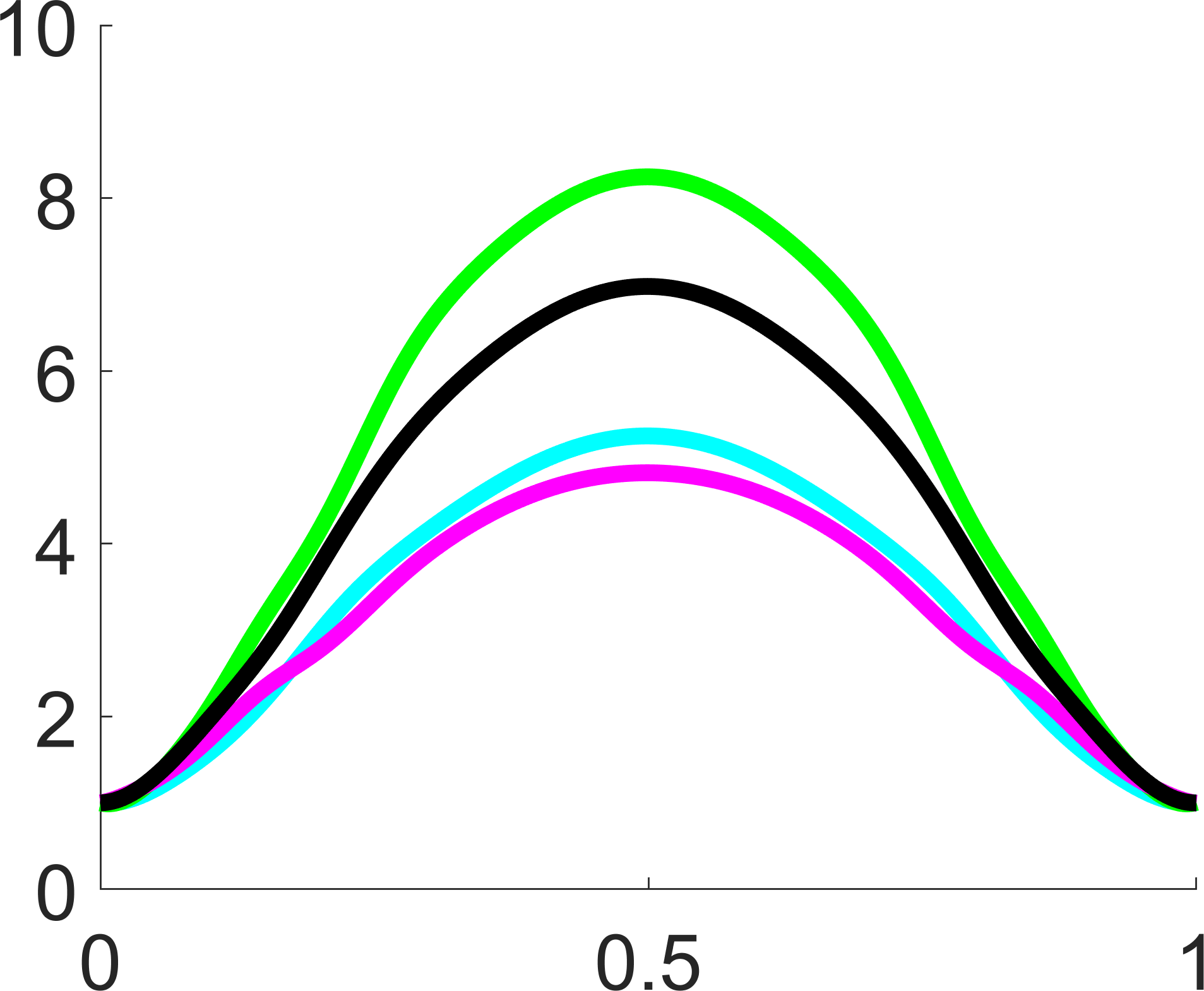}
  }
  \hfill
  \subfloat[$T=2$; normalized marker number in $\partial {\cal M}_i^n$
  ]{
    \includegraphics[width=0.222\linewidth]{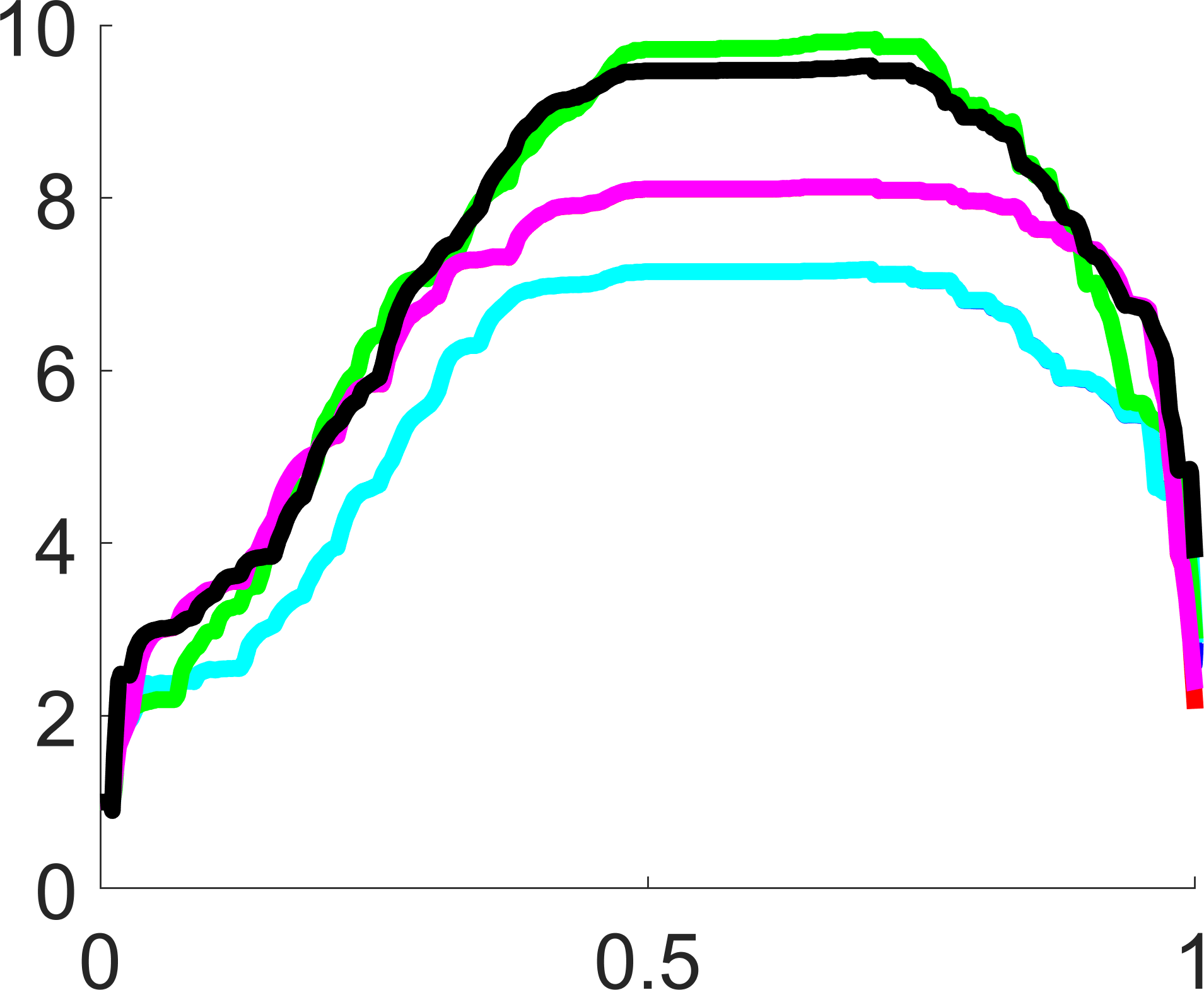}
  }
  \hfill
  \subfloat[$T=4$; normalized length of $\partial {\cal M}_i(\frac{t}{T})$]{
    \includegraphics[width=0.222\linewidth]{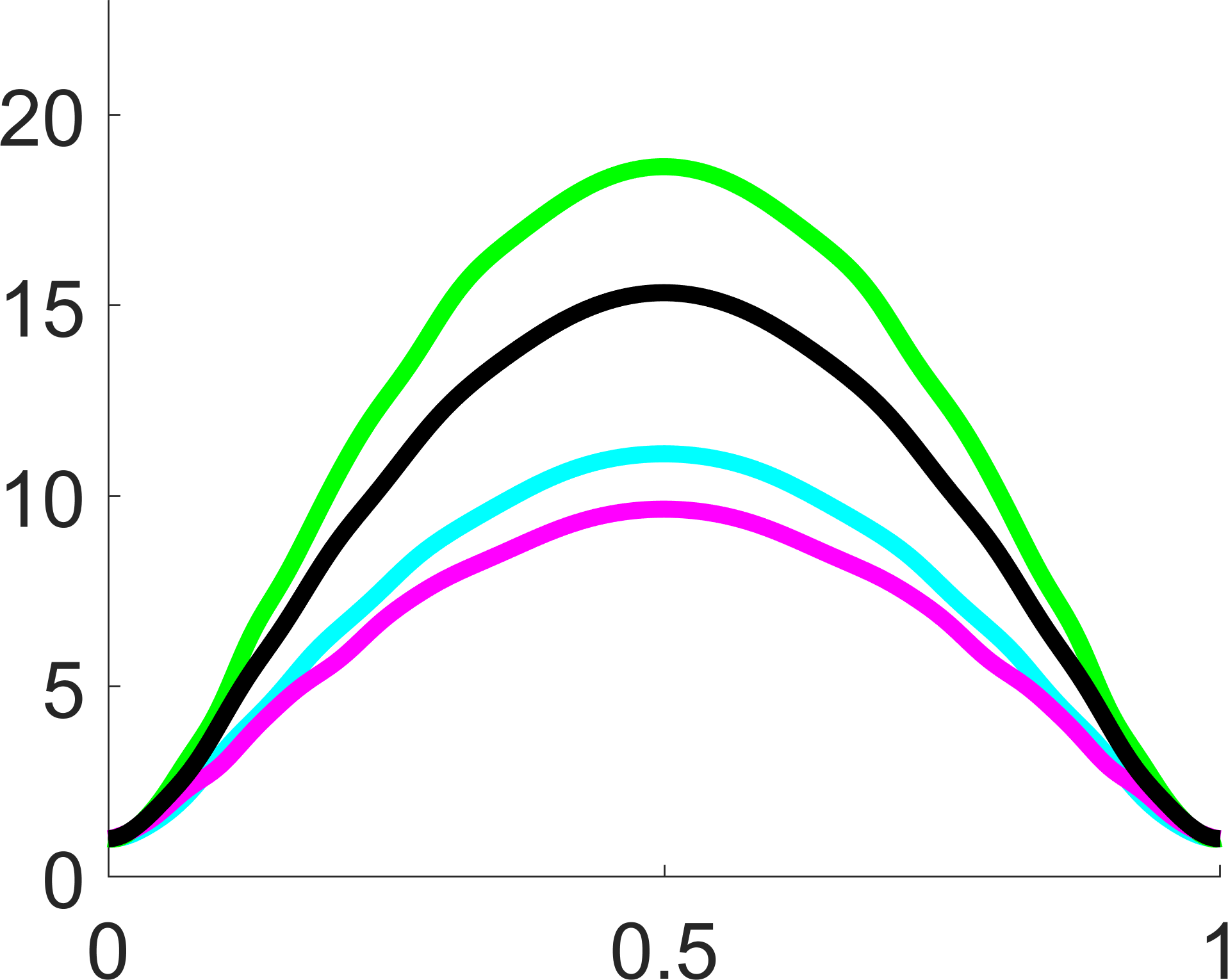}
  }
  \hfill
  \subfloat[$T=4$; normalized marker number in $\partial {\cal
    M}_i^n$
  ]{
    \includegraphics[width=0.222\linewidth]{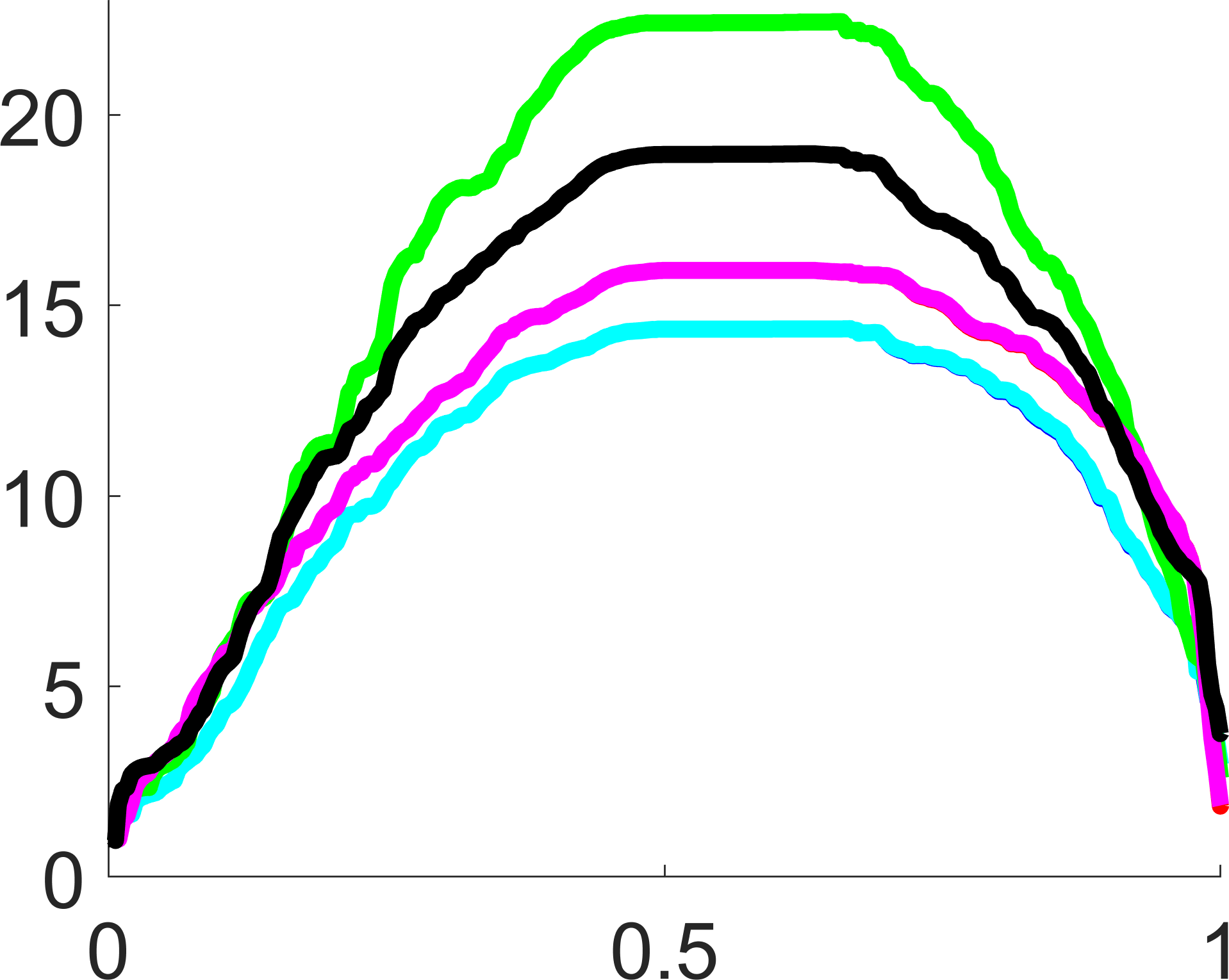}
  }
  \caption{Solutions of the cubic MARS method with $h_L^c = 0.2 h$
    for the deformation test
    of $T = 2$ and $T = 4$ on the Eulerian grid of $h = \frac{1}{32}$.
    Subplots (a)--(f) are snapshots of the solution
    at key time instants. 
    In subplots (g)--(j),
    each phase is represented by a curve of the same color
    except that the white phase (the unbounded complement of the circle)
    is represented by the black curve.
    Due to symmetry,
    the red and blue curves may not be visible.
    The initial distances between markers are $\frac{1}{2}h_L(\rho)$,
    where $h_L$ is defined in (\ref{eq:curvature-based-arms}).
    See Table \ref{tab:deformation-curvature-arms}
    for values of other parameters.
  }
  \label{fig:test-deformation}
\end{figure}

\subsection{Deformation of a circular disk divided into five phases}
\label{sec:testDeformation}

  


The flow map of this test is determined
 in the same way as that in Sec.~\ref{sec:testVortex}, 
 with the stream function as
\begin{equation}
  \label{eq:deformation-velocity}
  \begin{array}{l}
  \psi(x,y) = - \frac{1}{n_{\mathrm{v}}\pi}
  \sin \big(n_{\mathrm{v}} \pi (x + 0.5)\big)
  \cos \big(n_{\mathrm{v}} \pi (y + 0.5)\big)
  \cos \left( \frac{\pi t}{T} \right),
  \end{array}
\end{equation}
where $T = 2,4$ and the number of vortices is $n_{\mathrm{v}} = 4$.
At $t = \frac{T}{2}$, 
 the temporal factor reverses the velocity field so that 
 the exact solution $(\mathcal{M}_i(t))_{i = 1}^6$ at \mbox{$t = T$}
 is identical to the initial condition
 $(\mathcal{M}_i(t_0))_{i = 1}^6$ at $t_0 = 0$.
As shown in Fig.~\ref{fig:test-deformation} (a), 
 the five colored Yin sets constitute a circular disk
 with its radius as $r = 0.15$ and its center at $[0.5, 0.5]^T$ 
 while the last Yin set 
 is the unbounded complement of the circle.

The above IT problem of six phases
 is numerically solved by the cubic MARS method
 with the curvature-based ARMS strategy
 (\ref{eq:curvature-based-arms}) specified by
 \begin{equation}
   \label{eq:curvatureStrategyConstants}
   \begin{array}{l}
   r_{\text{tiny}} = 0.05; \
   (\rho_{\min}^c, \rho_{\max}^c ) = (10^{-5}, 1); \
   r_{\min}^c = 0.1, 0.05;\
   \sigma^c(x) = x.
   \end{array}
 \end{equation}

Different from those in (\ref{eq:vortex-test-parametersCB}), 
 parameter values in (\ref{eq:curvatureStrategyConstants})
 are more suitable for the deformation tests.
For example,
 the higher value of $\rho^c_{\max}$
 accounts for the much larger percentage of low-curvature arcs
 in Fig.~\ref{fig:test-deformation}
 than that in Fig.~\ref{fig:test-quartered-circular-disk-T4}.
As $T$ increases from 2 to 4,
 we decrease $r_{\min}^c$ from $0.1$ to $0.05$
 to accommodate the larger deformations.
The maximum density-increase ratio $R_{\max}$,
 which is also the maximum ratio of the longest chordal length
 over the shortest one,
 equals 200 and 400 for $T=2$ and $T=4$, respectively. 
 
\begin{table}
  \caption{Errors and convergence rates
    of the multiphase cubic MARS method
    with $k = \frac{1}{8}h$ and curvature-based ARMS
    in solving the deformation test of $T = 2, 4$.
  }
  \centering
  \renewcommand{\arraystretch}{1.2} 
  \input{{\TABDIR}table6.3}
  \label{tab:deformation-curvature-arms}
\end{table}

\begin{table}
  \caption{Errors of individual phases and convergence rates
    of the multiphase cubic MARS method
    in solving the deformation test of $T = 2$.
    Errors in the second tabular are at $t = \frac{1}{2}T$
    and are calculated by Richardson extrapolation.
  }
  \centering
  \renewcommand{\arraystretch}{1.2} 
  \input{{\TABDIR}table6.4}
  \label{tab:compare-result}
\end{table}

As shown in Fig.~\ref{fig:test-deformation}, 
 the interface graph $G_{\Gamma} = (V_{\Gamma}, E_{\Gamma},\psi_{\Gamma})$
 that represents the interface topology
 remains the same during the entire simulation: 
 $E_{\Gamma}$ always consists of the ten edges that connect
 the five T junctions located on the circle 
 and the junction of degree 5 at the disk center. 
By Algorithm \ref{alg:partitionGraph},
 $C_S$ contains a single circuit formed by the five T junctions 
 while $T_S$ has five trails
 that correspond to the five disk radii.
Despite the tremendous deformations,
 each phase remains connected, 
 demonstrating the capability of our method
 in preserving topological structures.

As shown in Fig.~\ref{fig:test-deformation}(g,i),
 the boundary lengths increase and then decrease.
Accordingly, in Fig.~\ref{fig:test-deformation}(h,j),
 the number of interface markers for each phase 
 first increases, then stagnates roughly as a constant,
 and finally decreases.
For the number of markers, the increase and decrease 
 are clearly driven by those of the boundary length
 while the stagnation 
 follows from the finite width of the interval $[r_{\Tiny}h_L, h_L]$
 and the fact that it takes time 
 for the distances of adjacent markers
 to decrease from $h_L$ to $r_{\Tiny}h_L$.
At the end of the simulation,
 the number of markers for each phase 
 is roughly twice as much as that at the initial time.
This ratio being around 2,
 together with the similarity between
 subplots (g,i) and (h,j) in Fig.~\ref{fig:test-deformation},
 illustrates the versatility and effectiveness of ARMS
 in managing the regularity of interface markers.

\begin{table}
  \caption{Accuracy comparison of the proposed method
    ($h_L^c = 0.2 h$; $k=\frac{1}{8}h$) with some VOF methods
    based on the total IT error $\sum_{i = 1}^3 E_i$
    for the three-phase deformation test
    in \cite[Sec. 3.6]{sijoy2010volume}.
    The curvature-based ARMS is specified by
    (\ref{eq:curvatureStrategyConstants}) and $r_{\min}^c = 0.1$.
    The errors in the last three lines
    are taken from \cite[Tab. 8]{sijoy2010volume},
    with the two-letter acronyms defined in the caption
    of Table \ref{tab:comparisonVortex}.
  }
  \centering
  \renewcommand{\arraystretch}{1.2} 
  \input{{\TABDIR}table6.5}
  \label{tab:comparisonDeformation}
\end{table}

The IT errors and convergence rates of the proposed MARS method
 are listed in Table \ref{tab:deformation-curvature-arms}, 
 where convergence rates
 are close to 4, 6, and 8
 for the choices of $h_L^c$ being $O(h)$, $O(h^{\frac{3}{2}})$, 
 and $O(h^2)$, respectively.
In Table \ref{tab:compare-result},
 we compare results of different ARMS strategies.
Convergence rates of constant ARMS
 vary from one phase to another, 
 indicating that,
 even with $r_{\Tiny}=0.005$,
 the computation has not yet reached the asymptotic range.
In contrast,
 convergence rates of curvature-based ARMS 
 are more phase-independent,
 indicating that high-curvature arcs
 have been well resolved by the more efficient marker distributions.
Finally, an increase of the time step size 
 from $k=\frac{1}{8}h$ to $k=h$
 yields convergence rates very close to five, 
 implying the dominance of temporal discretization errors
 over the spatial approximation errors. The last tabular in 
 Table \ref{tab:compare-result} contains errors of phases
 at $t = \frac{1}{2}T$ by Richardson extrapolation, 
 indicating that the converge rates of our method
 are independent of periodicity of the velocity field.

Results of MARS and VOF methods
 in solving the deformation test in \cite[Sec. 3.6]{sijoy2010volume}
 are compared in Table \ref{tab:comparisonDeformation},
 where the proposed method 
 is shown to be more accurate
 than these VOF methods by many orders of magnitude.

\begin{table}
  \caption{\revise{CPU time in seconds
    of the multiphase cubic MARS method
    in solving the deformation test 
    in Sec.~\ref{sec:testDeformation} on 
    an Intel Core i7-12700 at 4.2GHz.
    The parameters are identical 
    to those in Table~\ref{tab:deformation-curvature-arms}}.
  }
  \centering
  \renewcommand{\arraystretch}{1.2} 
  \input{{\TABDIR}tableCpu1}
  \label{tab:deformation-curvature-arms-Cputime}
\end{table}

\begin{figure}[t]
  \centering
  \subfloat[$t = 0$]{
    \includegraphics[width=0.3\linewidth]{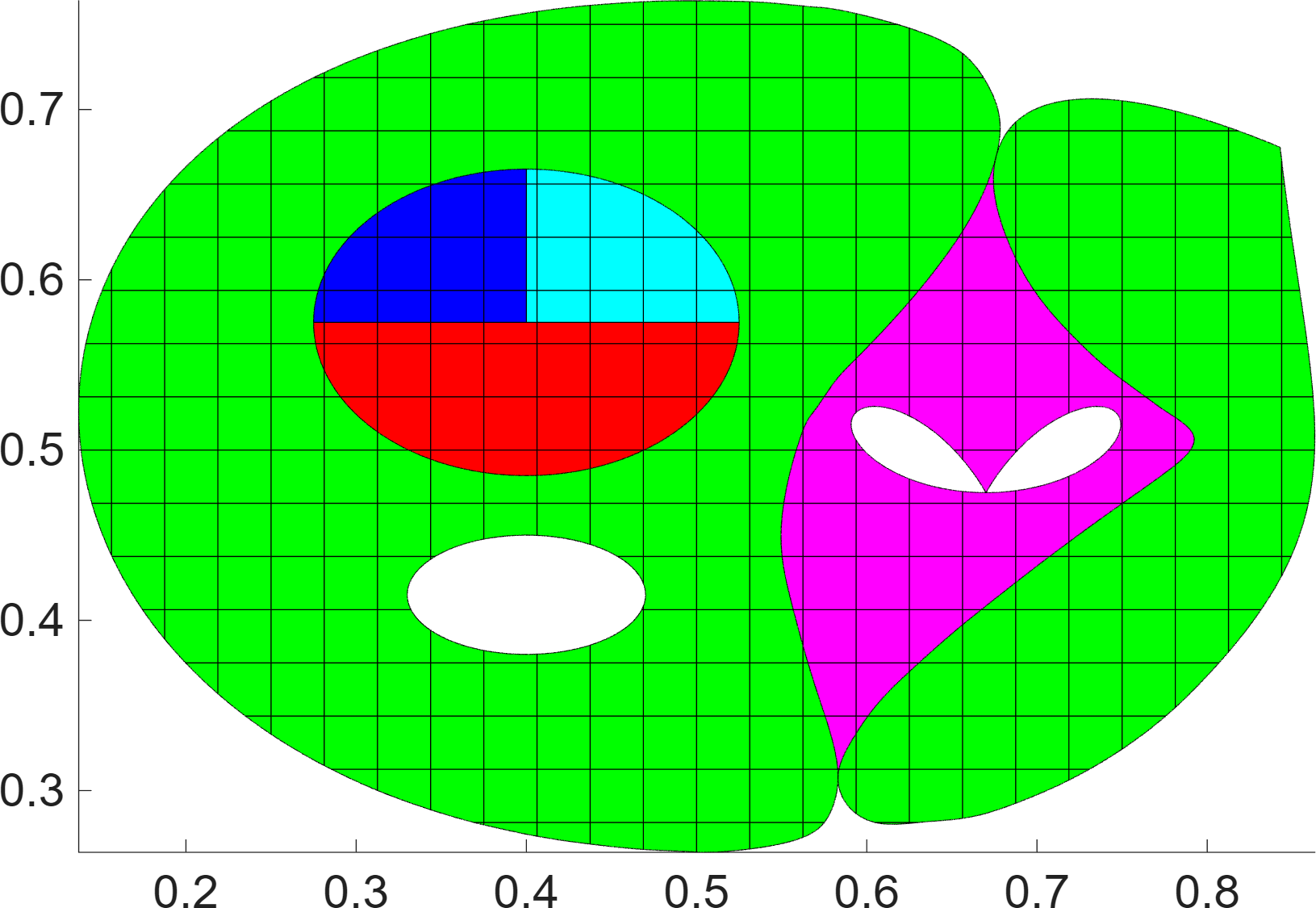}
  }
  \hfill
  \subfloat[$t = \frac{1}{4}T_v$; vortex]{
    \includegraphics[width=0.285\linewidth]{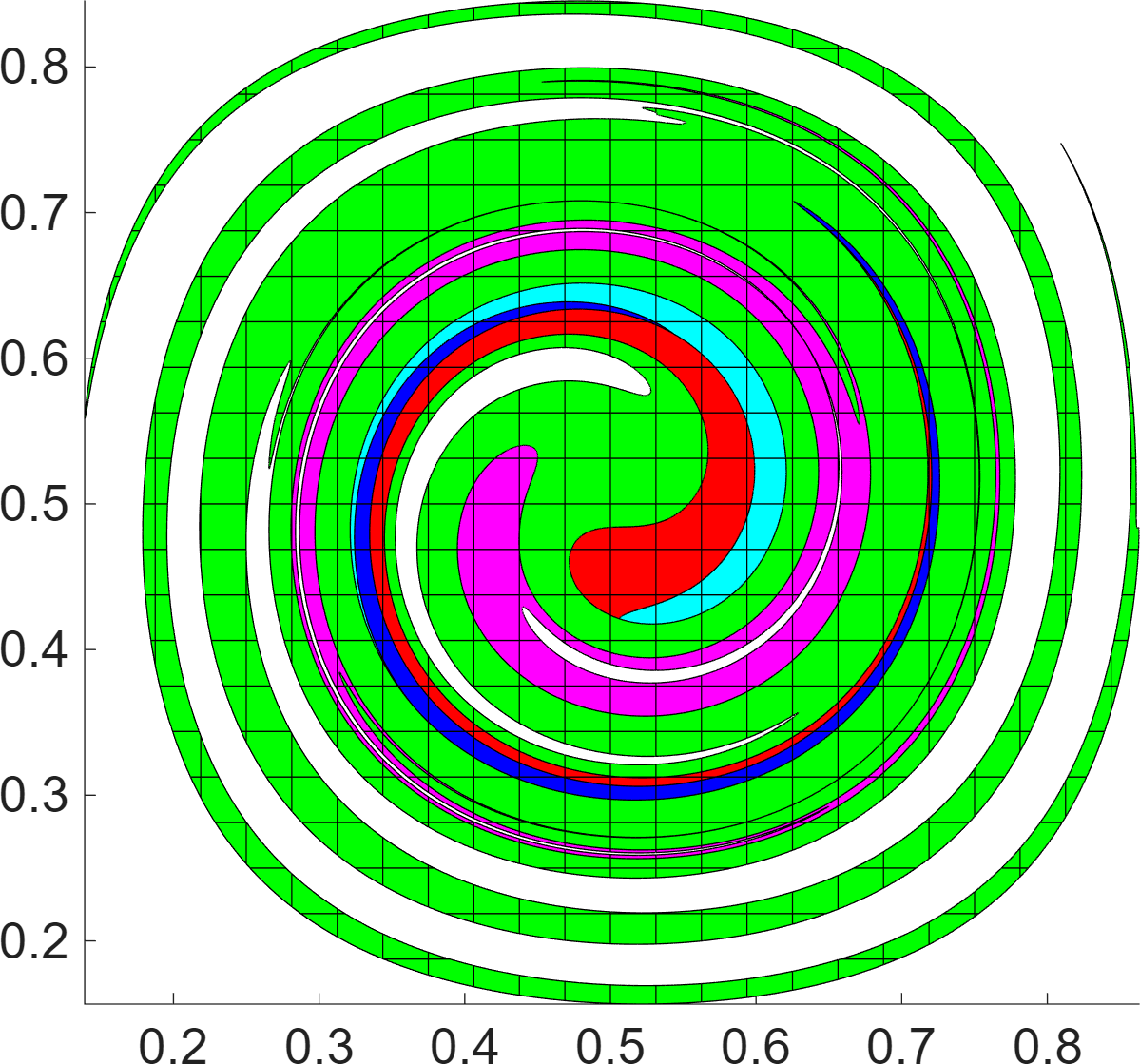}
  }
  \hfill
  \subfloat[$t = \frac{1}{4}T_d$; deformation]{
    \includegraphics[width=0.34\linewidth]{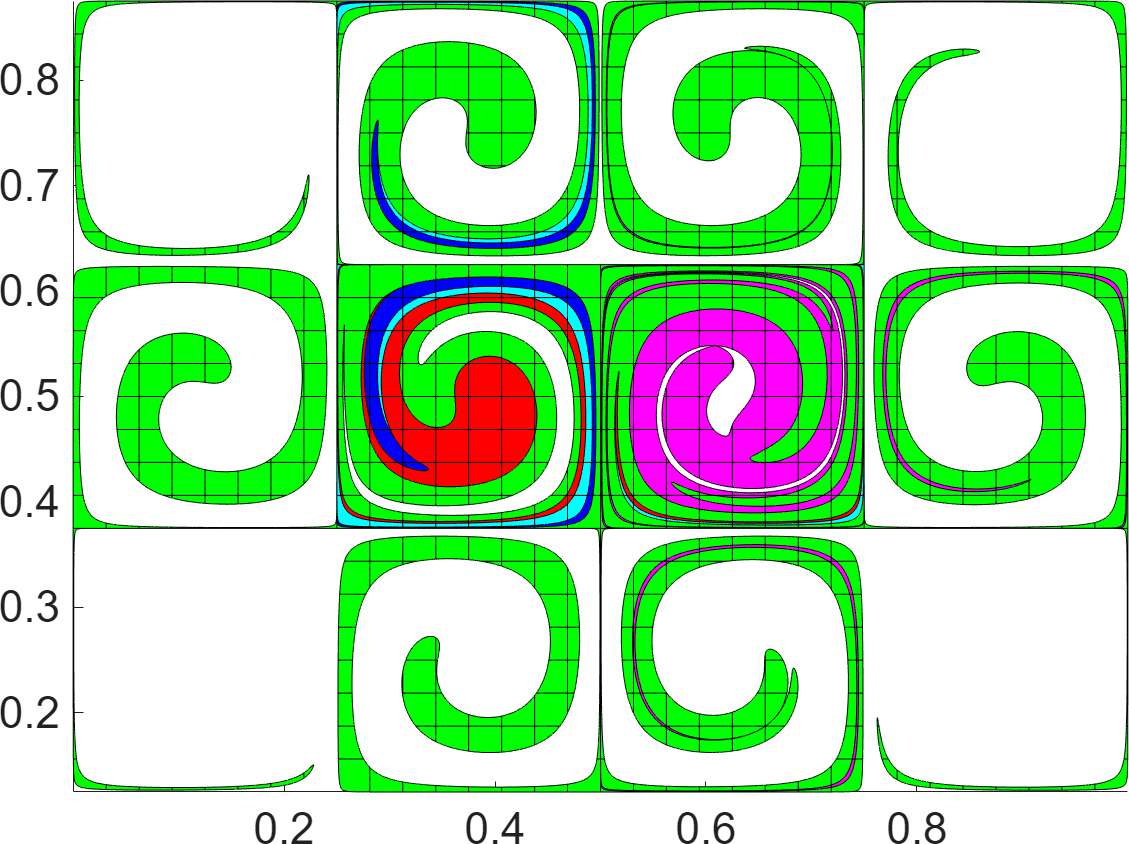}
  }

  \subfloat[$t = \frac{1}{2}T_v$; vortex]{
    \includegraphics[width=0.285\linewidth]{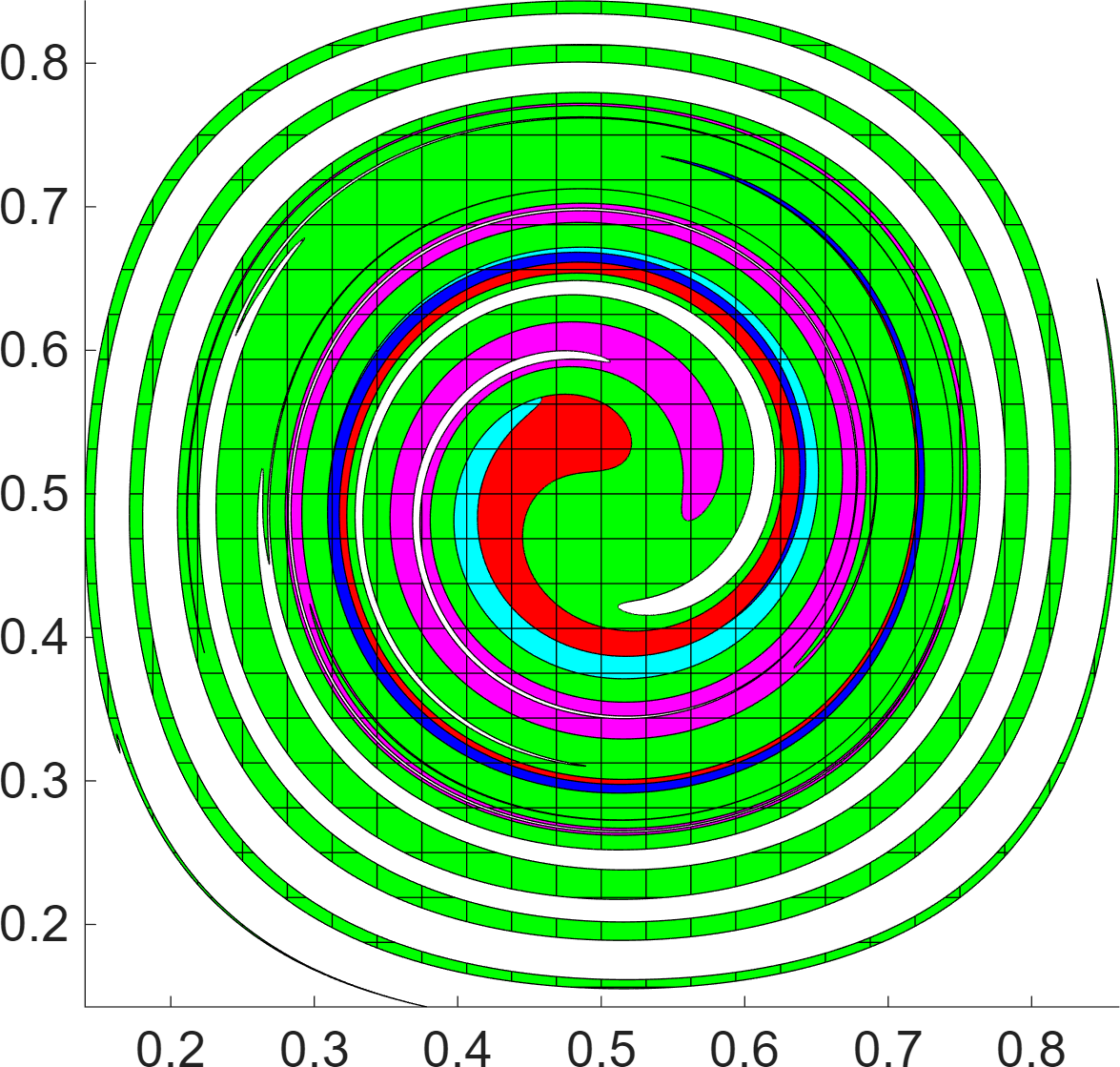}
  }
  \hfill
  \subfloat[$t = \frac{1}{2}T_d$; deformation]{
    \includegraphics[width=0.34\linewidth]{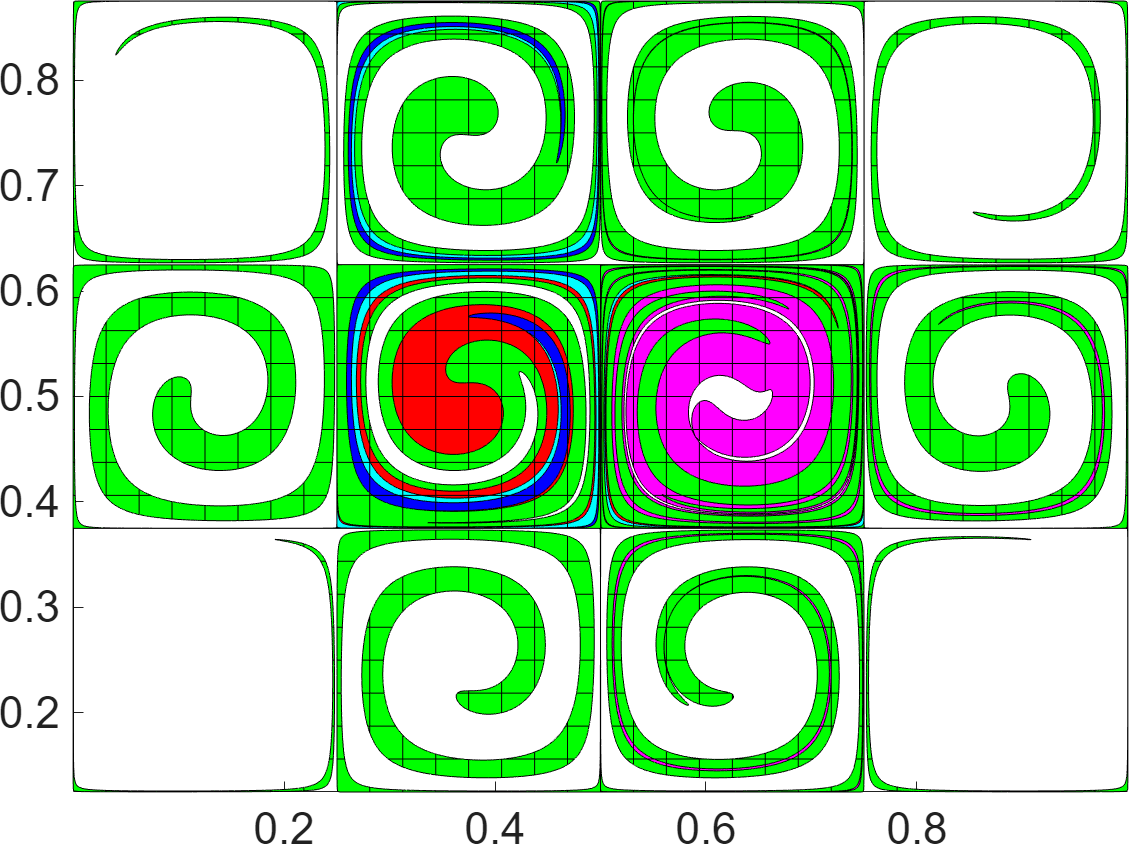}
  }
  \hfill
  \subfloat[the end of simulation]{
    \includegraphics[width=0.3\linewidth]{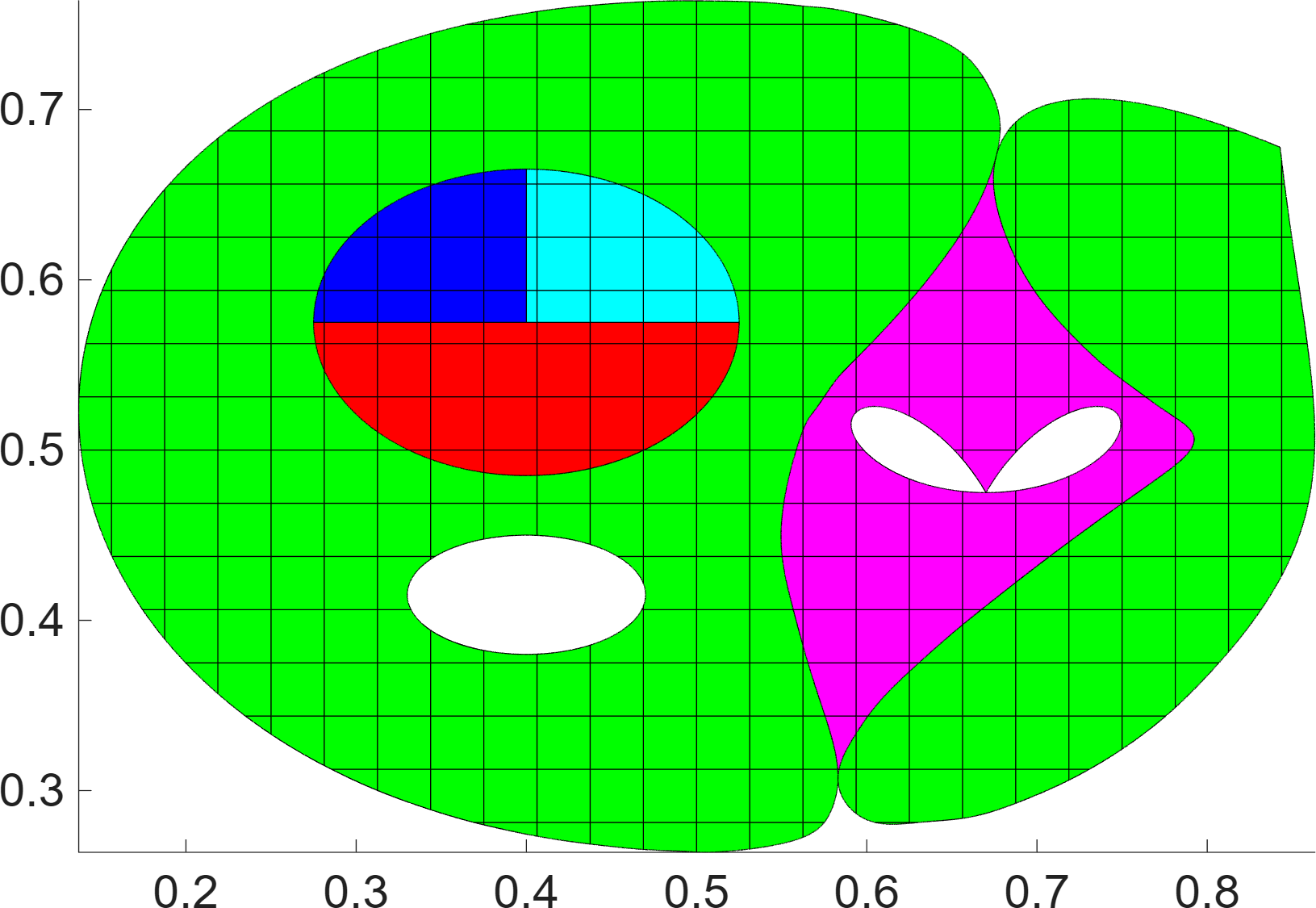}
  }

  \subfloat[vortex: normalized length
  of $\partial {\cal M}_i(\frac{t}{T_v})$]{
    \includegraphics[width=0.22\linewidth]{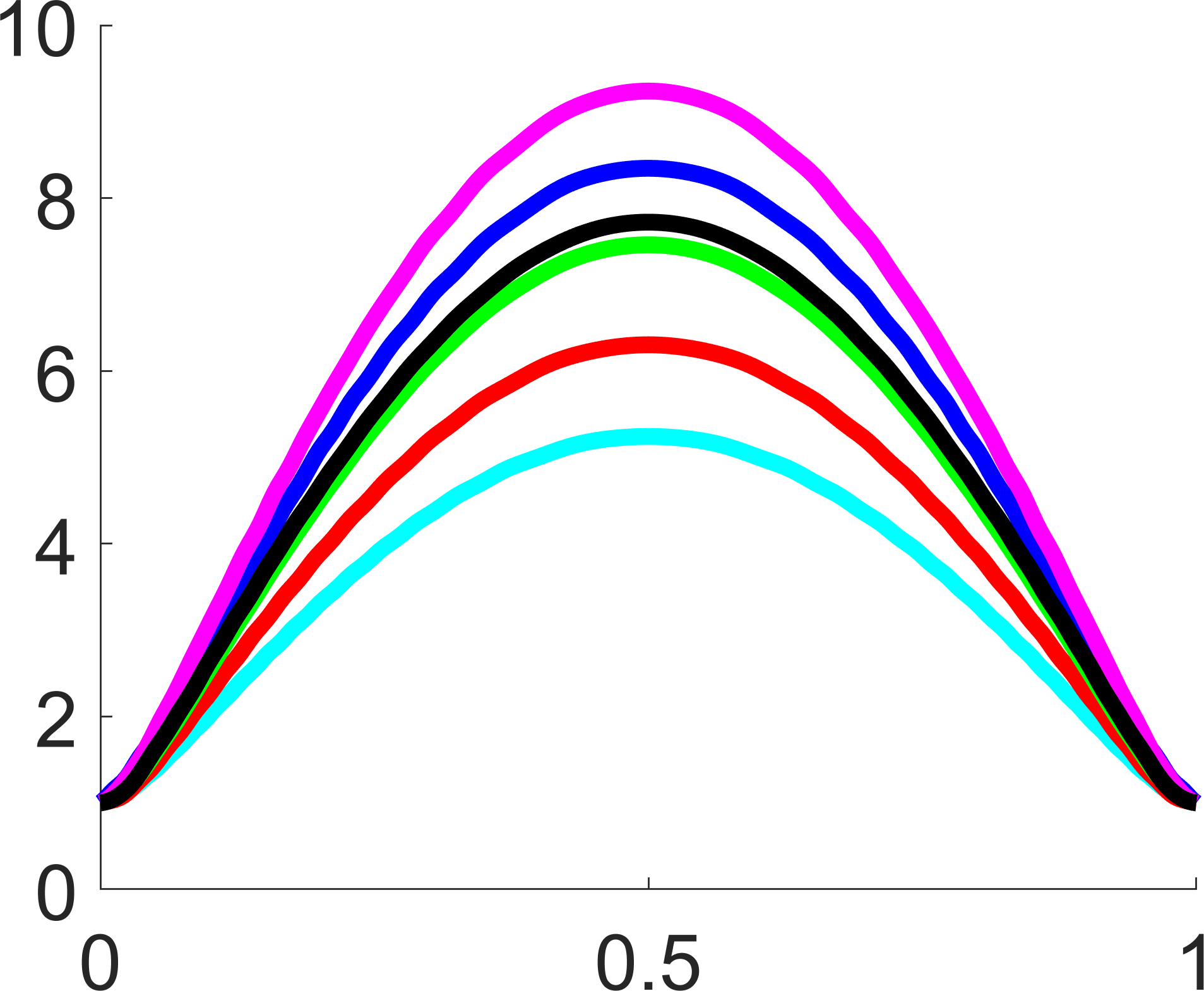}
  }
  \hfill
  \subfloat[vortex: normalized marker number in $\partial {\cal M}_i^n$]{
    \includegraphics[width=0.22\linewidth]{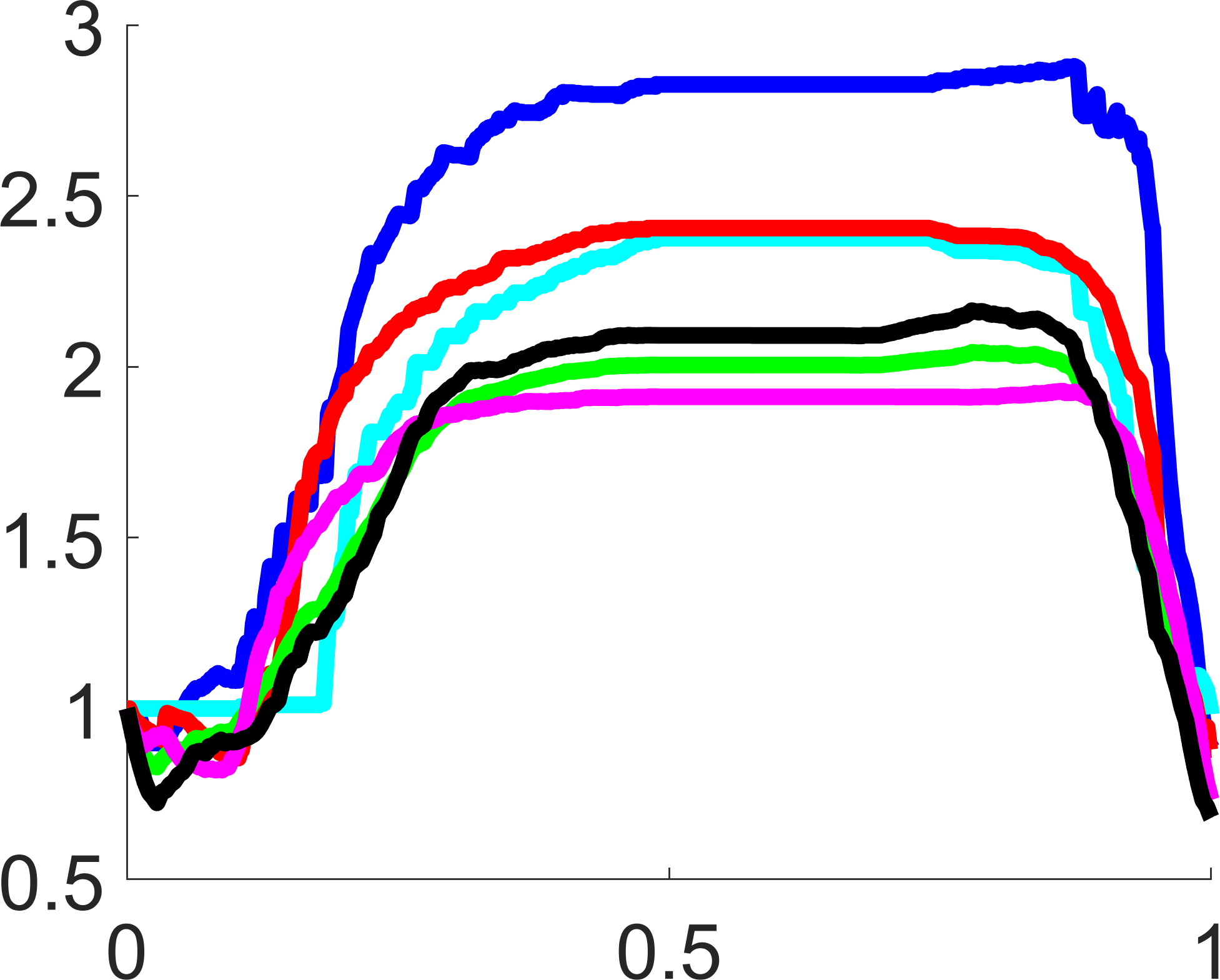}
  }
  \hfill
  \subfloat[deformation: normalized length
  of $\partial {\cal M}_i(\frac{t}{T_d})$]{
    \includegraphics[width=0.22\linewidth]{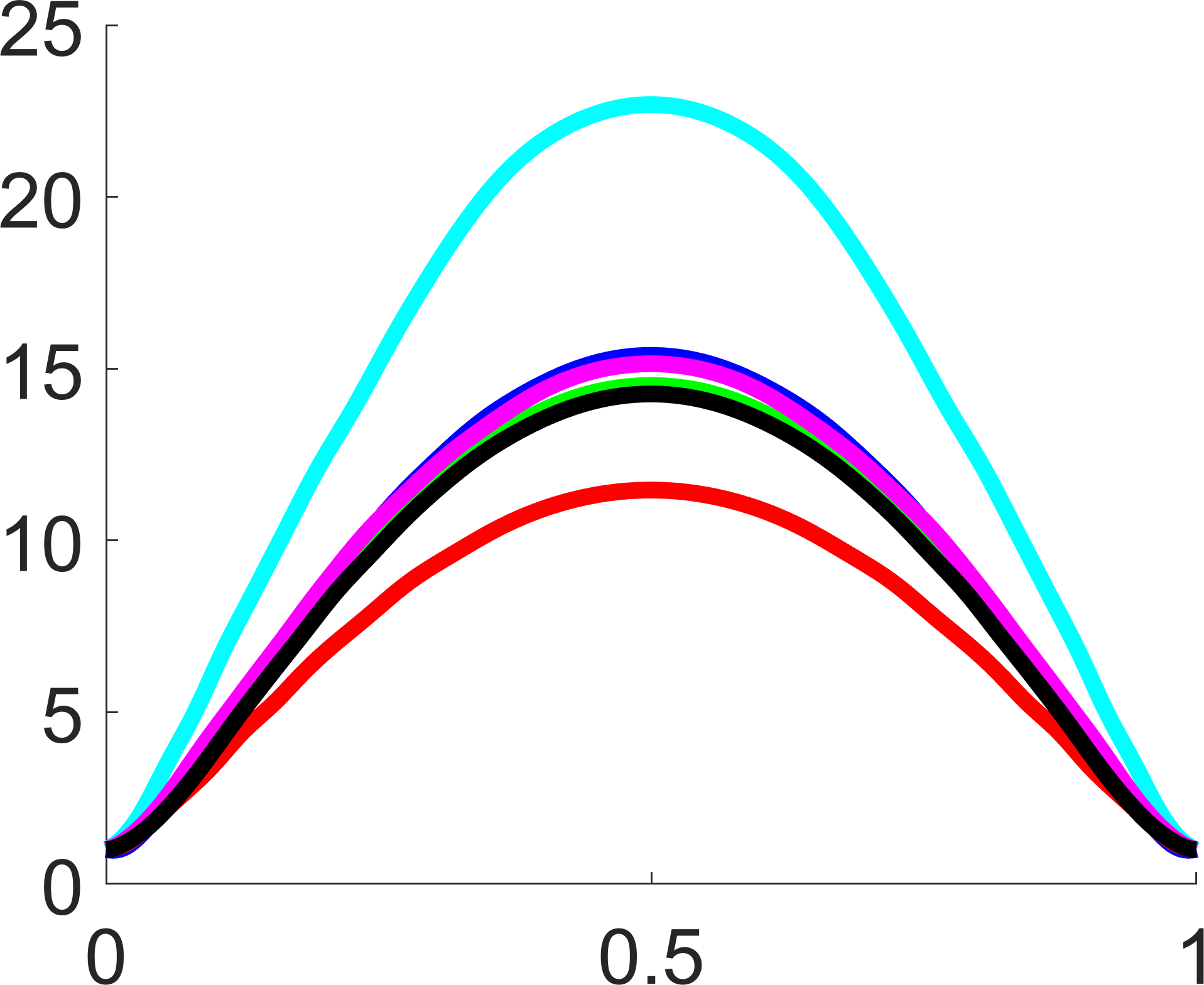}
  }
  \hfill
  \subfloat[deformation: normalized marker number
  in $\partial {\cal M}_i^n$]{
    \includegraphics[width=0.22\linewidth]{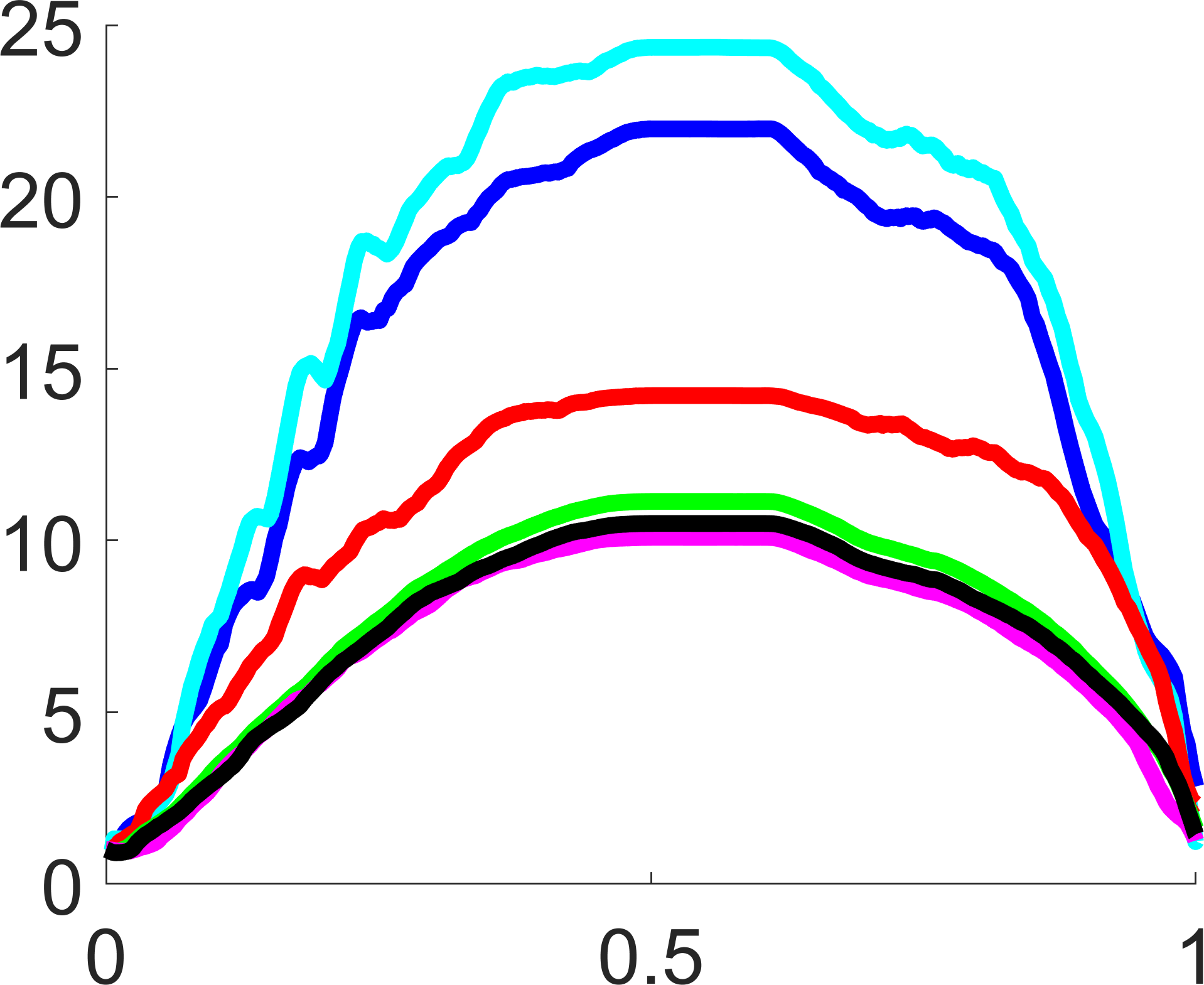}
  }
  \caption{Solutions of the cubic MARS method with $h_L^c=0.2h$
    for the vortex shear test ($T_v = 16$)
    and the deformation test ($T_d = 4$) 
    on the Eulerian grid of $h = \frac{1}{32}$.
    Subplot (a) shows the initial Yin sets, 
    which are the same 
    as those in Fig.~\ref{fig:multiphase-representation}(a),
    and subplot (b) presents the final solutions of the two tests,
    which are visually indistinguishable.
    Subplots (c)--(f) are solution snapshots at key time instances.
    In subplots (g)--(j),
    each phase is represented by a curve of the same color
    except that the unbounded white phase
    is represented by the black curve.
    The initial distances between markers are $\frac{1}{2}h_L(\rho)$,
    where $h_L$ is defined in (\ref{eq:curvature-based-arms}).
}
  \label{fig:test-multiphase-representation}
\end{figure}

\revise{
Finally, 
 the computational cost of the proposed method
 for solving the deformation test
 is presented
 in Table \ref{tab:deformation-curvature-arms-Cputime}. 
For the choices of $h_L^c$
 being $O(h)$, $O(h^{\frac{3}{2}})$, and $O(h^2)$,
 the rate of growth of the CPU time
 is slightly over 2, 2.5, and 3;
 these results confirm Lemma \ref{lem:complexityOfARMS}
 that the complexity of the multiphase MARS method per time step 
 is linear with respect to the number of markers.
} 

\subsection{Vortex shear and deformation
  of the six phases
  in Fig.~\ref{fig:multiphase-representation}(a)}
\label{sec:testMultiphase-representation}

Flow maps of this test are 
 (\ref{eq:vortex-velocity}) and (\ref{eq:deformation-velocity}).
The initial conditions are the five Yin sets
 $(\mathcal{M}_i)_{i=1}^5$
 shown in Fig.~\ref{fig:multiphase-representation}(a), 
 whose boundaries are approximated to sufficient accuracy
 by $\mathcal{C}^4$ quintic splines, 
 elliptical arcs, linear segments, and
 \revise{a partial rose curve $r = a \sin(3 \theta)$
   where $r$ denotes the radial distance from the origin, 
   $\theta\in \left[\frac{2 \pi}{3}, \frac{4 \pi}{3}\right]$
   the polar angle, 
   and $a$ a scaling parameter for the length of the two petals.
   As indicated in Fig.~\ref{fig:multiphase-representation}(a),
   although $e_{10}$ and $e_{11}$ are smoothly connected, 
   $e_{15}$ and $e_{16}$ are not. 
   Consequently, the maximally connected smooth curve
   $(e_{15},e_{10},e_{11},e_{16})$
   is identified as a trail instead of a cycle
   and is approximated by a not-a-knot spline
   instead of a periodic spline.
 }

\begin{table}
  \caption{Errors and convergence rates of 
    the multiphase cubic MARS method
    with $k = \frac{1}{8}h$ and curvature-based ARMS 
    in solving the vortex shear test at $T_v = 4, 8, 12, 16$ 
    and the deformation test at $T_d = 2, 4$.
    The initial Yin sets are 
    shown in Fig.~\ref{fig:test-multiphase-representation}(a).}
  \centering
  \small
  \renewcommand{\arraystretch}{1.2} 
  \input{{\TABDIR}table6.6}
  \label{tab:test-graph41}
\end{table}

Parameters of the curvature-based ARMS strategy 
 for the vortex shear are 
 \begin{equation}
  \label{eq:curvatureStrategyGraphVortex}
  \begin{array}{l}
  r_{\text{tiny}} = 0.1; \quad
  (\rho_{\min}^c, \rho_{\max}^c ) = (10^{-5}, 0.2); \quad
  \sigma^c(x) = x
  \end{array}
\end{equation}
while those for the deformation are
\begin{equation}
  \label{eq:curvatureStrategyGraphDeformation}
  \begin{array}{l}
  r_{\text{tiny}} = 0.05; \quad
  (\rho_{\min}^c, \rho_{\max}^c ) = (10^{-5}, 1); \quad
  \sigma^c(x) = x,
  \end{array}
\end{equation}
where the parameter $\rho_{\max}^c$
 is selected based on the characteristics of flow fields.
As for $r_{\min}^c$,
 we choose it as a monotonically decreasing function of the period;
 see Table~\ref{tab:test-graph41}
 for its value of each test case.

As shown in Fig.~\ref{fig:test-multiphase-representation}, 
 the interface graph 
 that represents the interface topology
 remains the same during the entire simulation
 for both the vortex shear test and the deformation test:
 $C_S$ always contains three circuits 
 and $T_S$ always has four trails.
Despite the extremely large deformations,
 each phase remains connected. 
These invariants demonstrate the capability of
 the multiphase cubic MARS method
 in preserving topological structures.

The evolution of the boundary length and number of markers
 for each test and each phase
 is shown in Fig.~\ref{fig:test-multiphase-representation}
 (g)--(j), 
 demonstrating the effectiveness and versatility of ARMS
 in maintaining the $(r,h)$-regularity 
 even for geometrically and topologically complex interfaces with 
 high curvature ($r_{\min}^c \leq 0.01$).

Finally,
 we list, for this test,
 the IT errors and convergence rates in Table \ref{tab:test-graph41}, 
 which clearly demonstrate
 the fourth-, sixth-, and eighth-order convergence rates
 for $h_L^c=O(h)$, $O(h^{\frac{3}{2}})$, and $O(h^2)$, respectively.

\revise{
\subsection{Vortex shear of a piggy}
\label{sec:testPig}
}

\begin{figure}
  \centering
  \subfloat[$t = 0$]{
    \includegraphics[width=0.23\linewidth]{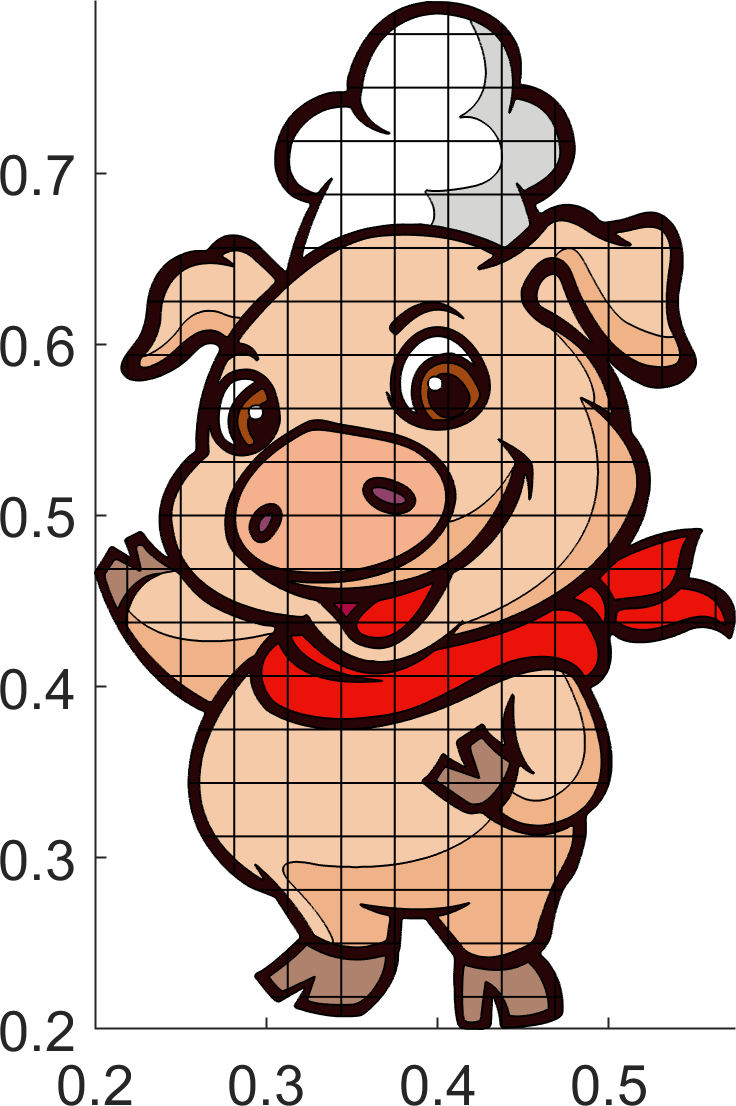}
  }
  \hfill
  \subfloat[$T = 4; t = \frac{1}{4}T$]{
    \includegraphics[width=0.29\linewidth]{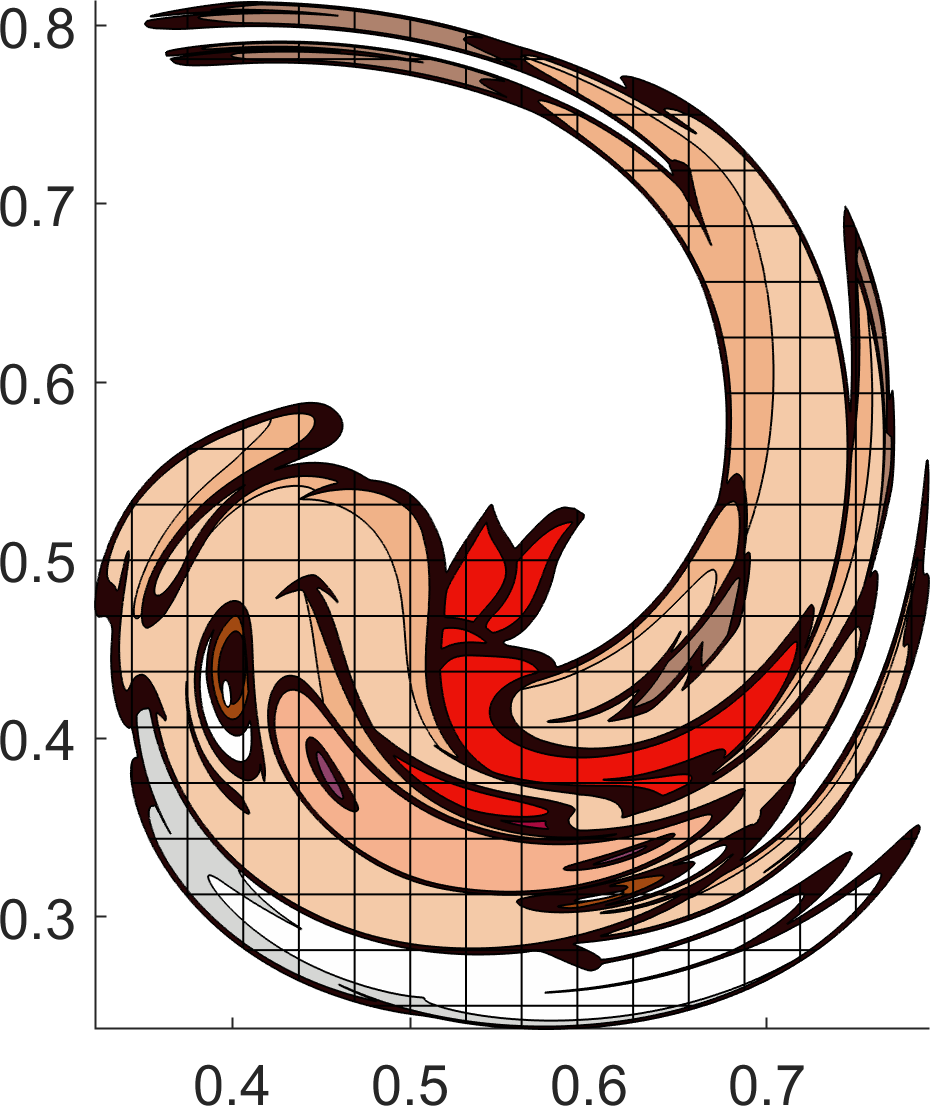}
  }
  \hfill
  \subfloat[$T = 8; t = \frac{1}{4}T$]{
    \includegraphics[width=0.37\linewidth]{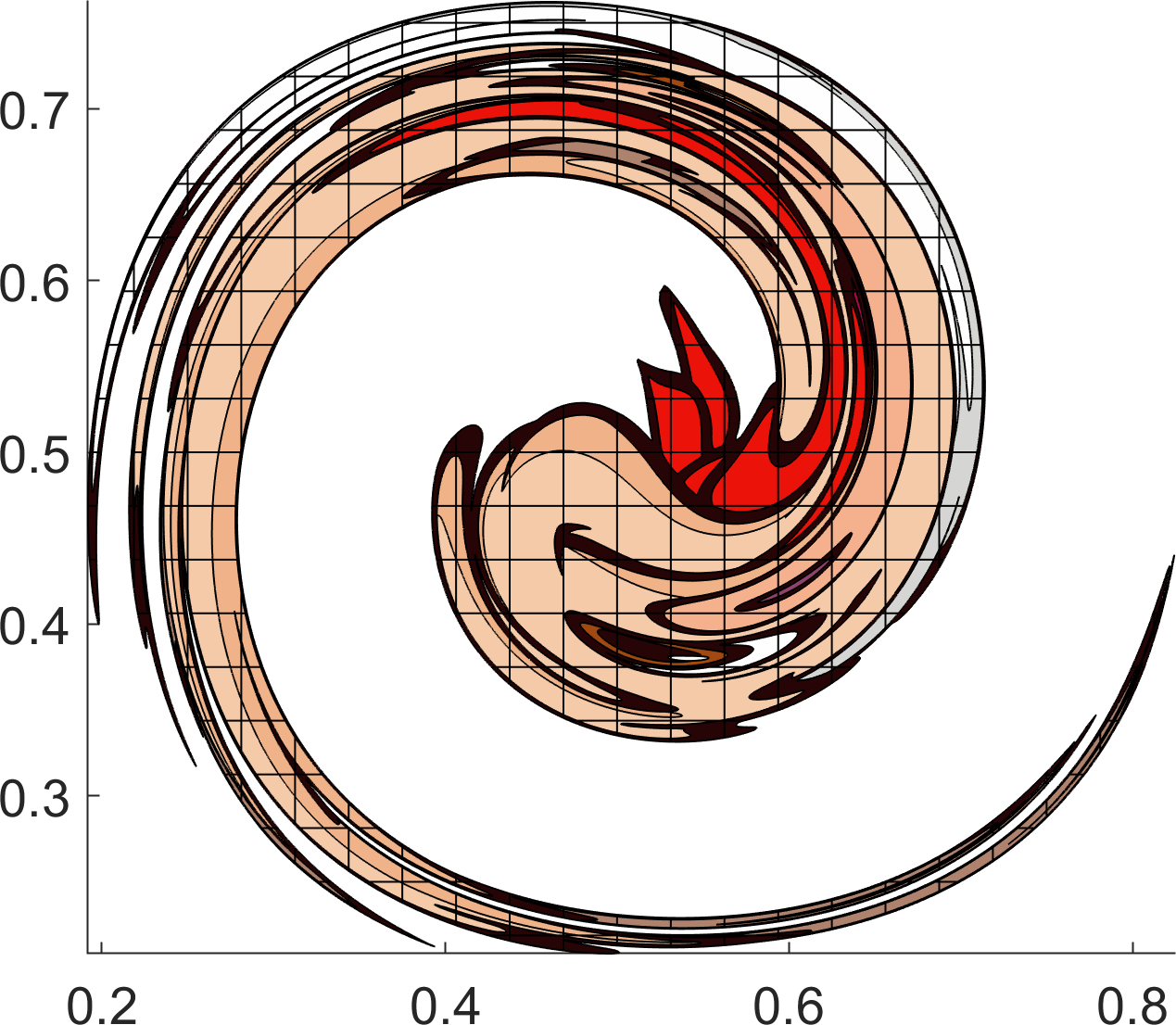}
  }

  \subfloat[$T = 4; t = \frac{1}{2}T$]{
    \includegraphics[width=0.335\linewidth]{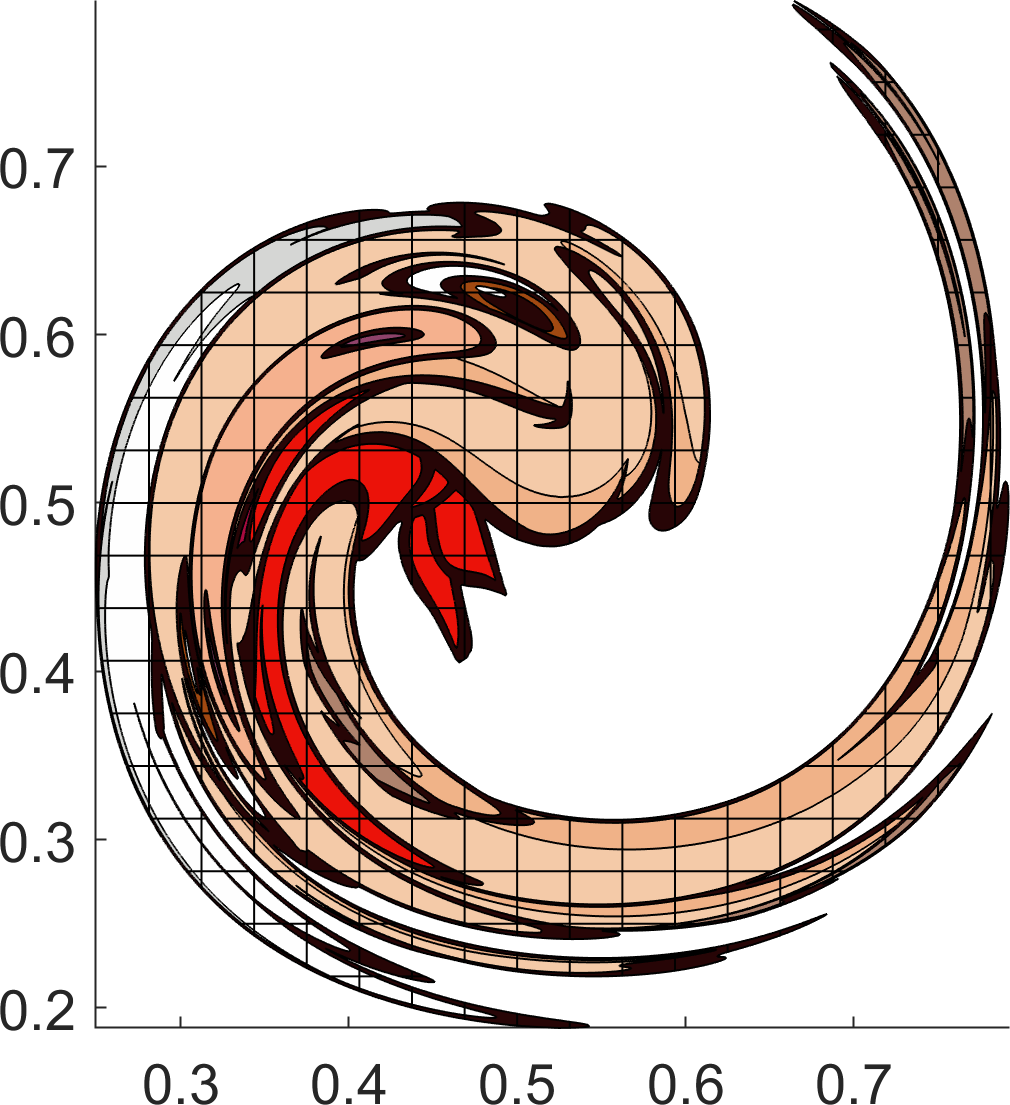}
  }
  \hfill
  \subfloat[$T = 8; t = \frac{1}{2}T$]{
    \includegraphics[width=0.33\linewidth]{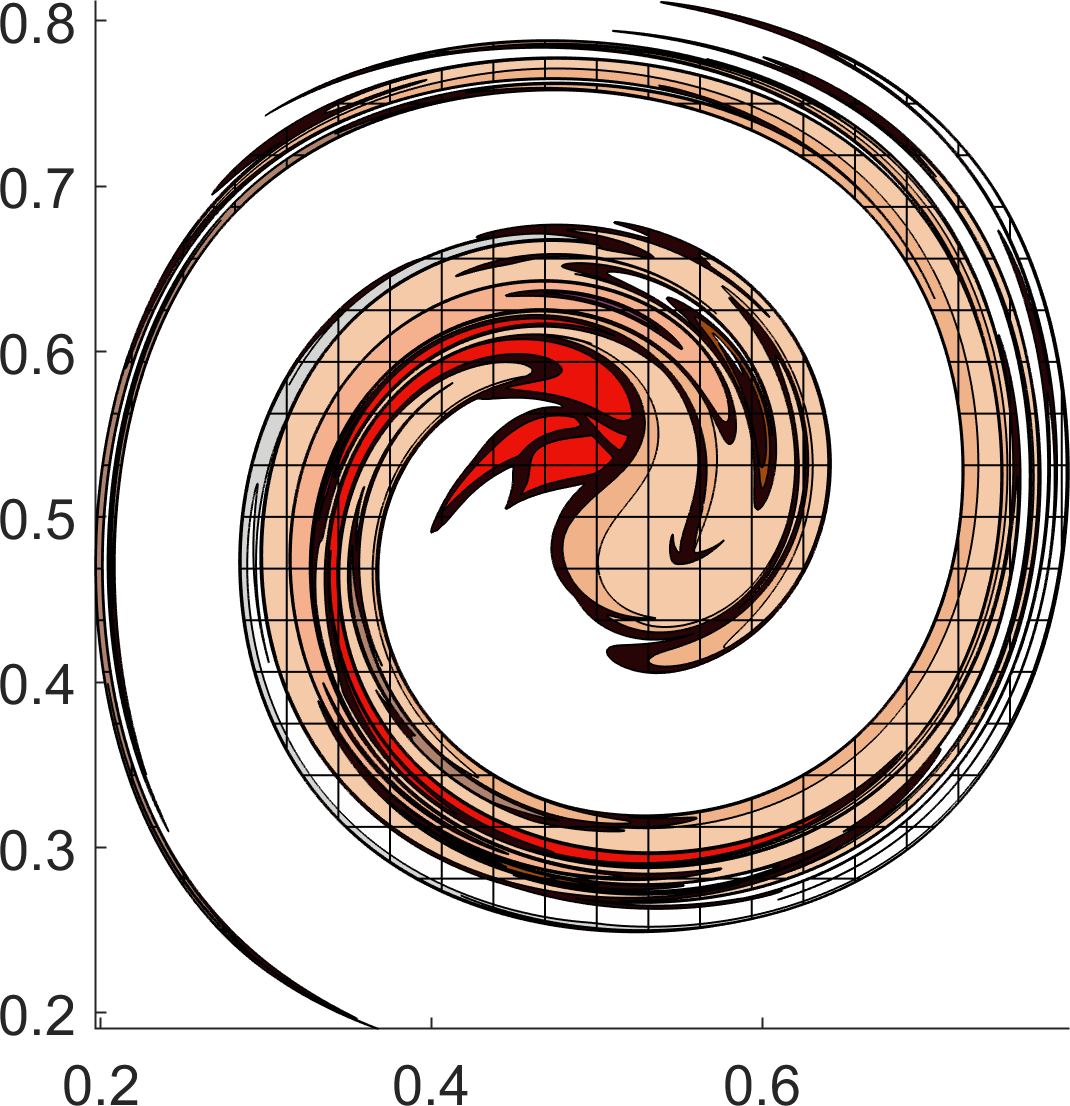}
  }
  \hfill
  \subfloat[the end of simulation]{
    \includegraphics[width=0.25\linewidth]{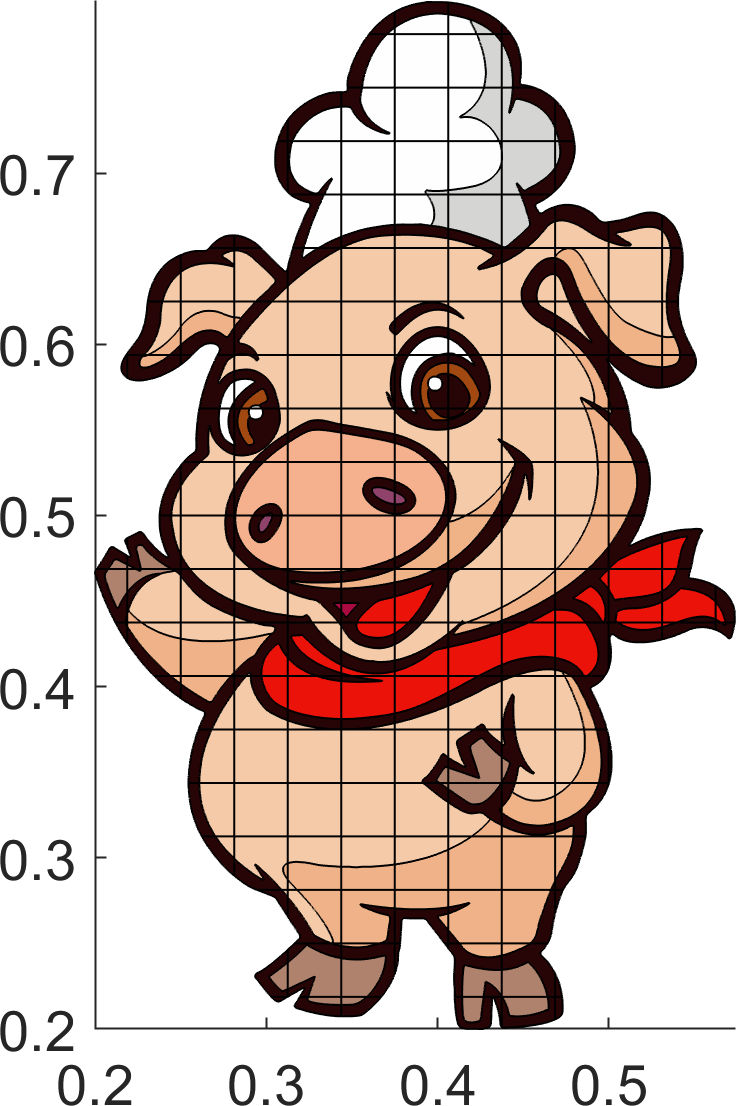}
  }
  
  \caption{Solutions of the multiphase cubic MARS method
    with $h_L^c = 0.1 h$
    for the vortex test of a piggy with $T = 4, 8$ 
    on the Eulerian grid of $h = \frac{1}{32}$.
    Subplots (a)--(f) are snapshots of the solution
    at key time instants. 
    The initial distances between markers are $\frac{1}{2}h_L(\rho)$,
    where $h_L$ is defined in (\ref{eq:curvature-based-arms}).
    The final solutions for $T=4,8$ overlap and are visually indistinguishable.
  }
  \label{fig:test-pig}
\end{figure}

The flow map of this test is the same as that 
 in (\ref{eq:vortex-velocity}).
As shown in Fig.~\ref{fig:test-pig}(a),
 the piggy at the initial time
 is constituted by 15 phases, 
 whose boundaries are approximated to very high accuracy
 by ${\cal C}^4$ quintic splines
 fitted through a sufficiently large number of breakpoints.

In setting up the initial topology for each run, 
 we include in the vertex set of the interface graph
 all junctions of the interface.  
On the other hand, we set $K_{\Gamma}=\emptyset$, 
 i.e., any non-smooth point in Definition \ref{def:kink}
 is treated as a smooth point with high curvature.
The purposes of the different treatments of junctions and non-smooth points are
 to be worry free about the special handling of kinks
 and to test the ability of curvature-based ARMS
 in resolving high-curvature arcs.
For each different choice of $h_L^c$ in Table \ref{tab:pig-arms}, 
 the initial interface is re-sampled from the finest representation
 in the previous paragraph
 so that all breakpoint sequences at the initial time
 satisfy the corresponding $(r_{\Tiny},h_L)$-regularity.
 
\begin{table}
  \caption{Errors and convergence rates
    of the multiphase cubic MARS method
    with $k = \frac{1}{8}h$ and curvature-based ARMS
    in solving the vortex shear of a piggy for $T = 4, 8$.
    The initial condition is shown in Fig.~\ref{fig:test-pig}(a)
    and the parameter values in
    \eqref{eq:curvatureStrategyGraphVortex} are adopted.
  }
  \centering
  \renewcommand{\arraystretch}{1.2} 
  \input{{\TABDIR}tablePig}
  \label{tab:pig-arms}
\end{table}
As shown in Fig.~\ref{fig:test-pig}, 
 the proposed MARS method maintains the constant interface topology 
 throughout the entire simulation. 
In Table~\ref{tab:pig-arms}, 
 we list convergence rates of the proposed MARS method,
 which are close to 4, 6, and 8 for 
 $h_L^c=O(h)$, $O(h^{\frac{3}{2}})$, and $O(h^2)$, respectively. 

The treatment of non-smooth points as smooth breakpoints via curvature-based ARMS
 is a solution to both static and dynamic kinks, 
 providing a flexible mechanism to balance accuracy and convenience. 
First, although kinks cause the slight deficit of convergence rates,
 the IT errors are already small enough in Table \ref{tab:pig-arms}.
Second, it is well-known that an initially smooth arc
 might develop new kinks,
 even under the flow map of a ${\cal C}^{\infty}$ velocity.
It is highly difficult 
 to detect the exact location of these dynamic kinks, 
 but their neighborhoods can be well resolved with curvature-based ARMS.
Lastly, this treatment still preserves of all phase adjacency,
 i.e., it generates neither overlaps nor vacuums.
However, all junctions must be included
 as vertices of the interface graph; 
 otherwise it would create overlaps and vacuums
 of adjacent phases. 
\newline

\revise{
\subsection{Vortex shear and deformation of a raccoon}
\label{sec:testRaccoon}
}

\begin{figure}
  \centering
  \subfloat[$t = 0$]{
    \includegraphics[width=0.26\linewidth]{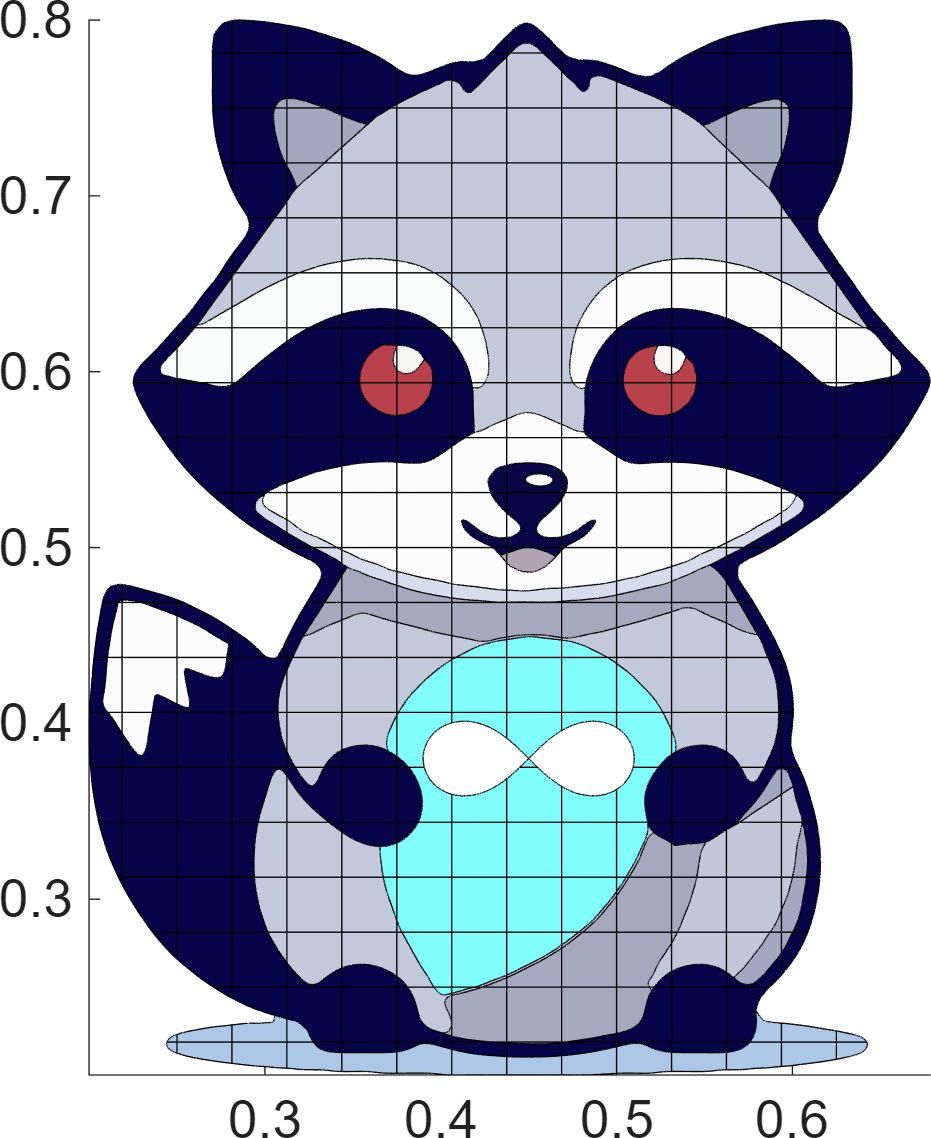}
  }
  \hfill
  \subfloat[$t = \frac{1}{4}T_v$; vortex]{
    \includegraphics[width=0.325\linewidth]{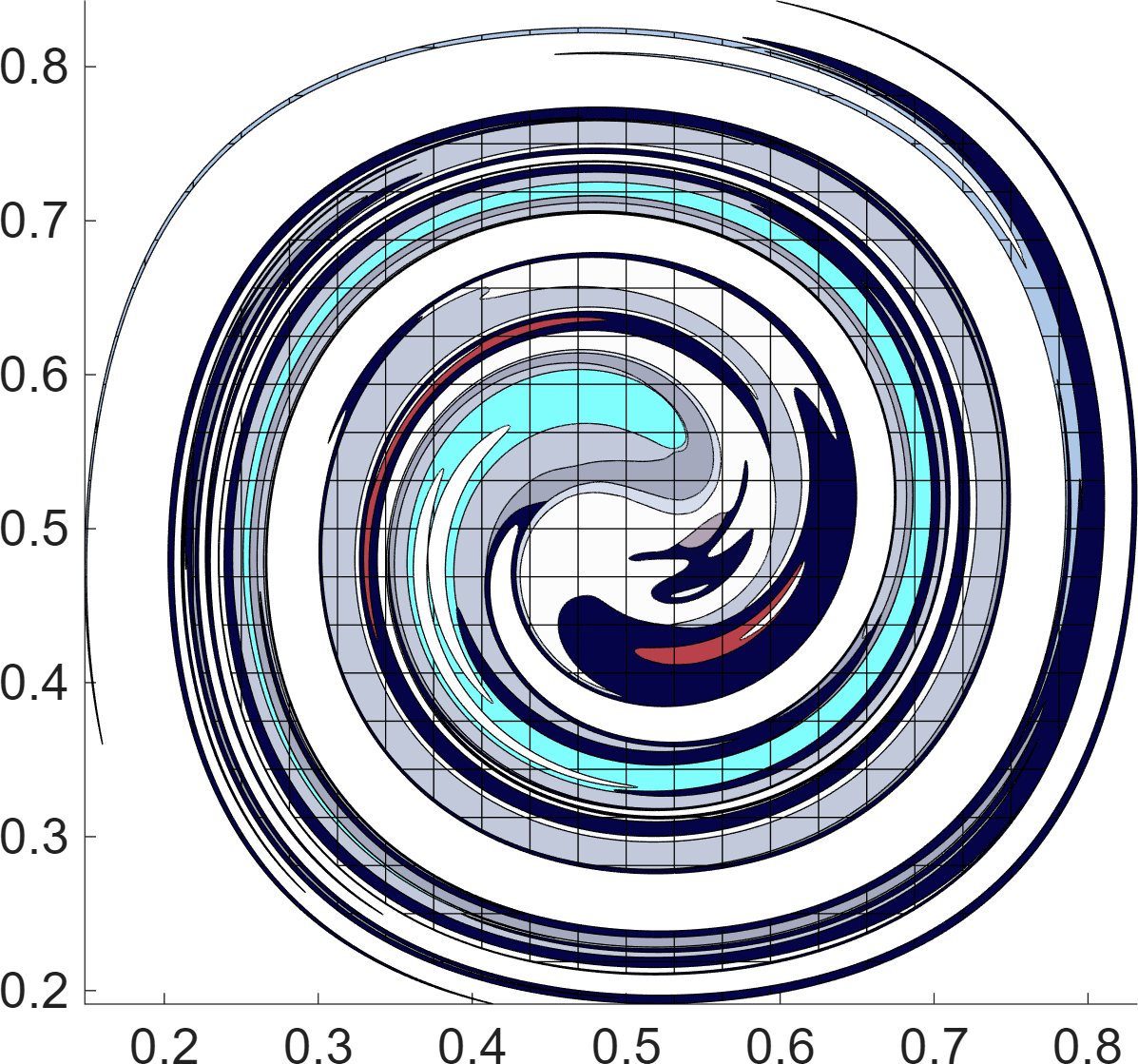}
  }
  \hfill
  \subfloat[$t = \frac{1}{4}T_d$; deformation]{
    \includegraphics[width=0.32\linewidth]{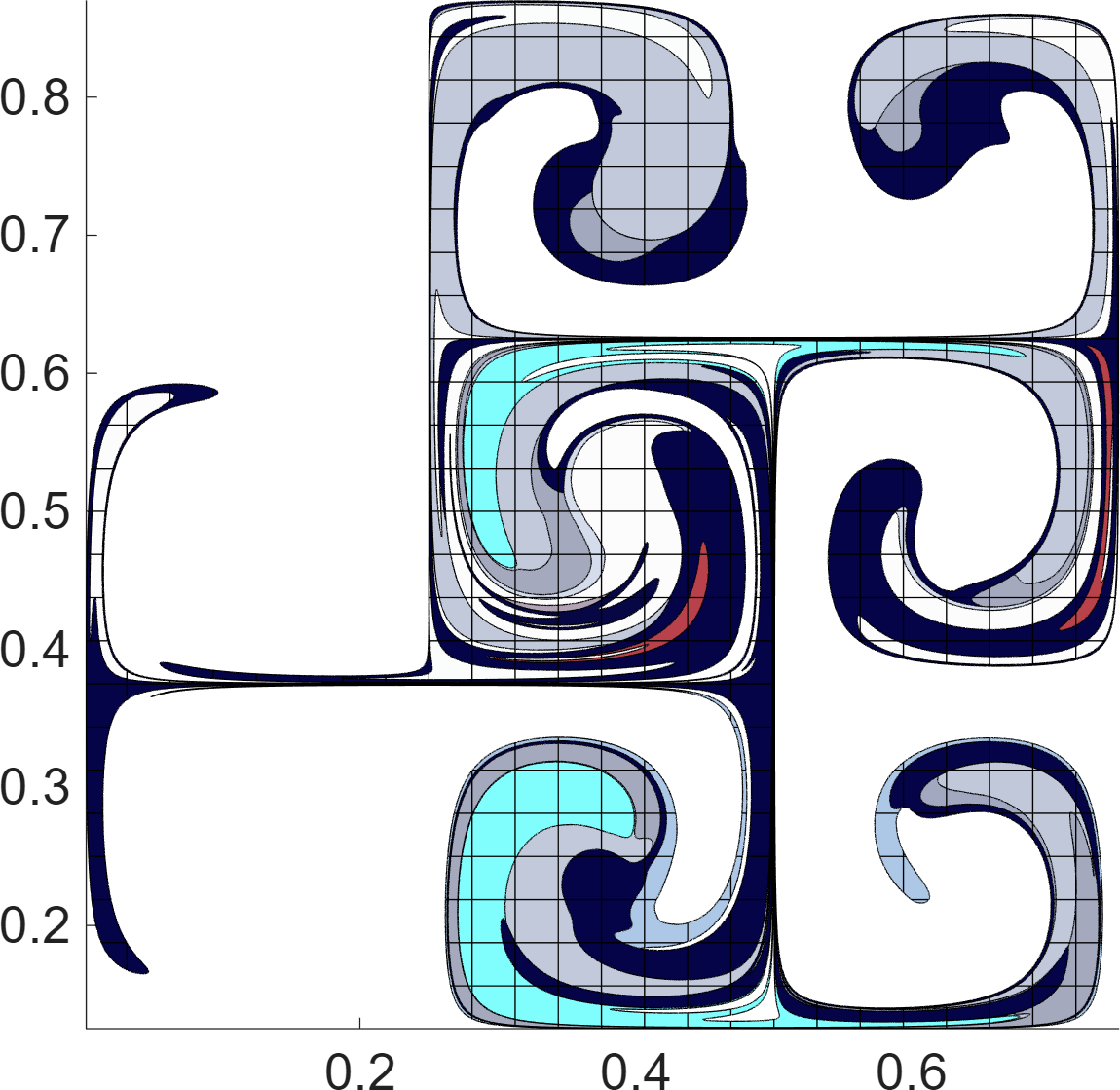}
  }

  \subfloat[$t = \frac{1}{2}T_v$; vortex]{
    \includegraphics[width=0.325\linewidth]{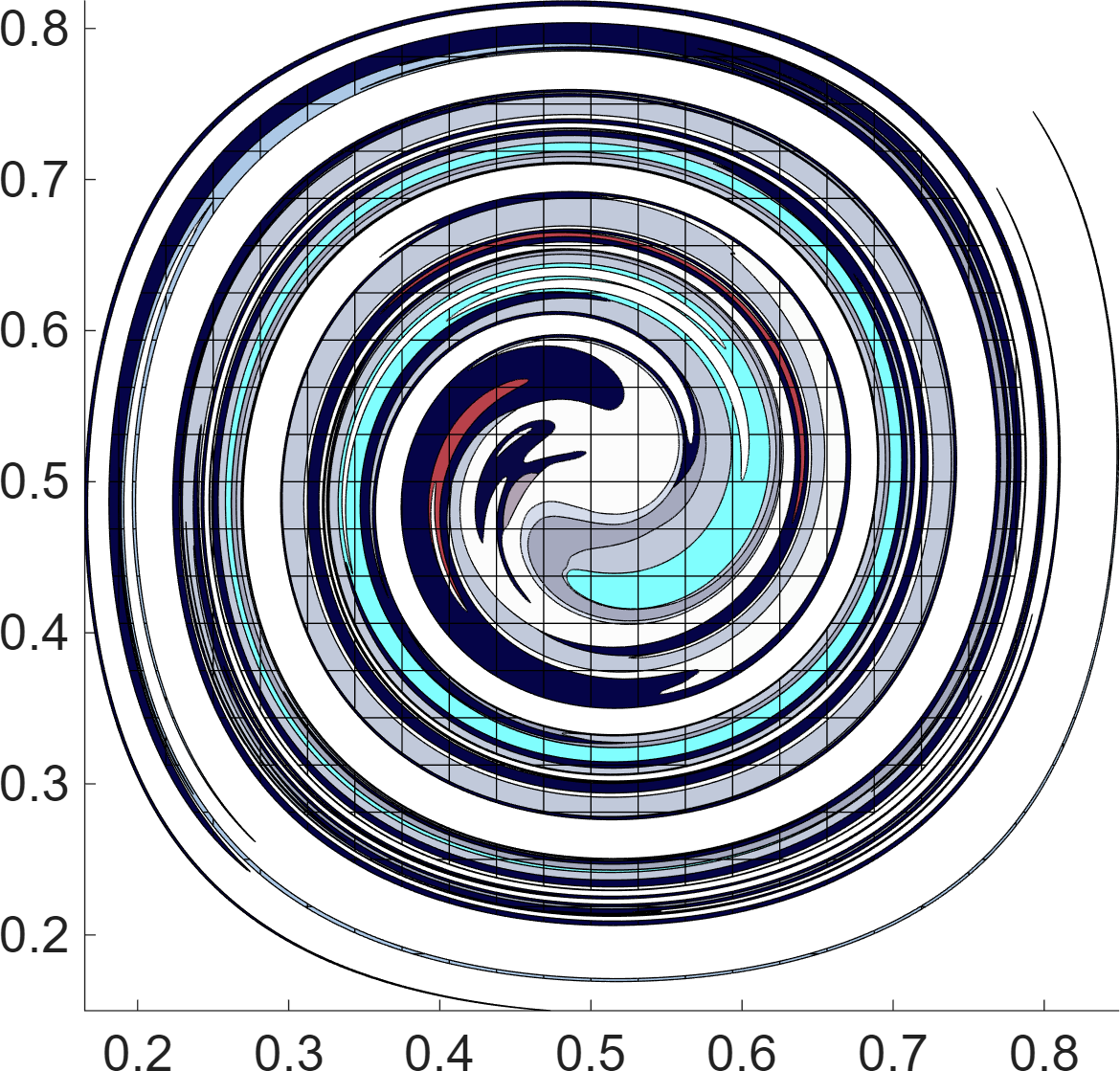}
  }
  \hfill
  \subfloat[$t = \frac{1}{2}T_d$; deformation]{
    \includegraphics[width=0.32\linewidth]{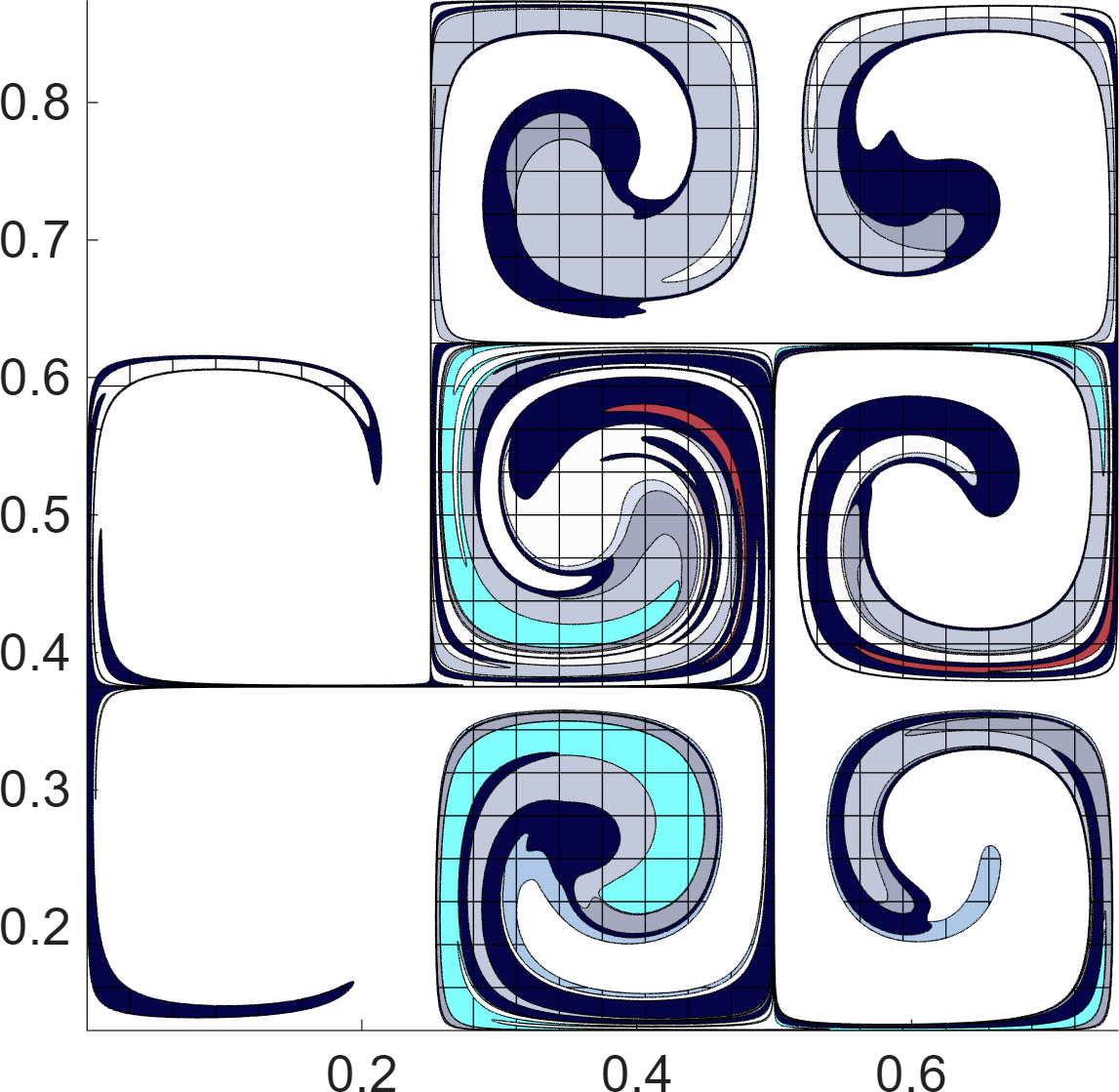}
  }
  \hfill
  \subfloat[the end of simulation]{
    \includegraphics[width=0.26\linewidth]{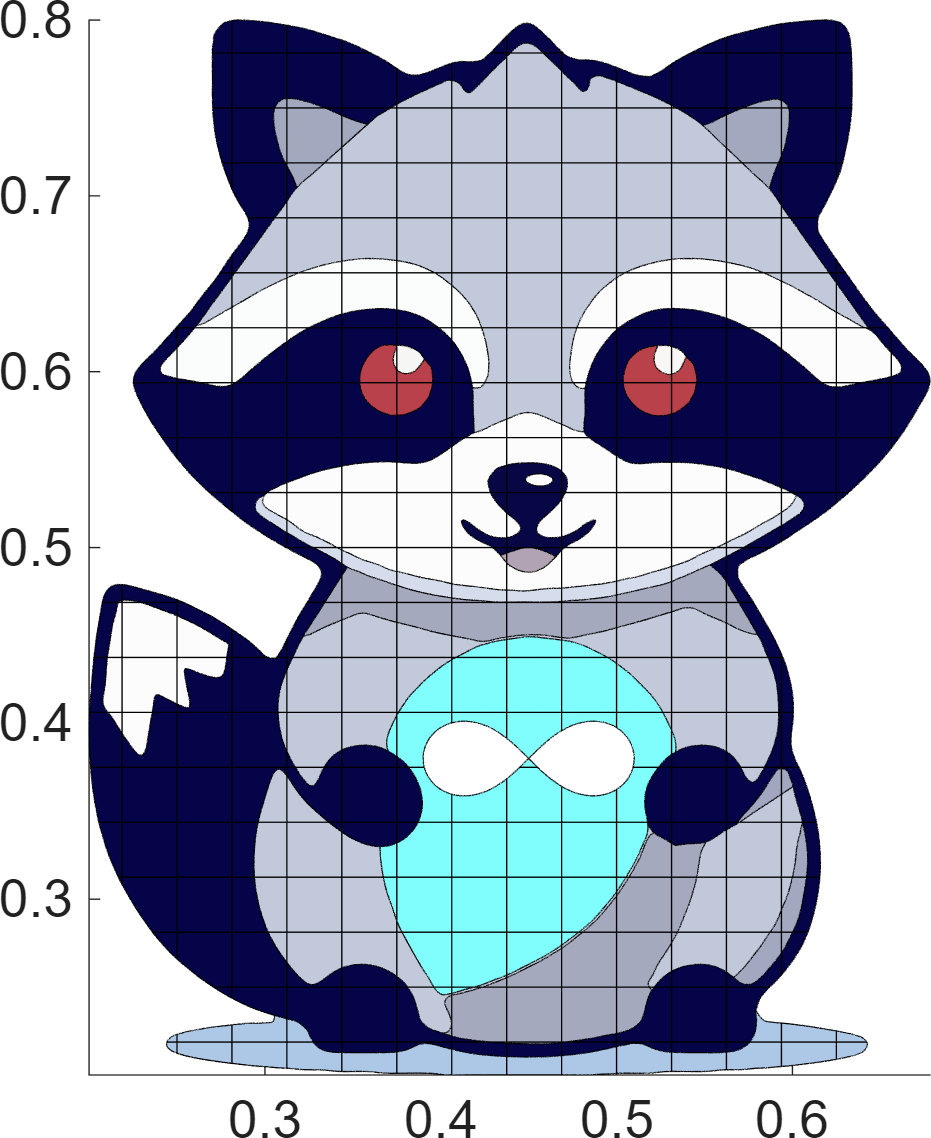}
  }
  
  \caption{Solutions of the multiphase cubic MARS method
    with $h_L^c = 0.2 h$
    for the vortex test ($T_v = 16$) and 
    the deformation test ($T_d = 2$)
    on the Eulerian grid of $h = \frac{1}{32}$.
    Subplots (a)--(f) are snapshots of the solution
    at key time instants. 
    The initial distances between markers are $\frac{1}{2}h_L(\rho)$,
    where $h_L$ is defined in (\ref{eq:curvature-based-arms}).
  }
  \label{fig:test-raccoon}
\end{figure}

The flow maps of this test are the same as those
 in (\ref{eq:vortex-velocity}) and (\ref{eq:deformation-velocity}).
As shown in Fig.~\ref{fig:test-raccoon}(a), 
 there are 22 tracked phases
 constituting a raccoon at the initial time.
On the chest of the raccoon,
 the ``$\infty$'' shape is represented by a lemniscate of Bernoulli
 $r^2 = a^2 \cos(2 \theta)$
 where $r$ denotes the radial distance from the origin, 
 $\theta\in \left[-\frac{\pi}{4}, \frac{\pi}{4}\right] 
 \cup \left[\frac{3\pi}{4}, \frac{5\pi}{4}\right]$
 the polar angle, 
 and $a$ the scaling parameter that
 determines the half-width of the curve.
Boundary curves of all phases at the initial time are represented
 by ${\cal C}^4$ quintic splines.

\begin{table}
  \caption{Errors and convergence rates
    of the multiphase cubic MARS method
    with $k = \frac{1}{8}h$ and curvature-based ARMS
    in solving the vortex test with $T_v = 12, 16$ and 
    the deformation test with $T_d = 2$.
    The initial condition is shown in Fig.~\ref{fig:test-raccoon}(a).}
  \centering
  \renewcommand{\arraystretch}{1.2} 
  \input{{\TABDIR}table6.7}
  \label{tab:raccoon-arms}
\end{table}

\begin{table}
  \caption{CPU time in seconds
    of the multiphase cubic MARS method
    in solving the raccoon test 
    in Sec.~\ref{sec:testRaccoon}
    on an Intel Core i7-12700 at 4.2GHz.
    The parameters are identical to those 
    in Table~\ref{tab:raccoon-arms}.}
  \centering
  \renewcommand{\arraystretch}{1.2} 
  \input{{\TABDIR}tableCpu2}
  \label{tab:raccoon-arms-Cputime}
\end{table}

Different from the rose curve in the previous section,
 the lemniscate of Bernoulli is smoothly connected
 for both pairs of arcs at its self-intersection
 and is therefore approximated by a periodic cubic spline.
The treatment of junctions and kinks in the previous subsection
 is adopted: 
 none of the non-smooth points are included in the interface graph
 while all junctions are.
 
Parameters of the curvature-based ARMS strategy 
 for the flow maps are the same as those 
 in (\ref{eq:curvatureStrategyGraphVortex}) and 
 (\ref{eq:curvatureStrategyGraphDeformation}),
 with values of $r_{\min}^c$ 
 listed in Table~\ref{tab:raccoon-arms},
 where the convergence rates of our MARS method
 are close to 4, 6, and 8 for 
 $h_L^c=O(h)$, $O(h^{\frac{3}{2}})$, and $O(h^2)$, respectively.
It is also clear in Fig.~\ref{fig:test-raccoon}
 the proposed MARS method preserves
 the topological structures at all time instances.
For both flow maps, 
 their final solutions at $t=T_v,T_d$ 
 are visually indistinguishable
 and are shown in Fig.~\ref{fig:test-raccoon}. 

Finally,  
 the computational cost of the proposed method
 with curvature-based ARMS
 for solving the raccoon test 
 are presented 
 in Table \ref{tab:raccoon-arms-Cputime}, 
 where the choices of $h_L^c$
 being $O(h)$, $O(h^{\frac{3}{2}})$, and $O(h^2)$
 result in growth rates of the CPU time
 being close to 2, 2.5, and 3. 
These timing results confirm Lemma \ref{lem:complexityOfARMS}: 
 the multiphase cubic MARS method
 has a linear complexity in each time step 
 with respect to the number of markers
 and thus its efficiency is optimal.
 


\vspace{5mm}
\section{Conclusion}
\label{sec:conclusion}
We have developed a cubic MARS method 
 with a curvature-based ARMS strategy
 for fourth- and higher-order IT of multiple materials.
The geometry of the interface
 is approximated by cubic splines
 while the topology 
 represented by directed cycles. 
For homeomorphic flow maps,
 the separation of the topology from the geometry
 leads to simple, efficient, and accurate algorithms: 
 topological structures are determined
 from the initial condition once and for all
 so that evolving the interface only entails
 advancing the periodic and not-a-knot cubic splines.
The superior generality, efficiency, and accuracy of our method
 are demonstrated by results of
 an array of benchmark tests.

A number of future research prospects follow.
First, we will prove the fourth- and higher-order accuracy
 of the proposed method under the MARS framework.
Second, 
 the Yin space, the MARS framework,
 the static representation of multiple Yin sets,
 and the general ARMS strategy
 pave the way to
 the study of topological changes of multiple phases; 
we will develop theoretical characterizations 
 and design highly accurate algorithms
 to tackle topological changes.
Third, we will couple this work
 with the PLG algorithm \cite{zhang2024PLG} and
 fourth-order projection methods \cite{zhang2016gepup}
 to form a generic fourth-order finite-volume solver
 for simulating incompressible multiphase flows on moving domains.
We also plan to apply this solver to study real-world applications
 such as wetting and spreading \cite{bonn09:_wettin}.



\begin{acknowledgements}
  \revise{The authors are grateful to one anonymous referee, 
    whose insightful comments lead to substantial improvements 
    of the original manuscript.
    We also acknowledge helpful suggestions
    from Shaozhen Cao, Shuang Hu, Junxiang Pan, Yunhao Qiu, 
    Jiatu Yan, and Chenhao Ye,} 
    graduate students in the School of Mathematical Sciences
    at Zhejiang University.
\end{acknowledgements}

\section*{Declarations}

\subsection*{Data Availability}

The datasets generated during and/or analysed during the current study
are available from the corresponding author on reasonable request.

\subsection*{Conflict of interest}

The authors declare that they have no conflict of interest.

\bibliographystyle{spmpsci}
\bibliography{bib/MARS-n2D}




\end{document}